\documentstyle[12pt,amsfonts,amscd,amssymb]{amsart}

\newtheorem{thm}{Theorem}[section]
\newtheorem{lm}[thm]{Lemma}

\newtheorem{prop}[thm]{Proposition}

\def\ul{\underline}
\def\la{\longrightarrow}
\def\ni{\noindent}

\def\a{\alpha}
\def\b{\beta}
\def\d{\delta }
\def\e{\epsilon}
\def\l{\lambda}

\def\I{\cal I}
\def\O{\cal O}

\def\X{\cal X}
\def\Y{\cal Y}

\def\C{\Bbb C}
\def\F{\Bbb F}

\def\P{\Bbb P}

\def\Z{\Bbb Z}

\def\mg-{\overline M_g}
\def\mult{\text{mult}}
\def\Pic{\text{Pic}}

\def\th{{\text{th}}}

 \def\ul{\underline}
 \def\la{\longrightarrow}
 \def\ni{\noindent}

\def\g{\gamma}

\def\Prod{(D_1\cdot D_2)}

\def\ND{N(D)}
\def\NDi{N_i(D)}
\def\NDj{N_{i_j}(D_j)}
\def\rDi{r_0^i(D)}

\def\VDm{V_{\underline m}(D)}
\def\VDi{V_i(D)}
\def\VDj{V_{i_j}(D_j)}
\def\rDu{r_0^{i_1}(D_1)}
\def\rDt{r_0^{i_2}(D_2)}
\def\rDT{r_0^{i_3}(D_3)}

\def\rD{r_0(D)}

\def\VD{V(D)}

\def\mg{\overline{M}_{0,4}}

\begin{document}

\

\centerline{\bf{PARAMETER SPACES FOR CURVES ON SURFACES}}
\centerline{\bf{AND ENUMERATION OF RATIONAL CURVES}}

\

\

\

\

\noindent {\bf Lucia Caporaso}

\noindent Mathematics department, Harvard University,

\noindent 1 Oxford st., Cambridge MA 02138, USA

\noindent caporaso{\char'100}abel.harvard.edu

\

\noindent {\bf Joe Harris}

\noindent Mathematics department, Harvard University,

\noindent 1 Oxford st., Cambridge MA 02138, USA

\noindent harris{\char'100}abel.harvard.edu

\

\

\tableofcontents

\section{Introduction}

In this paper we will be concerned with the geometry of families of rational
curves on a  surface $S$.

Let us state the main problem. Let $S$ be a nonsingular, rational surface,
and let
$D$ be an effective divisor class in $S$. We  denote by $|D|$
 the set of all effective divisors linearly equivalent to $D$;
 this is a projective space whose dimension we will denote by
$r(D)$. Inside $|D|$, we want to consider the locus of
rational curves: we let
$$
\tilde V(D) \; = \; \{ [X] \in |D| \  \mbox{such that $X$ is an irreducible,
rational curve}\} .
$$
This is a locally closed subset of the projective space $|D|$; we let $V(D)
\subset |D|$ be its closure. We call $V(D)$ the {\it Severi variety} of
rational curves
associated to
the divisor class $D$ on $S$, and we denote its dimension by
$r_0(D)$. We have in general
$$
r_0(D) \; \ge \; r(D) - p_a(D)
$$
with equality holding in all the cases that we shall study.

The particular aspect of the geometry of $V(D)$ of concern to us here is
 its degree, which we will denote by
$N(D)$.
This can also be characterized
directly: it is  the number of
irreducible rational curves that are linearly equivalent to
$D$ and that pass through $r_0(D)$ general points of $S$.
 The principal results of this
paper will be the computation of  $N(D)$ in
some cases.
For simplicity, we will assume that $N(D)$ is zero if
$V(D)$ is empty.

\subsection{The general strategy: the  cross-ratio method}

There are various approaches to the calculation of degrees  of Severi
varieties  (see for example \cite{CH}). The one we take in this paper is
based on ideas of
Kontsevich and Manin. In their paper \cite{KM} they describe a beautiful
formula, found
initially by Kontsevich, for the number  of plane rational curves of given
degree passing
through the appropriate number of points (another proof of it was given
independently by
Ruan and Tian in \cite{RT}). These methods have also been used to give formulas for the
degrees of genus 0 Severi varieties on other rational surfaces; see \cite{CM},
\cite{DI} and
\cite{KP}.

 The method that we are going to
describe   was suggested to us by the   ``First Reconstruction
Theorem" of \cite{KM}. It is based on the
analysis of a one-parameter family of curves:  we will consider
 the family of irreducible rational curves
 in the linear system
$|D|$ passing through
$r_0(D) - 1$ general points of $S$. By suitably marking four points on each
of these curves, we can (possibly after a base change) associate to the
family  a
cross-ratio function on its base.
Moreover, we can do this in such a way
that the $N(D)$ curves in the family passing through a given
$r_0(D)^\th$ point of $S$ are among the zeroes of the cross-ratio function,
and all the other zeroes and poles of the cross-ratio function
 will occur at reducible
fibes. Thus, to determine the number of curves in the family
passing through the last point, we need to describe the set of its
reducible elements  (and the multiplicity of the cross-ratio at
each).

This gives us, in principle, a way of solving the problem
recursively.  If $D$ and $D'$ are divisor classes on $S$,  we  say that
$D'<D$ if $D-D'$ is effective and non zero. Then
finding the number of reducible curves
$X$ in
$V(D)$ passing through $r_0(D) - 1$ general points involves
   curves (the components $X_i$ of
such curves $X$) in divisor classes $D' < D$,
which we may consider known inductively.
In fact, in
simple cases ($S = \P^2, \  \P^1 \times \P^1$ or $\F_1$), this works quite
smoothly:
every union $X = \cup X_i$ of rational curves $X_i \subset S$ is a limit of
irreducible rational curves, and we end up with an expression for the number
$N(D)$ of rational curves in the linear series $|D|$ through $r_0(D)$ points
simply in terms of the numbers
$N(D')$ for divisor classes $D' < D$. By contrast, a more delocate
situation arises when we consider other ruled surfaces $\F_n$: here the
components $X_i$ of a union
$X= \cup X_i$ have to satisfy additional conditions in order for the point $[X]
\in |D|$ to be in $V(D)$. In the following chapter, we will analyze exactly this
situation, and in the final chapter we will apply the results of this
analysis to
derive recursive formulas for $N(D)$ on $\F_n$.

Let us describe more precisely the set-up. Fix two
irreducible curves $C_3$, $C_4 \subset S$ having positive intersection
number with each other and with $D$, and intersecting transversely.
 Let
$q_1,\dots, q_{r_0(D)-1}
\in S$ be general points and denote by $H_{q_i}$  the hyperplane in $|D|$
parametrizing curves through $q_i$. Then let
$$
\Gamma \; = \; V(D) \cap H_{q_1} \cap \dots \cap H_{q_{r_0(D)-1}}
$$
be the corresponding linear section of $V(D)$---equivalently, the closure in
$|D|$ of the set of irreducible rational curves passing through $q_1,\dots,
q_{r_0(D)-1}$. Now, for a general point
$p$ in $S$, we can interpret
the degree of $V(D)$ as
 the number of points $[X] \in \Gamma$ corresponding to curves $X$  that
pass through
$p$. Let $\X \subset \Gamma \times S$ be the universal family over
 $\Gamma$, that is, the subscheme of $\Gamma \times S$ whose fiber over each
point $[X] \in \Gamma$ is simply $X$; let $f  :
\X \to \Gamma$ be the projection. By construction, $f : \X \to \Gamma$ is a
flat  family of  curves, whose general fiber is an irreducible
rational curve.

Next, we introduce a family whose general fiber is the normalization of the
corresponding fiber of $\X \to \Gamma$. To do this, we first take
$\Gamma^\nu \to \Gamma$ the normalization of $\Gamma$, pull the family
$\X \to \Gamma$ back to $\Gamma^\nu$, and take $\X^\nu$ to be the
normalization of the total space of this pullback: that is, we set
$$
\X^\nu \; = \; (\X \times_\Gamma \Gamma^\nu)^\nu .
$$ The composite map $\X^\nu \to \Gamma^\nu$ (which we will again
denote by $f$) is then a flat family,
 with general fiber
isomorphic to $\P^1$.

Now, we want endow this family with four sections
$$p_i \; : \; \Gamma^\nu \; \la \; \X^\nu , \  \  \  i=1,2,3,4\  \text{ such that
}\  f\circ p_i = id_{\Gamma^\nu},
$$
so as to define a
cross-ratio function.  To do that,  we will take the first two sections to be
given by the first two base points $q_1$ and $q_2$, that is, for general
$[X] \in \Gamma$ we set
$$
p_i([X]) \; = \; ([X], q_i) \quad \text{for} \quad i = 1, 2 ;
$$
this then extends to a regular map on all of $\Gamma^\nu$.  (The extension
follows from the fact that $\Gamma^\nu$ is a smooth curve. The reader may
wonder about  the points $[X]$ corresponding to the
curves in our family with nodes at
$q_1$ and
$q_2$, where it may at first appear the section $p_i$ can't be well-defined;
but in fact
$\Gamma$ will have ordinary nodes
and those points
and the two branches of the normalization $\Gamma^\nu$ over those points
exactly correspond to the choices of value for $p_i$.)

For the remaining two sections, we want to pick out a point of intersection of
each curve $X$ of our family with each of the curves $C_3$ and $C_4$, that
is, we want to choose for general $[X] \in \Gamma$,
$$
p_j([X]) \; = \; ([X], r_j) \quad \text{with} \quad r_j \in X \cap C_j \quad
\text{for}
\quad j = 3,4 .
$$
This requires another base change. To be precise, let $\pi_2 :
\Gamma \times S \to S$ be the projection, and for
$j = 3,4$ set
$$
{\cal C}_j \; = \; \X \cap \pi_2^{-1}(C_j) \; \subset \; \Gamma \times S .
$$
The projection ${\cal C}_j \to \Gamma$ will have degree $d_j = (D \cdot
C_j)$.  We let ${\cal C}_j^0$ be
the  union of those components of ${\cal C}_j$ of
positive degree over $\Gamma$, and finally let $B$ be the normalization of
their product over $\Gamma$, that is,
$$
B \; = \;  ({\cal C}_3^0 \times_\Gamma {\cal C}_4^0)^\nu .
$$
Then if we let
$$
\X ' \; = \; \X^\nu \times_{\Gamma^\nu }B\; \longrightarrow \; B
$$
be the pullback family, we  define two further sections $p_3, p_4 : B \to
\X "$  by
$$
p_j(b) \; = \; (b, \pi_j(b))
$$
for general $b \in B$,
where $\pi_j : B \to {\cal C}_j$ is the projection.

We now have a family $\X ' \to B$ of curves over a smooth base $B$, with
general fiber $X'$  isomorphic to $\P^1$. We may accordingly define the
cross-ratio map
$$
\phi \; : \; B \; \longrightarrow \; \P^1
$$
by letting, for general $b \in B$, $\phi(b)$ be the cross-ratio of the four
points
$p_1(b), p_3(b), p_4(b), p_2(b) \in X' \cong \P^1$---that is, in terms of any
affine coordinate on $X'\cong \P^1$ we set
$$
\phi(b) \; = \; {(p_1 (b)-p_2(b)) (p_3 (b)-p_4(b))\over  (p_1 (b)-p_3(b))(p_2
(b)-p_4(b))}.
$$
Equivalently,   the family $\X ' \to B$ will determine a family
of stable $4$-pointed rational curves over  $B$, and hence a canonical
morphism  $\tilde\phi$ from $B$ to  $\mg$. (It may not be apparent that a
further base change is not necessary to arrive at such a family; but this
is true
and will emerge in the subsequent analysis.) We can then choose an
identification
$\eta$ of
$\mg$ with $\P ^1$ so that  the three points  of $\mg$
corresponding to singular (that is, reducible) curves

\

\vspace*{1.4in}

\ni \hskip-.3in \special {picture crossratio}

\

\ni map to the points $0, \infty $
and
$1$. The
composition $\eta \circ \tilde\phi$ is then the map $\phi$ above.

Now we want to compute the number of zeroes and poles of $\phi$ to obtain a
formula for
$N(D)$.
Observe now that the cross-ratio
function
$\phi$ can have a zero or a pole at a point $b \in B$ only under one of two
circumstances: if two of the points $p_i(b)$ actually coincide; or if the fiber
$X_b$ of $\X ' \to B$ over $b$ is reducible. As for the first of these, clearly
$p_1(b)$ can never equal $p_2(b)$; and since the points $q_i$ are general,
they cannot lie on $C_3$ or $C_4$. The first
case will thus occur exactly when $p_3(b) = p_4(b)$, which in turn can happen
only when the curve $X_b$ contains one of the points $p$ of intersection of
$C_3$ with
$C_4$. Conversely, for every curve in the original family $\X \to \Gamma$
containing a point $p \in C_3 \cap C_4$, there will be a unique point $b \in B$
lying over $[X] \in \Gamma$ with $p_3(b) = p_4(b)$, namely the point $(p,p)
\in {\cal C}_3^0 \times_\Gamma {\cal C}_4^0$ (the fact that this is a
smooth point of ${\cal C}_3^0 \times_\Gamma {\cal C}_4^0$, so
that there will be a unique point of $B$ over it, follows from the fact that
every curve in the original family $\X \to \Gamma$
containing a point $p \in C_3 \cap C_4$ will intersect $C_3$ and $C_4$
transversely, which will be a consequence  of
Proposition
\ref{dimensioncount} below). Thus,  every curve of our original family
passing through a point of $C_3 \cap C_4$  will contribute a zero to the
function $\phi$. Once we verify that these are
all simple zeroes, we conclude that $\phi$  has a total of $(C_3 \cdot
C_4)N(D)$  zeroes at points $b \in B$  corresponding to irreducible curves
$X_b$.

It remains to describe the zeroes and poles of $\phi$ coming from reducible
curves $X$ in the family. Here there is a fundamental difference between the
``simple" cases of surfaces $S = \P^2$, $\P^1 \times \P^1$,  or $\F_1$ and the
general surfaces $\F_n$. The difference is this: for surfaces such as
$\P^2$, the
dimension $r_0(D)$ of the family of rational curves in a given linear series
$|D|$ is affine linear in the class $D$, without exception: it is given by
the formula
$$
r_0(D) \; = \; -(K_S \cdot D) - 1 .
$$
It follows that for any two effective divisor classes $D_1$, $D_2$ with $D_1
+ D_2 = D$ we have
$$
r_0(D_1) + r_0(D_2) \; = \; r_0(D) - 1,
$$
for three divisor classes $D_1$, $D_2$, $D_3$ with $D_1 +D_2 +D_3 = D$ we
have
$$
r_0(D_1) + r_0(D_2) + r_0(D_3) \; = \; r_0(D) - 2,
$$ and so on.  As a consequence, we see that  there are no curves $X \in
|D|$ with three or more  components, all of which are rational,  passing through
$r_0(D)-1$ general points of
$S$;  and there are exactly
$$
\sum_{D_1+D_2=D}{r_0(D)-1 \choose r_0(D_1)}N(D_1)N(D_2)
$$
 such curves with two components: for every division of the $r_0(D)-1 =
r_0(D_1) + r_0(D_2)$ points $q_i$ into subsets of $r_0(D_1)$ and
$r_0(D_2)$, there will be exactly
$N(D_1)N(D_2)$ pairs $X_1 + X_2$ with $X_1$ containing the first set and
$X_2$ containing the second. By a naive dimension count, moreover, all such
curves will be limits of irreducible
rational curves. In these cases, then, we can inductively enumerate all
reducible elements of our family, and say which contribute to the zeroes and
which to the poles of
$\phi$; in the end, equating the degrees of the divisors of zeroes and poles of
$\phi$ we are led to a recursive formula for the number $N(D)$.

For the general $\F _n$ instead, a dimension count will not be enough to
describe all possible types of degenerations
occuring in codimension 1. This is due to the presence of the exceptional
curve $E$ on $\F _n$.
We shall see that
in a family of generically  irreducible nodal curves there will be special
fibers containing $E$ as a component.
The description of these types of degeneration will in fact require a
delicate analysis.

\subsection{Two sample calculations}  To give an idea  of what sort of
information
we need  to
 carry out our project, we will  derive formulas for two relatively simple
examples. The first of these is the formula found by Kontsevich
for the degrees of the Severi varieties of rational plane curves. The
second is a direct
calculation of the degree of the Severi variety associated to the linear series
$|2C|$ on $\F_2$. We will not justify   the assertions made about the
types of degenerate fibers occurring in our families, or about the
multiplicities
of the corresponding zeroes and poles of the cross-ratio function, but we will
indicate when such justifications  are necessary,
and refer to the
 results of the following chapters.

\subsubsection{Kontsevich's formula} The cross-ratio method is well
 illustrated in the case of $S=
\P^2$. We will derive here Kontsevich's formula for the number $N(d)$ of
rational
curves passing through $3d-1$ general points in the plane.
(Note that
the dimension
 of the variety of rational curves in the linear series $|\O_{\P^2}(d)|$ is
$3d-1$. )

To proceed, we  take
$C_3$ and $C_4$ distinct lines, and choose $q_1,\dots,q_{3d-2}$
general points of the plane. We let $\Gamma \subset |\O_{\P^2}(d)|$ be the
closure of the locus of irreducible rational curves of degree $d$ passing
through the points $q_i$, and carry out the remaining steps of the construction
described above to arrive at a family $\X ' \to B$ whose base $B$ is a finite
cover of the curve $\Gamma$ of degree $d^2$, with four sections $p_i : B \to
\X '$ coming from the two base points $q_1$ and $q_2$ and the points of
intersection of the curves in the family with the lines $C_3$ and $C_4$.

As already explained, our formulas will be obtained by
 equating the degrees of the
divisors of zeroes and poles of the cross-ratio function $\phi$ on $B$. To begin
with, the only way $\phi$ can have a zero or pole at a point $b \in B$
corresponding to an irreducible curve $X_b$ is if the curve $X_b$
 passes through the point $p = C_3 \cap C_4$ of intersection of the
two lines.

Now, since the $3d-1$ points $q_1,\ldots,q_{3d-2},p$ are general, there
will be exactly $N(d)$ curves $X$ of our original family $\X \to \Gamma$ passing
through $p$. Moreover, these curves will all correspond to smooth points
$[X] \in
\Gamma$, so that there will again be exactly $N(d)$ points of $\Gamma^\nu$
corresponding to such curves. Next, since each such curve
$X$ intersects
$C_3$ and $C_4$ transversely at $p$, the point $(p,p)
\in {\cal C}_3^0 \times_\Gamma {\cal C}_4^0$ will be a smooth  point of
$\Gamma^\nu$ and there will be a unique point of $B$ lying over it.
Finally, again as
a consequence of the analysis in Chapter 2, each such point of $B$ will be
a simple
zero of the cross-ratio function; so that in sum  we have a total of exactly
$N(d)$ zeroes of
$\phi$ at points $b \in B$ with $X_b$ irreducible.

It remains to count the zeroes and poles of $\phi$ occurring at points
corresponding to reducible fibers. We need
to make one remark here in general, about the geometry of $\Gamma$ near a
point $[X]$ corresponding to a reducible curve $X = X_1 \cup X_2$. If the degree
of the components $X_i$ is $d_i$, then each $X_i$ will have ${d_i-1 \choose 2}$
nodes, and in addition there will be $d_1d_2$ transverse points of
intersection of
$X_1$ with $X_2$; thus $X$ will have
$ {d-1 \choose 2} + 1$
nodes in all. Now, as we approach $[X]$ along any branch of $\Gamma$, we see
that the limiting positions of the ${d-1 \choose  2}$ nodes of the general
curve of
our family will consist of the ${d_1-1 \choose 2}$ nodes of $X_1$, the ${d_2-1
\choose 2}$ nodes of
$X_2$, and all but one point $r$ of the $d_1d_2$ points of
intersection of the two components (so that when we take the normalization
of the
total space along that branch, what we see is the normalizations of the $X_i$,
joined along the one point of each lying over the point $r$). Conversely, if we
choose any point $r \in X_1 \cap X_2$, there exist deformations of $X$ smoothing
the node $r$, preserving the remaining nodes, and of course continuing to pass
through the points $q_1,\dots,q_{3d-2}$. Thus, in sum, we see that the curve
$\Gamma$ will have exactly $d_1d_2$ branches at the point $[X]$.

With that said, let us proceed to count the remaining zeroes and poles of
$\phi$.
We have to ask: how many reducible curves $X$ are there in our family
with $q_1$ and $q_2$ on one component, and $p_3$ and $p_4$ on the other (by
the  dimension count made above in general, $X$ can have only two components);
and how many are there with $q_1$ and $p_3$ on one component, and $q_2$ and
$p_4$ on the other? To start, we can certainly count the number of points
$[X] \in
\Gamma$ corresponding to reducible curves $X = X_1 \cup X_2$ with $q_1, q_2
\in X_1$: such a curve will consist of a component $X_1$ of degree
$d_1$ passing through $q_1, q_2$ and $3d_1-3$ of the other points $q_i$, and a
curve $X_2$ passing through the remaining $3d-4-(3d_1-3)=3d_2-1$ of the points
$q_3,\dots,q_{3d-2}$. To specify such a curve, then, we have first to break
up the
set
$\{q_3,\dots,q_{3d-2}\}$ into subsets of $3d_1-3$ and $3d_2-1$ points; then
choose $X_1$ any one of the $N(d_1)$ rational curves of degree $d_1$ through
$q_1$,
$q_2$ and the first set, and $X_2$ any one of the $N(d_2)$ rational curves of
degree $d_1$ through the second set. The total number of such points on
$\Gamma$ is thus
$$
N(d_1)N(d_2){3d-4 \choose 3d_1-3} .
$$
Next, by the remark above, for each such point $[X] \in \Gamma$, there will be
$d_1d_2$ points of $\Gamma^\nu$ lying over it; thus we have
$$
N(d_1)N(d_2){3d-4 \choose 3d_1-3}d_1d_2 .
$$
such points of $\Gamma^\nu$. Moreover, of the $d^2$ points of $B$ lying
over each
such point of $\Gamma^\nu$, exactly $d_2^2$ will correspond to triples
$([X],p_3,p_4)$ with $p_3, p_4 \in X_2$; so we have a total of
$$
N(d_1)N(d_2){3d-4 \choose 3d_1-3}d_1d_2^3 .
$$
 points of $B$ corresponding to unions $X = X_1 \cup X_2$ of curves of degrees
$d_1$ and $d_2$.

Finally, we have to check that each of the corresponding points of $B$ is a
simple
zero of the cross-ratio function $\phi$, which will  follow from
Lemma~\ref{crossratiomult}.
With this verified, we see in sum
 we have a total of exactly
$$
\sum_{d_1+d_2=d} N(d_1)N(d_2){3d-4 \choose 3d_1-3}d_1d_2^3 .
$$
 zeroes of
$\phi$ at points $b \in B$ with $X_b$ reducible. This now accounts for all the
zeroes of $\phi$, so that we have
$$
\deg \phi \; = \; N(d) + \sum_{d_1+d_2=d} N(d_1)N(d_2){3d-4 \choose
3d_1-3}d_1d_2^3 .
$$

The poles of $\phi$ are  counted in the same fashion as the zeroes coming
from reducible fibers. In fact, the only difference is that now we have to
look at
curves $X = X_1 \cup X_2$ consisting of a component $X_1$ of degree
$d_1$ passing through $q_1$ and $3d_1-2$ of the other points $q_i$, and a
curve $X_2$ passing through $q_2$ and the remaining $3d-4-(3d_1-2)=3d_2-2$ of
the points
$q_3,\dots,q_{3d-2}$. To specify such a curve, then, we have first to break
up the
set
$\{q_3,\dots,q_{3d-2}\}$ into subsets of $3d_1-2$ and $3d_2-2$ points; so that
the total number of such points on
$\Gamma$ is now
$$
N(d_1)N(d_2){3d-4 \choose 3d_1-2} .
$$
Again, for each such point $[X] \in \Gamma$, there will be
$d_1d_2$ points of $\Gamma^\nu$ lying over it; and now, of the $d^2$ points of
$B$ lying over each such point of $\Gamma^\nu$, exactly $d_1d_2$ will
correspond to triples
$([X],p_3,p_4)$ with $p_3 \in X_1$ and $p_4 \in X_2$; so we have a total of
$$
N(d_1)N(d_2){3d-4 \choose 3d_1-2}d_1^2d_2^2 .
$$
 points of $B$ corresponding to unions $X = X_1 \cup X_2$ of curves of degrees
$d_1$ and $d_2$.

As before, we have to check that each of the corresponding points of $B$ is a
simple pole of $\phi$; once we have done this, since these are all the
poles of $\phi$
we conclude that
$$
\deg \phi \; = \;  \sum_{d_1+d_2=d} N(d_1)N(d_2){3d-4 \choose
3d_1-2}d_1^2d_2^2 .
$$
Equating the two values of the degree of $\phi$, we arrive at Kontsevich's
formula:
$$  N(d) =
\sum _{d_1 + d_2 =d} N(d_1)N(d_2)d_1d_2\left[
\left( \begin{array}{c} 3d-4 \\ 3d_1-2
\end{array} \right) d_1d_2 -
\left( \begin{array}{c} 3d-4 \\ 3d_1-3
\end{array} \right) d_2^2 \right] .
$$

\

\subsubsection{The linear series $|2C|$ on $\F_2$}
\label{2C2}  Let $C\in \Pic (\F _2)$ be as defined in
   \ref{term}.
The linear series $|2C|$ on $\F_2$
has dimension 8 and arithmetic genus $p_a(2C) = 1$; the Severi variety $V(2C)
\subset |2C|$ correspondingly has dimension 7, with general member $X$ a curve
with one node. To carry out the calculation, then, we select two general curves
$C_3$, $C_4$, which we  choose to be linearly equivalent to $C$,
and 6 general points $q_1,\dots,q_6 \in S$. We let $\Gamma
\subset V(2C)$ be the locus of curves in $V(2C)$ passing through
$q_1,\dots,q_6$ and
proceed from there as before to construct the family $\X ' \to B$, where
$B$ will be in
this case a $(2C \cdot C_3)(2C \cdot C_4) = 4\cdot 4 = 16$-sheeted cover of the
normalization $\Gamma^\nu$ of $\Gamma$.

Again, the elements of the family $\X \to \Gamma$ passing through the two
points of
intersection of $C_3$ and $C_4$ each contribute a simple zero to $\phi$,
and again
these are the only singularities of $\phi$ at points $b \in B$ with $X_b$
irreducible.

As
for the reducible fibers, we can readily describe those not containing the curve $E$.
Since  $(2C \cdot E) = 0$ any curve $X$ in the linear system not
containing $E$ will be disjoint from it; and since the only curves on
$\F_2$ disjoint
from $E$ are those linearly equivalent to a multiple of $C$, any reducible
element of
$|2C|$ not containing $E$ must be of the form $X = X_1 + X_2$ with $X_1
\sim X_2 \sim
C$. Now, the  curves of this type form a six-dimensional subvariety of
$|2C|$, and it is
not hard to see that they are all in $V(2C)$, that is, they are all limits
of irreducible
singular curves in $|2C|$. In fact, a general such curve $X = X_1 + X_2$
has two nodes,
coming from the intersection $X_1 \cap X_2$, each of which may be the limit
of the
node of a nearby irreducible rational curve. The base $\Gamma$ of our
original family
will thus have two branches at the point $[X]$ corresponding to such a
curve, and the
normalization $\Gamma^\nu$ two points lying over $[X]$.

We count the number of zeroes and poles of $\phi$ at points  $b \in B$
lying over such
points as before. A zero arises when the two base points $q_1, q_2$ lie on
the same
component of $X$, which we will designate $X_1$; thus, to specify such a
fiber of $\X
\to \Gamma$ we have to pick a curve $X_1 \in |C|$ passing through $q_1$,
$q_2$ and
one point $q_j$ of the four remaining points $q_3,\dots q_6$; and a curve
$X_2 \in |C|$
passing through the remaining three of the points $q_3,\dots q_6$. $X$ will thus
be determined simply by the choice of $j \in \{3,4,5,6\}$. There will be
two points of
$\Gamma^\nu$ lying over each such point $[X] \in \Gamma$, as remarked; and
a total
of $(C \cdot C_3)(C \cdot C_4) = 2 \cdot 2 = 4$ points of $B$ lying over each,
corresponding to the choice of $p_3 \in X_2 \cap C_3$ and $p_4 \in X_2 \cap
C_4$.
Finally, each such point of $B$ will be a simple zero of $\phi$, so we have
a total of
$$
{4 \choose 1}\cdot 2 \cdot 4 \; = \; 32
$$
zeroes of $\phi$ of this type.

The number of poles is described analogously: a pole arises when we have
$q_1$ and
$p_3$ on one component of $X = X_1 + X_2$ and $q_2$ and $p_4$ on the other;
so to
specify such a curve $X$ we  choose a subset $\{q_i,q_j\} \subset
\{q_3,\dots,q_6\}$
and take $X_1$ the (unique) curve in $|C|$ through $q_1$, $q_i$ and $q_j$,
and $X_2$
the curve in $|C|$ through $q_2$ and the remaining two points of
$\{q_3,\dots,q_6\}$.
Again, there are two points of $\Gamma^\nu$ over each such $[X] \in \Gamma$, and
four points of $B$ over each of those at which $\phi$ has a  pole. These
poles are
simple, and so we have a total of
$$
{4 \choose 2}\cdot 2 \cdot 4 \; = \; 48
$$
poles of $\phi$ of this type.

It remains to describe the curves $X$ in our family $\X \to \Gamma$ that
contain $E$.
The curves in the linear
series $|2C|$ that contain $E$ form a hyperplane $\Sigma_E = E + |C+2F|
\cong \P^7
\subset |2C|
\cong \P^8$. Thus, even though a general such curve $X = E + X_1$ looks
exactly like
one of the curves discussed in the last case---two smooth rational
components meeting
transversely at two points---for purely dimension-theoretic reasons it
cannot be a
limit of irreducible rational curves in the series $|2C|$. The key question
is, then:
which curves $X = E + X_1$ are limits of irreducible rational curves? Or,
equivalently, what is the intersection
$V(2C)
\cap
H _E$?

To answer this,  let us  see why a general such curve $X = E + X_1$
cannot be such a limit.  Suppose we had a one-parameter family $\X
\subset \Gamma \times \F_2 \to
\Gamma$ of irreducible rational curves specializing to such a curve $X$. After
normalizing the base and total space of such a family we would arrive at a
family
$\X^\nu \to \Gamma^\nu$ whose general fiber was smooth, and whose special fiber
had one component $\tilde E$ mapping to $E$ and one component $Y$ mapping to
$X_1$ via the projection $\pi : \X^\nu \to \X \to \F_2$. Since the
arithmetic genus of
the special fiber must be 0, $\tilde E$ and $Y$ will meet at one point
(this will be the
point lying over the node of $X$ that is not a limit of nearby curves in
the family).

\

Now,
the inverse image
$\pi^{-1}(E) \subset \X^\nu$ does not meet the general fiber of $\X^\nu \to
\Gamma^\nu$; and since it must have pure codimension 1 it can only consist of
$\tilde E$ itself:
$$
\pi^{-1}(E) \; = \; \tilde E .
$$
This is impossible:  since $Y$ maps onto $X_1$, it has two points
mapping to $E$; but $\pi^{-1}(E) \cap Y = \tilde E \cap Y$ will consist of
just one point.
Thus $X$ cannot be a limit of irreducible rational curves.

This analysis also suggests which reducible curves containing $E$ are such
limits.
Basically, the contradiction above derived from two hypotheses: that $X_1$
met $E$
transversely in two points; and that $Y$ met $\tilde E$ in just one point.
In fact, we
will see in the following chapter, as a very special case of Propositions
\ref{codimensionone} and \ref{singtotalspace}, that conversely if either of
those fails, the curve
$X$ is such a limit. In other words, the intersection $V(2C) \cap H_E$
consists of the union of
two loci:

\

\ni $\bullet$ \thinspace The locus of curves $X = E + X_1$ with $X_1$
tangent to $E$;
and

\

\ni $\bullet$ \thinspace The locus of curves $X = E + X_1 + F$, with $X_1
\sim C+F$ and
$F$ a fiber of $\F_2$.

\

In the latter case, the apparent contradiction above is resolved because
$Y$ will have
two components, each meeting $\tilde E$ at a point, so that the picture of
$\X^\nu \to
\Gamma^\nu$ is

\

\vspace*{2.2in}

\ni \special {picture F2type2}

\

\ni In other words, the limiting position of the nodes of nearby fibers of a
one-parameter family of irreducible rational curves specializing to $X$
must be the
node $X_1 \cap F$ of $X$. In particular, at a point $[X]$ corresponding to
such a curve in
our original family $\X \to \Gamma$, $\Gamma$ will have a single branch,
and there
will correspondingly be a unique point of $\Gamma^\nu$ lying over $[X] \in
\Gamma$.

We can at last count the remaining zeroes and poles of $\phi$. First of
all, the curves
of the form $X = E + X_1$ with $X_1$ tangent to $E$ do not contribute at
all: for such a
curve, all four points $p_i$ must be distinct points of the component
$X_1$. Secondly,
since a fiber $F$ can pass through at most one of the two points $q_1,
q_2$, the curves
of the form
$X = E + X_1 + F$, with
$X_1
\sim C+F$ and
$F$ a fiber of $\F_2$ will contribute a zero only if the curve $X_1$
contains $q_1$,
$q_2$ and three of the four points $q_3,\dots ,q_6$, $F$ is the fiber of $\F_2$
through the remaining point of $q_3,\dots ,q_6$, and $p_3$ and $p_4$ are
chosen to be
the unique points of intersection of $F$ with $C_3$ and $C_4$ respectively.
There are
thus exactly four such points of $\Gamma$, and over each there will be a
unique point
of $\Gamma^\nu$, and a unique point of $B$ at which $\phi$ has a zero.
There is one
new wrinkle here: as we will see in Propositions \ref{singtotalspace} and
\ref{crossratiomult} below,
each of these points is a double zero of $\phi$. Thus we have a total of 8
zeroes of $\phi$ at
such points.

The poles of $\phi$ at such points are counted similarly. Here there will
be exactly two
curves $X$ contributing: we could take $F$ the fiber through $q_1$, and
$X_1$ the
curve in $|C+F|$ through $q_2,\dots,q_6$, and then take $p_3 = F \cap C_3$
and choose
$p_4 \in X_1 \cap C_4$; or we could reverse the roles of 1 and 2, and 3 and
4. In each
case, there will be a unique point of $\Gamma^\nu$ lying over $[X] \in
\Gamma$ and,
since $(X_1 \cdot C_3) = (X_1 \cdot C_4) = 3$, three points of $B$ lying
over it. Finally,
these points are similarly double poles of $\phi$, so we have a total of 12
poles of
$\phi$ at such points.

The calculation is now complete: equating the number of zeroes and poles of
$\phi$, we
find that
$$
2 \cdot N(2C) + 36 + 8 \; = \; 48 + 12
$$
and hence
$$
N(2C) \; = \; 10 .
$$

Actually, there are many (easier) ways to calculate this number. But this
example serves to illustrates the questions we must answer in order to apply the
cross-ratio method to calculate degrees of rational Severi varieties in
general, and also
the sort of answer we may find. In the following
chapter, we will present a series of results describing exactly what sort of
reducible fibers we may expect to find, in codimension $1$,
 for general linear series $|D|$ on $\F_n$; and the
geometry of the varieties $\Gamma$ and $\X$  in a neighborhood of the
corresponding points and curves. In the last chapter, we will apply the results
of this analysis to derive recursive formulas for the degrees of  some
Severi varieties on these surfaces.

\subsection{Notation and Terminology}
\label{term}

Our base field will be the field of complex numbers.

Let
$\F_n =\P (\O _{\P ^1}\oplus \O _{\P ^1}(n))$  be a rational ruled surface.
The Picard group of
$\F_n$ has rank $2$, and we choose generators as follows:
$$
\Pic (\F_n )=  \Z\cdot C \oplus  \Z \cdot F
$$
where $C ^2 = n$, $\  F^2=0$
and $F\cdot C =1$. We denote by $E$ the unique curve of negative self
intersection, so that $E^2=-n$ and
$E \sim C-nF$.

\

Let
$D$ be any divisor class on the surface $S = \F_n$ other than $E$, and let
$\ul{m} := (m_1,m_2,\dots,m_k)$ be any sequence of positive integers with
$\sum m_i =
(D \cdot E)$. We define the locally closed subvariety $\tilde V_{\ul{m}}(D)
\subset V(D)$
be the locus of irreducible rational curves $X$ such that, if $\nu :
\P^1 \to X$ is the normalization of $X$, the
pullback  divisor
$$
\nu ^*(E) = \sum m_i\cdot q_i
$$
for some collection of distinct points $q_1,\ldots,q_k \in \P^1$,
and we let $V_{\ul{m}}(D) \subset V(D)$ be its closure; for example, as we
will see, if
$\ul{m} = (1,1,\ldots, 1)$, then $V_{\ul{m}}(D) =V(D)$.  When $\ul{m}$
contains  a single integer $i$ greater than $1$
(i.e. $\ul{m} = (i,1,1,\ldots ,1)$), we will denote these by $\tilde
V_{i}(D)
$ and
$\VDi $ respectively.
We set
$$\rDi  = \dim (\VDi)$$
and
$$\NDi  = \deg \VDi$$
We have $\VD$ for $V_1(D)$,   $\ND = N_1(D)$ and $r_0(D) = r_0^1(D)$.
We define $\NDi$ to be zero if $\VDi$ is empty.

Similarly, let $\Omega = \{p_1,\ldots,p_k\} \subset E \subset \F_n$ be any
collection of
$k$ distinct points. We let $\tilde W^\Omega_{\ul{m}}(D) \subset V(D)$  be
the locus of
irreducible rational curves $X$ such that, if $\nu :
\P^1 \to X$ is the normalization of $X$, them for some collection of
distinct points
$q_1,\ldots,q_k \in \P^1$ we have
$$
\nu(q_i) = p_i
$$
and
$$
\nu ^*(E) = \sum m_i\cdot q_i
$$
and again let $W^\Omega_{\ul{m}}(D) \subset V(D)$ be its closure.

\subsection{Summary of  results}
Here we give a list of some known formulas
  including all the ones that we prove in this paper.
These
recursions  are very similar from a  formal point of view.
In what follows, we will state them in a way that
highlights the analogies.

To begin with, fix any  rational surface $S$, a divisor class
$D$ on it, and two curves $C_3$ and $C_4$. For any pair of divisor classes
$D_1$ and $D_2$ we introduce the function

$$\gamma (D_1,D_2):=$$

$$N(D_1) N(D_2)
\left[
\left( \begin{array}{c}
r_0(D) - 3 \\ r_0(D_1)-1
\end{array} \right) (D_1\cdot C_3)(D_2\cdot C_4) -
\left( \begin{array}{c}
r_0(D) -3 \\ r_0(D_1)-2
\end{array} \right)  (D_2\cdot C_3)(D_2\cdot C_4)
 \right]
$$

Using this notation, we state the following results

 \

\ni  {\bf Recursion for $\P ^2$} \thinspace (\cite{KM}) {\it Let}  $C_3$
{\it and }
$C_4$ {\it be two fixed lines in the plane, then}
$$N(D) =
\sum _{D_1 + D_2 =D} \gamma (D_1,D_2) (D_1 \cdot D_2).$$

\

\ni  {\bf Recursion for $\P ^1 \times \P ^1$} \thinspace  (\cite{KM},
\cite{DI}, \cite{KP}) {\it Let
} $C_3$ {\it and } $C_4$ {\it be two fixed elements of the two
 distinct rulings, then}
$$N(D) =
\sum _{D_1 + D_2 =D} \gamma (D_1,D_2) (D_1 \cdot D_2).$$

\

The first new result of this paper is a similar recursion formula for the
degrees of
 Severi
varieties of rational curves on the ruled surface $\F_2$. Note the slightly
different form of the
recursion: the presence of an extra term not analogous to those in the two
preceding formulas is
due to the contribution of degenerate curves containing $E$.

\

\ni  {\bf Recursion for $\F_2$}  \thinspace  (Theorem \ref{F2}) {\it Let }
$C_3$ {\it
and } $C_4$ {\it be two fixed elements of the  class $C$, then}
$$2N(D) =
\sum _{D_1 + D_2 =D} \gamma (D_1,D_2) (D_1 \cdot D_2)\;  + \; 2\sum _{D_1 +
D_2 =D-E} \gamma (D_1,D_2) (D_1 \cdot E)(D_2\cdot E).$$

\

Now, we will see that in $\F_n$, the
general reducible curves
$X =
\cup X_i \in |D|$ that are limits of irreducible rational curves and
contain $E$
 have the property that each component
$X_j$ may have a point of tangency of order
$i_j$ with $E$---that is, will belong to $V_{i_j}(D_j)$, where $D_j$ is the
divisor class of $X_j$.
Accordingly, we shall define later  (Section~\ref{Fn}) a generalized
version of the number
$\gamma (D_1,D_2)$;  this will be a function
$\gamma _{i_1, i_2,\dots,i_t}(D_1,D_2,\dots,D_t)$
depending  recursively on the  degrees $\NDj$.
In these terms, we give a formula expressing the degree $N(D)$ of
$V(D)$ on $\F_n$ in terms of the degrees of the tangential Severi varieties
of smaller divisor
classes:

\ni  {\bf A sample formula for $\F_n$}  \thinspace (Theorem \ref{Fn})

$$
n\ND =
\sum _{D_1 + D_2 = D} \Prod \g _{1,1}(D_1,D_2) +
      $$

$$
+ \sum _{t=2}^{n} \; \sum _{D_1+D_2+\dots +D_t=D-E} \; \sum
_{i_1,\dots,i_t}  \; \prod_{j
:i_j=1}(E\cdot D_j)
\g_{i_1,\dots,i_t}(D_1,\dots, D_t)
$$

\

The difference here is that in case $n \ge 3$ this does not give a complete
recursion: to be able to
enumerate rational  curves on such surfaces, we would need formulas for the
degrees of the
``tangential" Severi  varieties as well, that is, we need formulas for
$\NDi$. The first case for which this occurs is that of $\F _3$. Very possibly
a complete recursion  could still be obtained using the
cross-ratio method,  although  the level of difficulty seems to us to get
very high.
Instead we found a different technique that we successfully applied in a
few cases;
for example, we obtained a complete set of recursions for the surface $\F _3$.
This  different method is the subject of another paper of ours (cf.
\cite{CH}); it also is  heavily
based on the deformation theory results that are developed in the second
chapter of this paper.

Finally, we obtain a simple closed formula for the class $2C$ on any ruled
surface
$\F_n:$

\

\ni {\bf Closed formula for 2C on $\F_n$} \thinspace (Theorem \ref{2C})

$$N(2C)=\sum _{k=0}^{n-1}(n-k)^2 {2n+2 \choose k}$$

\newpage

\section{Degenerations of rational curves}

In this chapter we  prove the  results on
degenerations of rational curves that we will need to obtain
our formulas.

\subsection{The basic set-up }\label{setup}

We start with the complete linear system
$|D|$ associated to a divisor class
$D$ on the ruled surface $S = \F_n$, and with the Severi variety $V(D)
\subset |D|$. We then choose
$ r_0(D) - 1$ general points $q_1,\ldots,q_{r_o (D) -1} \in S$, and let
$\Gamma$ be the intersection of
$V(D)$ with the linear subspace of curves in
$|D|$ passing through $q_1,\ldots, q_{r_o (D) -1}$; we let $\X \subset \Gamma
\times S$ be the corresponding family of curves over $\Gamma$.

Next, we let $\Gamma^\nu \to \Gamma$ be the normalization of the base $\Gamma$,
and

$$
\X^\nu \; = \; \left( \X \times_\Gamma \Gamma^\nu \right)^\nu \; \longrightarrow
\; \Gamma^\nu
$$
the normalization of the pullback of the family to $\Gamma^\nu$, so that
$\X^\nu \to \Gamma^\nu$ is a family whose general fiber is a smooth rational
curve.  If $X$ is a fiber of $\X \to \Gamma$, the notation
  $X^\nu$ will  be used for a corresponding fiber of the family
$\X^\nu
\to
\Gamma^\nu$,  which may differ from
the normalization of  $X$.

Then we fix two curves $C_3$ and $C_4$ in $\F _n$, which will be linearly
equivalent to $C$.
We
need to make a further base change
$B
\to
\Gamma^\nu$, so that the points of intersection of the curves
 in our family with

$C_3$ and $C_4$ become rational over the base. We  thus let $B \to
\Gamma^\nu$ be any finite cover, unramified at the points $b \in \Gamma^\nu$
with $X^\nu_b$ singular, and let $\X' \to B$ be the pullback of the family
$\X^\nu \to
\Gamma^\nu$ to $B$. (By Propositions
\ref{dimensioncount} and \ref{codimensionone}, the map $B \to
\Gamma^\nu$ introduced in Chapters 1 and 3 in order to define the sections
$p_i$
 will indeed be
unramified at the points of $B$ corresponding to the singular fibers of
$\X^\nu
\to B$.) Because the results of this chapter are all local in the base of
our family,
however, we will not need to introduce this extra step in the construction.  For
the remainder of this chapter, accordingly, we will take $B = \Gamma^\nu$;
and all of
the results of the chapter describing the map $\X^\nu \to \Gamma^\nu$ will
still hold after the
base change $B \to
\Gamma^\nu$.

Next we  introduce the {\it nodal reduction}
of the family
$\X' \to B$. That is to say, after making a base change $\tilde B \to B$ and
blowing up the pullback family $\X  '\times _B \tilde B \to \tilde B$, we
arrive at a
family $\Y \to \tilde B$ such that

\

\begin{enumerate}
\item $\Y \to \tilde B$ is a family  all of whose fibers of $\Y \to
\tilde B$ are reduced curves having only nodes as singularities;
\item the total space $\Y$ is smooth;
\item $\Y$ admits a regular birational map $\Y
\to \X'
\times_B \tilde B$ over $\tilde B$.
\end{enumerate}

In fact, most of our concerns with this definition will turn out in the end
to be unnecessary: we will see below as a corollary of Propositions
\ref{describegamma} and \ref{singtotalspace} that in fact
$\X' \to B$ is already a family of nodal curves.
Thus, in practice, we will not have to make a base change at all at this
stage, and $\Y$ will be simply the minimal desingularization of
$\X'$. For this reason (and because $B$ is itself already an arbitrary finite
cover of the normalization $\Gamma^\nu$ of our original base $\Gamma$) we will
abuse notation slightly and omit the tilde in
$\tilde B$, that is, we will speak of the family
$\Y
\to B$.

One further remark: in the applications we will have four sections of the
family $\Y \to B$ and will correspondingly want to consider this as a
family of four-pointed nodal curves. For this reason, we may want to
make further blow-ups at points where these sections cross.
By  Propositions
\ref{dimensioncount} and \ref{codimensionone}, however, the sections in question
will cross only at smooth fibers of $\Y \to B$ and so this will not affect our
descriptions of the singular fibers of the family.

The final construction is one that we will use only in the following
chapter, but we  mention it here just to have all the definitions in one
place. After arriving as above at a family $\Y
\to B$ of nodal curves with four disjoint sections $p_i$, we  may
then proceed to blow down ``extraneous" components of fibers $Y$ of
$\Y \to B$:  that is, any component of
$Y$ that meets the other components of $Y$ in only one point, and that
meets at most one of the sections $p_i$. Iterating this process until
there are no extraneous components left, we arrive at what we will call
the {\em minimal smooth semistable model} of our family: that is, a
family
${\cal Z}
\to B$ such that ${\cal Z}$ is smooth, the fibers  are nodal, the
sections
$p_i$ are disjoint and ${\cal Z} \to B$ is minimal with respect
to these properties. Note that the special fiber
$Z$ of ${\cal Z}$
must be a chain of rational curves $G_0,\dots,G_\ell$ with two of
the sections meeting each of the two end components:

\vspace*{2.8in}

\hskip.6in \special {picture semistable}

\

\noindent (the case $\ell = 0$ is simply the case where $Z$ is
irreducible). Finally, we can blow down the intermediate components
$G_1,\dots,G_{\ell - 1}$ in this chain to arrive at a family ${\cal W}
\to B$ of 4-pointed stable curves, called the {\it stable model} of our
family. The
special fiber
 of this family will have just two components (or one, if $\ell = 0$), with a
singularity of type $A_\ell$ at the point of their intersection.

In sum, we have the diagram of families and maps:

\vspace*{4in}

\hskip.5in \special {picture stablediagram}

\

\subsection{The main results from deformation theory}

We  give here a summary of the main results to be proved in this chapter.

\

\ni $\bullet$ \thinspace The first is Proposition~\ref{dimensioncount} in which
we consider the  Severi varieties  $V(D)$ and $\VDm$, compute their dimension and
 describe the geometry of their general point.
In particular, we characterize the general fiber of the
family
$\X  \to \Gamma$.   The results are  unsurprising: for example, the
general point $[X]$ of  $V(D)$ corresponds to a curve $X$
with only nodes as singularities; general points $[X]$,
$[X']$ of,  respectively, $V(D)$ and $V(D')$ correspond to curves $X$,
$X'$ that intersect transversely..

\

\ni $\bullet$ \thinspace  Then, in
 Proposition
\ref{codimensionone},
 we study the geometry of the general point of the boundary of
$V(D)$.  We  do that by listing all types of reducible fibers that occur in the
family $\X \to\Gamma$.
This result is not  predictable on the basis of a simple dimension count;
as we have seen in example \ref{2C2}, in most linear systems
$|D|$ on $\F_n$ the subvariety corresponding to  reducible rational
curves containing $E$   is larger-dimensional than $V(D)$; so the question
of which points of the
former lie in the closure of the latter does not have an immediate  answer.

\

\ni $\bullet$ \thinspace The third result is Proposition ~\ref{describegamma},
which is specifically about the family $\X \to \Gamma$. We describe  the
geometry of the base $\Gamma$ in a neighborhood of each point $[X]
\in
\Gamma$ corresponding to a degenerate fiber $X$.  In particular, we
say how many branches $\Gamma$ has at $[X]$ and say how the
nodes of the nearby irreducible fibers approach the singularities of
$X$ as we approach
$[X]$ along each branch of $\Gamma$.

\

\ni $\bullet$ \thinspace Finally we have Proposition ~\ref{singtotalspace},
describing the singularities of the total space of the
families $\X \to \Gamma$ and
$\X^\nu \to \Gamma^\nu$.  This will be a crucial ingredient in calculating the
multiplicities of zeroes of the cross-ratio function on the base of our
family.

\

One word of warning is in order. Many of both the statements and
proofs of these propositions are just routine verifications of statements
easily guessed on
the basis of naive dimension counts.  At the same time, mixed in with
these largely predictable statements are some  interesting
phenomena . These are described in the second parts of Propositions
\ref{codimensionone}, \ref{describegamma} and \ref{singtotalspace}, in
which we describe the geometry of the one-parameter families  $X \to
\Gamma$ and $\X^\nu \to \Gamma^\nu$ in a neighborhood of the reducible fibers
containing $E$. Near such a curve, the local geometry of the universal
family over the Severi variety is,  to us, somewhat surprising.

 \subsection{The geometry of the Severi varieties}
Here is the first result about the varieties $\VDm$ defined in
section~\ref{term}.

\begin{prop}\label{dimensioncount} Let $|D|$ and $|D'| \ne |E|$ be any linear
series on the surface
$S = \F_n$; let $G \subset S$ be any fixed curve not containing $E$ and let
$P_1,P_2,\ldots \in
S$ be any given finite collection of  points.
Let $\ul{m} = (m_1,m_2,\dots)$ be any collection of positive
integers
with $\sum m_i = (D\cdot E)$.

\ni 1. If $V_{\ul{m}}(D) $ is nonempty, then it has pure dimension
$$
\dim(V_{\ul{m}}(D)) = -(K_S\cdot D)-1- \sum (m_i - 1) .
$$

\ni 2. A general point $[X]$ of any component of
$V_{\ul{m}}(D)$ corresponds to a curve $X
\subset S$ having only nodes as singularities, smooth everywhere along
$E$, intersecting $G$ transversely and not containing $P_i$ for any $i$.

\ni 3. If $[X]$ and $[X']$ are general points of irreducible components of
$V_{\ul{m}}(D)$ and
$V_{\ul{m}'}(D')$ respectively, then $X$ and $X'$ intersect transversely,
and none of their points of intersection lie on $G$ or $E$.
\end{prop}

\

\ni \underbar{Remark}. Many of the techniques necessary to prove this statement
are in
\cite{H}. In fact, many of these assertions are proved there, but unfortunately
with slightly different hypotheses: they are proved first on a general
rational surface
$S$, but only for $V(D)$, that is, without the tangency condition
(Proposition (2.1) of [H]); and then with a single tangency condition, but
only with respect to a line in the plane (Lemma (2.4) of [H]).

\

\begin{pf}  We start with the dimension statement. The assertion that
the dimension of
$V(D)$ is everywhere equal to $-(K_S\cdot D)-1$ is standard deformation
theory (and is well known; c.f. \cite{K}).
To see it, observe first that if
$[X]
\in
\tilde V(D)$ is any point and $\nu : X^\nu \to X \subset S$ the normalization of
the corresponding curve, the first-order deformations of the map $\nu$ are
given by sections of the pullback
$\nu^*(T_S)$ of the tangent bundle to $S$. Now, the tangent bundle to
the ruled surface $S=\F_n$ is generated by its global sections everywhere except
along
$E$; since $X$ doesn't contain $E$, it will likewise be true that the pullback
$\nu^*(T_S)$ will be generically generated by its global sections. Since $X^\nu
\cong \P^1$, it follows in turn that $h^1(X^\nu,
\nu^*(T_S)) = 0$. The deformations of the map $\nu$ are thus unobstructed,
from which it follows that the space of such deformations is smooth of
dimension
\begin{equation*}
\begin{split} h^0(X^\nu, \nu^*(T_S)) &= \deg(\nu^*(T_S)) + 2 \\ &= -(K_S
\cdot D) + 2 .
\end{split}
\end{equation*}
If we mod out by automorphisms of the domain $\P^1$,
 we see that the space of deformations of the image curve $X
\subset S$ as a rational curve has dimension
\begin{equation*}
\begin{split} h^0(\P^1, \nu ^*(T_S)) - 3 = -(K_S \cdot D) - 1 .
\end{split}
\end{equation*}
which is the same as the dimension of $T_{[X]}V(D)$.

We next establish the
\proclaim Claim.  The dimension of $\tilde V_{\ul{m}}(D)$, and
hence of
$V_{\ul{m}}(D)$, is everywhere at least
$r_0(D)-\sum (m_i-1)$.

To see this, set $l = (D \cdot E)$. Let $[X] \in
\tilde V(D)$ be any point, $U$ an analytic neighborhood of $[X]$ in $\tilde
V(D)$, $\X \subset U \times S \to U$ the universal family of
curves over
$U$, and $\X^\nu$ and $U^\nu$ the normalizations of $\X$ and $U$; we may
assume that
the map $\tau :
\X^\nu
\to U^\nu$ is smooth. Now let $\X^\nu_l$ the $l^\th$ symmetric fiber
product of $\X^\nu
\to U^\nu$. We then have a map
\begin{equation*}
\begin{split}
\rho \; : \; U \; &\la \; \X^\nu_l \cr [X] \; &\longmapsto \;
\psi_{[X]}^* \nu_{[X]}^*(E)
\end{split}
\end{equation*}
Now, inside the symmetric product
$\X^\nu_l$, the locus $\Gamma_{\ul{m}}$
 of divisors having  points of multiplicities $m_i$ or
more is irreducible of codimension
$\sum (m_i-1)$; since
$\tilde V_{\ul{m}}(D) \cap U$ is an open subset of the inverse image
$\rho^{-1}(\Gamma_{\ul{m}})$, it follows that it must have dimension at
least
$\dim(V(D)) -
\sum (m_i-1)$ everywhere.

\

Note that an analytic neighborhood $U$ of any point of $\tilde V_{\ul{m}}(D)$
admits a map to $E^k$, sending $[X] \in U$ to the images $q_i = \nu(p_i)$;
the fibers of
this map are analytic open sets in the varieties $W^\Omega_{\ul{m}}(D)$. In
particular, we have
$$
\dim(V_{\ul{m}}(D)) \; \le \; \dim(W^\Omega_{\ul{m}}(D)) + k
$$
so that in order to prove the opposite inequality $\dim(V_{\ul{m}}(D)) \le
r_0(D)-\sum (m_i-1)$, it is enough to show that the dimension of the variety
$W^\Omega_{\ul{m}}(D)$ is equal to
$r_0(D)-\sum m_i$ for any subset $\Omega =
\{p_1,\ldots,p_k\} \subset E$.

\

To prove the remaining parts of the Proposition requires a tangent space
argument. This comes in two parts: first, we will identify the projective
tangent space
to the
space of deformations of a given reduced curve $X$ preserving the geometric
genus of $X$; and
then the subspaces  corresponding to deformations that also preserve
singularities other
than nodes and/or tangencies with fixed curves. This is the part that is in
common with [H], and
for the most part we will simply recall here the statements of the relevant
results
(Theorem~\ref{Zariski} and Lemma~\ref{tangencycondition}). Then, to apply
these, we need to
estimate the dimension of these subspaces of
$|D|$; this is carried out in Lemma~\ref{indcons} and the following argument.

We may identify the tangent space to the linear series
$|D|$ at $[X]$ with the {\it characteristic series}
$$
H^0(X, \O_X(X)) = H^0(S,
\O_S(X))/\C\tau
$$
 where $\tau \in H^0(S,
\O_S(X))$ is the section vanishing along $X$ (this identification is natural
up to scalars; more precisely, the tangent space to $\P(H^0(S,
\O_S(X)))$ at $[X] = \C\tau$ is
$$
{\text{Hom}}(\C\tau, H^0(S,
\O_S(X))/\C\tau = (\C\tau)^* \otimes H^0(S,
\O_S(X))/\C\tau ).
$$

Now suppose that we are given any subvariety $W$ of the linear series
$|D|$ on $S$. Let
$[X] \in W$ be a general point of $W$. The following theorem of Zariski
(\cite{Z}, Theorems 1 and 2) characterizes the tangent space to $W$ at $[X]$:

\begin{thm}\label{Zariski} {\bf (Zariski's theorem)} In terms of the
identification of the
tangent space to the linear series
$|\O_S(D)|$ at $[X]$ with the characteristic series $H^0(X, \O_X(X))$,

1. The tangent space $T_{[X]}W$ is contained in the subspace $H^0(X,
{\cal I}(X))$ of $H^0(X,
\O_X(X))$, where
$\I \subset \O_S$ is the  {\em adjoint ideal} of $X$;

2. If $X$ has any singularities other than nodes, then   $T_{[X]}W$
is contained in a subspace $H^0(X, {\cal J}(X))$ where ${\cal J}
\subsetneqq \I$ is an ideal strictly contained in the adjoint ideal.
\end{thm}

This characterizes the tangent space to $V(D)$ at a general point $[X]$.
(If the fact that it
does is not clear, it will be after Lemma \ref{indcons} below.) Now, we have to
consider the additional information coming from the tangency with $E$.
To express this, note first that, if $\nu : X^\nu \to X$ is the
normalization of $X$ and ${\cal J} \subset \O_{X}$ is any ideal
contained in the adjoint ideal of $X$, then the pullback map gives a
natural bijection between ideals ${\cal J}
\subset {\cal I} \subset \O_X$ contained in ${\cal I}$
and ideals
$ \nu^*{\cal J}\subset \nu^*{\cal I} \subset \O_{X^\nu}$. We
will invoke this correspondence implicitly in our notation: if $p \in
X^\nu$ is any point, and ${\cal J} \subset \O_{X}$ any ideal
contained in the adjoint ideal of $X$, we will write ${\cal J}(-mp)
\subset
\O_{X}$ to mean the ideal in $\O_X$ whose pullback to $X^\nu$ is
$\nu^*{\cal J} \otimes \O_{X^\nu}(-mp)$. In these terms, we have
the following
\begin{lm}\label{tangencycondition}
Let
$G
\subset S$ be any fixed curve and $p \in G$ a smooth point of $G$. Let  $W$
be any
subvariety of
$|D|$. If the general point
$[X]$   of $W$ satisfies the condition:  there is a point $q \in X^\nu$
 such that $\nu(q) = p$ and
$$
\mult _q (\nu^*(G)) \; = \; m \, ,
$$
then the tangent space to $W$ at $[X]$ satisfies
$$
T_{[X]}W \; \subset \; H^0(X, {\cal I}(X)(-mp)) .
$$
Moreover, if $X$ has any singularities other than nodes, or is singular at
the point $p$, we
have
$$
T_{[X]}W \; \subset \; H^0(X, {\cal J}(X)(-mp))
$$
 where ${\cal J}
\subsetneqq \I$ is  an ideal strictly contained in the
adjoint ideal.
\end{lm}
\begin{pf}
We will  prove the Lemma by applying Zariski's
theorem to the proper transform of
$X$ on the surface
$\tilde S$ obtained by blowing up  $S = \F_n$ a total of $m$ times along
the curve $E$. To
carry this out, let
$S_1 \to S_0$ be the blow-up of $S_0 = S$ at the point $p$, $E_1 \subset
S_1$ the
exceptional divisor of the blow-up and $p_1 \in E_1$ the point of
intersection of  $E_1$ with
the proper transform of $E$ in $S_1$. Similarly, let $S_2 \to S_1$ be the
blow-up of $S_1$
at the point $p_1$, $E_2 \subset S_2$ the exceptional divisor of the
blow-up and $p_2 \in
E_2$ the point of intersection of  $E_2$ with the proper transform of $E$
in $S_1$, and so on,
until we arrive at the surface $\tilde S = S_m$; we will denote by $\pi :
\tilde S
\to S$ the composite of the blow-up maps, by $\tilde X$ the proper
transform of $X$ in $\tilde S$ and by $\tilde E_i$ the proper transform of
$E_i$ in $\tilde S$;
so that the pullback to $\tilde S$ of the divisor $E$ is given by
$$
\pi^*E \; = \; \tilde E + \sum i \cdot \tilde E_i \, .
$$
We  denote by $X'$ the branch of $X$ corresponding to the
point
$q
\in X^\nu$, that is, the image of an analytic neighborhood of
$q$ in $X^\nu$, by $\tilde X'$ its proper transform in $\tilde S$, and by
$\tilde p$ the
point of $\tilde X'$ lying over $p$.

Now, let $X_i$ be the proper transform of $X$ in $S_i$, and let $k_i$ be
the multiplicity of
$X_{i-1}$ at the point $p_{i-1}$; for each $j = 1,\ldots,m$ we will set
$$
l_j \; = \; k_1+k_2+\ldots+k_j \, .
$$
Thus, for example, we have the equality of divisors
$$
\pi^*X \; = \;  \tilde X +  \sum_{i=1}^m l_i \cdot \tilde E_i \, .
$$
Similarly, we let $X_i'$ be the proper transform on $X'$ in $S_i$, $k'_i$
the multplicity of
$X'_{i-1}$ at $p_{i-1}$ and $l'_j = k'_1+\ldots k'_j$. Note that $l_j \ge
l'_j$ for each $j$;
and the requirement that
$X'$ have intersection multiplicity
$m$ with
$E$ at
$p$ is equivalent to the assertion that
$$
mult_p(X' \cdot E) \; = \; (\pi^*X' \cdot \tilde E) \; = \; l'_m \; = \; m \, ,
$$
so that we have in particular $l_m \ge m$, with equality if and only if
(locally) $X = X'$.
We can also write the intersection number $m_p(X' \cdot E)$  as
$$
mult_p(X' \cdot E) \; = \; mult_q(\tilde X' \cdot \pi^*E) \; = \;
m_q(\tilde X' \cdot (\tilde E +
\sum j\cdot \tilde E_j))
$$
so we see that one of three things occurs: either

\

$\bullet$ \thinspace $X'$ is smooth, $k_i=1$ for all $i$, and $\tilde X'$
meets the last
exceptional divisor
$E_m$ transversely; or

\

$\bullet$ \thinspace $\tilde X'$ passes
through the point
$\tilde E_i \cap \tilde E_{i-1}$ for some $i < m$; or

\

$\bullet$ \thinspace for some $j < m$, $\tilde X'$ meets  the exceptional
divisor
$\tilde E_j$ at a point other than $\tilde
E_j \cap \tilde
E_{j-1}$ or $\tilde
E_j \cap \tilde
E_{j+1}$, and  has a point of intersection multiplicity
$m/j > 1$ with $\tilde E_j$.

\

We now compare the adjoint ideal $\I_X$ of $X$ with that of $\tilde
X$. The basic fact here is that if $C \subset S$ is any curve on a smooth
surface, $p \in C$ a
point of multiplicity $k$, and $\tilde C \subset \tilde S$ the proper
transform of $C$ in the
blow-up $\pi : \tilde S \to S$ of $S$ at $p$, the adjoint ideals of $C$ and
$\tilde C$ are
related by the formula
$$
\pi^*\I_C \; = \; \I_{\tilde C}(-(m-1)E)
$$
where $E$ is the exceptional divisor. Applying this $m$ times to
the curve $X$, we have
$$
\pi^*\I_X \; = \; \I_{\tilde X}(-\sum (l_j-j)\tilde E_j)\, .
$$

Now, $[X] \in W$ being general, any deformation of $X$ coming from the
family $W$
preserves the multiplicities $k_i$, and hence the decomposition $\pi^*X =
\tilde X + \sum l_i
\tilde E_i$. It also preserves the geometric genus of $\tilde X$, so that
identifying the space
$H^0(\tilde X, \O_{\tilde X}(\tilde X))$ of deformations of $\tilde X
\subset \tilde S$ with a
subspace of the deformations $H^0(X, \O_{X}(X))$  of
$X \subset S$ via the pullback map, we have
\begin{equation*}
\begin{split}
T_{[X]}W \; &\subset \; H^0(\tilde X, \I_{\tilde X}(\tilde X)) \\
&= \; H^0(\tilde X,
(\pi^*\I_{\tilde X})(\sum (l_j-j)\tilde E_j)(\pi^*X - \sum l_j \tilde E_j)) \\
&=  \; H^0(\tilde X, (\pi^*\I_{X})(\pi^*X-\sum j
\tilde E_j)) \\
&=  \; H^0(\tilde X, \pi^*(\I_{X}(X))(-l_mq)) \\
&= \; H^0(X, \I_X(-l_mp)) \\
&\subset \; H^0(X, \I_X(-mp)) \, .
\end{split}
\end{equation*}
Note that the inclusion in the last line of the above sequence is proper if
$X \ne X'$. Now,
suppose that
$X = X'$ is not smooth at
$p$. In this case, as we noted
$\tilde X'$ will  either be singular at $\tilde p$ or be tangent to $\tilde
E_i$ there, or else
will pass through the point $\tilde E_i \cap \tilde E_{i-1}$ for some $i$.
In the first case,
since
$\tilde X$ has a unibranch singularity, its deformations correspond to
sections of
$H^0(\tilde X, {\cal K}(\tilde X))$ for some ideal $\cal K$ strictly
contained in the adjoint
ideal $\I_{\tilde X}$; while in the latter two cases the deformations
correspond to sections
of
$H^0(\tilde X, \I_{\tilde X}(\tilde X))$ vanishing at $q$. In either case,
the inclusion in the
first line of the equation above is strict. Thus $T_{[X]}W \subset H^0(X,
\I_X(-(m+1)q))$
unless
$X$ is smooth at $p$, and  the remainder of the statement of the Lemma follows.
\end{pf}

To  conclude the proof of Proposition
\ref{dimensioncount} we need one more fact. To state it, let $X
\in |D|$ be any irreducible rational curve,
$\nu : X^\nu \to X$ the normalization and $p_1,p_2,\ldots
\in X^\nu$ any points; suppose that the  divisor $\nu^*(E)$ has
multiplicity
$m_i$ at $p_i$. Let
$\I \subset \O_S$ be the adjoint ideal of $X$, and set
$$
{\cal K} \; = \;  \I(-\sum m_i p_i)  \subset \O_X
$$
Let ${\cal K}'$ be any ideal of index 2 or less  in
${\cal K}$---that is, any ideal ${\cal K}' \subset {\cal K}$ with
$h^0({\cal K}/{\cal K}') \le 2$, or equivalently an ideal of the form
$$
{\cal K}' \; = \; {\cal K}(-q-r)
$$
for some pair of points $q, r \in X^\nu$. We will need
these ideals ${\cal K}' \subset
{\cal K}$ of index 2 in order to see, for example, that a general curve $X
\in V(D)$ does not
have a node on $E$. In these terms, our result is the
\begin{lm}\label{indcons}
 The ideal $ {\cal K}'$ imposes independent conditions
on the linear series $|\O_X(X)|$, i.e.,
$$
h^0(X, {\cal K'}(X)) \; = \; h^0(X, \O_X(X)) -
\dim_{\C } (\O_X/{\cal K'})
$$
In particular, ${\cal K}$ imposes independent conditions on $|\O_X(X)|$, that
is,
\begin{equation*}
h^0(X, {\cal K}(X)) \; = \;  r_0(D) - \sum m_i .
\end{equation*}
\end{lm}
\begin{pf}
By the adjunction formula we have
$$
K_{X^\nu} \; = \; \nu^*( K_S \otimes \O_S(X) \otimes {\cal I}) .
$$
 Thus,
$$
\nu^*(\O_S(X) \otimes {\cal K}) \; = \; K_{X^\nu} \otimes
\nu^*(\O_S(-K_S)) \otimes \O_{X^\nu}(-\sum m_i p_i) .
$$
Now, $\nu^*E - \sum m_i p_i \ge 0$, and on $S = \F_n$, we have
$$
K_S \; = \; \O_S(-C-E-2F)
$$
so that we have an inequality of divisor classes
$$
\nu^*(\O_S(X) \otimes {\cal K}) \; \ge \; K_{X^\nu} \otimes \nu^*\O_S(C+2F) .
$$
Moreover, the divisor class $C+2F$ has intersection number at least 3 with any
irreducible curve $X$ not linear equivalent to either $F$ or $E$, so it
follows that
$$
\deg(\nu^*(\O_S(X) \otimes {\cal K})) \; \ge \; -2 + 3 = 1 .
$$
Thus
$$
\deg(\nu^*(\O_S(X) \otimes {\cal K}')) \; \ge -1 .
$$
so that $h^1(X^\nu, \nu^*(\O_S(X) \otimes {\cal K}'))=0$, and the result
follows.
\end{pf}

We can now  complete the proof of Proposition \ref{dimensioncount}. We have
already
established, in the Claim above, that
$$
\dim(\tilde V_{\ul m}(D)) \; \ge \; r_0(D) - \sum (m_i-1) ;
$$
but applying Lemmas \ref{tangencycondition} and \ref{indcons} in turn we
see that for any subset $\Omega =
\{p_1,\ldots,p_k\} \subset E$,
\begin{equation*}
\begin{split}
\dim(\tilde W^\Omega_{\ul m}(D)) \; &\le \; h^0(X, {\cal K}(X)) \\
&= \; r_0(D) - \sum m_i
\end{split}
\end{equation*}
and hence
\begin{equation*}
\begin{split}
\dim(\tilde V_{\ul m}(D)) \; &\le \; \dim(\tilde W^\Omega_{\ul m}(D)) + k \\
&= \; r_0(D) - \sum (m_i-1)
\end{split}
\end{equation*}
so that equality must hold. Moreover, if a general point $[X] \in V_{\ul m}(D)$
corresponded to a curve $X$ with singularities other than nodes, the second
inequality above would be
strict; so $X$ must be nodal, and smooth at ots points of intersection with $E$.

\

We can eliminate all the other possible misbehaviors of our general curve
$X$ similarly.  If the point $p \in X^\nu$ is
mapped to one of the points $P_i$, we
would have
\begin{equation*}
\begin{split}
\dim(\tilde V_{\ul m}(D)) \; &\le \; h^0(X, {\cal K}(X)(-p)) \\
&< \; h^0(X, {\cal K}(X)) ;
\end{split}
\end{equation*}
and if the multiplicity of the pullback divisor $\nu^*(G)$ at $p$ were $m >
1$ we
would have
\begin{equation*}
\begin{split}
\dim(\tilde V_{\ul m}(D)) \; &\le \; h^0(X, {\cal K}(X)(-(m-1)p)) \\
&< \; h^0(X, {\cal K}(X)) .
\end{split}
\end{equation*}
Suppose next that $X$ had a node on $E$, with branches corresponding to a
pair of points $q, r \in X^\nu$ and the branch corresponding to $r$
transverse to
$E$. It would follow that
$$
h^0(X, {\cal K}(X)(-q-r)) \; = \; h^0(X, {\cal K}(X)) - 1 ,
$$
since a section of ${\cal K}(X)$ vanishing at $q$ but not at $r$ would
correspond to a deformation of $X$ in $\tilde V_{\ul m}(D)$ in which the two
branches would meet $E$ in distinct points.

Finally, to prove part 3 of  Proposition \ref{dimensioncount}, we simply
let $X'$ be a
general member of the family $\tilde V_{\ul m'}(D')$ and apply the above to $X
\in \tilde V_{\ul m}(D)$, including $X'$ in $G$ and its points of
intersection with
$G$ and $E$ among the points $P_i$.
\end{pf}

The next Proposition is
stated as a characterization of the reducible elements of the
one-parameter family $\X^\nu
\to \Gamma$, but in fact it is a characterization of the codimension one
components of the boundary
$V(D)
\setminus
\tilde V(D)$ of $V(D)$.
\begin{prop}\label{codimensionone} Let $X \subset S$ be any reducible
fiber of the family $\X
\to \Gamma$.

1.  If $X$ does not contain $E$, then $X$ has exactly two irreducible
components
$X_1$ and
$X_2$,  with $[X_i]\in V(D_i)$ and $D_1+D_2=D$.
Moreover $[X_i]$  is a general point in
 $V(D_i)$.

2.  If $X$ does  contain $E$, then $X$  has  irreducible components
$E$, $X_1,\ldots,X_k$,  with $[X_i]\in V(D_i)$ and $E+D_1+\ldots +D_k=D$.
 Moreover each  $X_i$ is general  in
$V_{m_i}(D_i)$ for some collection $m_1,\dots,m_k$ of positive integers
such that
$\sum (m_i-1) = n-k$.

\end{prop}

\

\ni \underbar{Remark}.
Notice that by Proposition~\ref{dimensioncount}, the above result says
that if $X$ does not contain $E$, then it has only nodes as singularities.
And, if $X$ contains $E$,  away from the $k$ points of tangency of $E$ with
the curves $X_i$, $X$ has
only nodes as singularities.

\
\ni
 \begin{pf}  Assume first that $X$ does not contain $E$. Write the divisor
$X$ as a
sum
$$ X = \sum_{i=1}^k a_i \cdot X_i
$$  where $a_i > 0$ and the $X_i$ are irreducible curves in $S$. We claim
first that since
$[X]
\in V(D)$, all the curves
$X_i$ must be rational. To see this, take any one-parameter family
${\cal X}
\to B$ of irreducible rational curves specializing to $X$.
Proceeding as in \ref{setup} we arrive at a family $\Y \to B$ of nodal
curves, with general fiber
$\P^1$, that admits a regular map $\Y \to {\cal X}$. Now, since the
fibers of
$\Y
\to B$ are reduced curves of arithmetic genus 0, every component of
every fiber of
$\Y$ must be a rational curve. Thus every component of $X$ is
dominated by a rational curve and so must be itself rational.

Thus $[X_i] \in V(D_i)$, where $D_i$ are divisor classes such that $\sum
a_iD_i = D$. On the other hand, since $X$ is a general member of an
$(r_0(D)-1)$-dimensional family, we must have
$$
\sum_{i=1}^k r_0(D_i) \; \ge \; r_0(D) - 1
$$ which yields
\begin{equation*}
\begin{split}
\sum_{i=1}^k \bigl( -(K_S \cdot D_i) - 1 \bigr) \; &\ge \; (-K_S \cdot D) - 2 \\
&= \;
\sum_{i=1}^ka_i (-K_S \cdot D_i) - 2 .
\end{split}
\end{equation*} Comparing the two sides, we see that
$$ 2 - k - \sum_{i=1}^k(a_i-1) (-K_S \cdot D_i) \; \ge \; 0 .
$$ But
$(-K_S\cdot D_i)\geq 2$ for  any curve $D_i$ on $S$ other than $E$; so we may
conclude that all $a_i = 1$ and that $k \le 2$. Moreover, if $k=2$ we
have equality in the above inequality, which says that the pair of curves
$(X_1, X_2)$ is general in
$V(D_1) \times V(D_2)$.

We come now to the case where $X$ contains $E$. The first thing we see here is
that the dimension-count argument we used above doesn't work: since
$$
(-K_S \cdot (X-aE)) \; = \; (-K_S \cdot
X) + a(n-2) ,
$$
the sums $\sum a_iX_i$ of rational curves $X_i \in |D_i|$ may well  move in a
larger-dimensional family than $X$ itself.

The key here is to look  at the
semistable reduction of a family of  curves in $ \tilde V(D)$
specializing to $X$. This will allow us to limit the number of points of
intersection
of the curves $X_i$ with $E$, that is to say, to show that in fact the
$X_i$ belong
to $V_{\ul{m}}(D_i)$ for suitable $\ul m$. This replaces the naive bound
above on
the dimension of the family of such curves $X$ with a stronger one, which turns
out to be sharp.

Consider then the family $\Y \la B$ obtained from $\X \la \Gamma$ as in
Section~\ref{setup}.
We can thus assume that the total space $\Y$ of the family is smooth and
every fiber
of $\Y$ will be a union of smooth rational curves meeting transversely, and
whose dual graph is a tree.

Now, let
$Y$ be the special fiber of
$\Y \to B$. We  decompose $Y$ into two parts: we let
$Y_E$ be the union of the irreducible components of $Y$ mapping to
$E$, and $Y_R$  the union of the remaining components. Next, we decompose
$Y_R$ further into $k$ parts, letting $Y_i$ be the union of the components
mapping to
$X_i$.
 Denote the connected components of
$Y_E$ by
$Z_i$, and for each $i$ let
$\a_i$ be the degree of the map $\mu|_{Y_i} : Z_i \to E$, so that $\sum
\a_i = a$. Similarly, let $\{Z_{i,j}\}_j$ be the
connected components of $Y_i$ and $\a_{i,j}$ the degree
of the restriction $\mu|_{Z_{i,j}} : Y_{i,j}
\to X_i$, so that $\sum_j \a_{i,j} = a_i$.

Note that    the inverse image of $E$ in $\Y$ is given by
$$
\pi^{-1}(E) \; = \; Y_E \cup \Gamma_1 \cup \dots \cup \Gamma_b .
$$
(Where $\pi : \Y \to S$ is, as usual,  the natural map.)

As we indicated, the essential new aspect of the argument in this case is
keeping
track of the number of points of intersection of the $X_i$ with $E$. To do
this, we
 note that, over any such point, there will be a point of intersection of a
component of $Y_i$ with the inverse image $\pi^{-1}(E)$; which by the
expression above for $\pi^{-1}(E)$ will be either a point of intersection
of $Y_i$
with $Y_E$ or one of the $b$ points of intersection of the $\Gamma_i$ with $Y$.

It thus remains to  bound
 the number $\e$ of points of intersection of $Y_E$ with the remaining parts
$Y_i$ of $Y$. This we can do by using the fact that the dual graph of $Y$
is a tree:
this says that the number of pairwise points of intersection of the connected
components $Z_{i,j}$ of $Y_i$ and the connected components $Z_i$ of $Y_E$ is
equal to the total number of all such connected components, minus one. Thus,
\begin{equation*}
\begin{split}
\e \; = \; \#(Y_R \cap Y_E) \; = \; &\#\{ \text{connected components of }
Y_E\} \cr
+ \; &\sum \#\{ \text{connected components of } Y_i\} .
\end{split}
\end{equation*}
Note that the degree $\a_i > 0$ on each component $Z_i$ of $Y_E$, so that
$$
\#\{ \text{connected components of } Y_E\} \; \le \; a
$$
and similarly
$$
\#\{ \text{connected components of } Y_i\} \; \le \; a_i .
$$
Thus we can deduce in particular that
$$
\e \; \le \; a + \sum a_i - 1 .
$$

Now, say $X_i \in \tilde V_{{\ul m}^i}(D_i)$ for each $i = 1,\ldots,k$. Let
$\nu_i : X_i^\nu \to X_i$ be the normalization map. Choose any irreducible
component
$X_i^0$ of
$Y$ dominating
$X_i$ (and hence dominating the normalization $X_i^\nu$), and let
$\pi_i : X_i^\nu
\to X_i$ be the restriction of $\pi$ to $X_i^\nu$. Trivially, the total
number of
points of the pullback $\nu_i^*(E)$ of $E$ to $X_i^\nu$ is
\begin{equation*}
\begin{split}
\# \nu_i^*(E) \; &\le \; \#\pi_i^*(E) \\
&= \; \#(X_i^0 \cap Y_E)
\end{split}
\end{equation*}
and hence
\begin{equation*}
\begin{split}
\sum \# \nu_i^*(E) \; &\le \; \sum \#(X_i^0 \cap Y_E) \\
&\le \; \#(Y_R \cap Y_E) \\
&= \; \e
\end{split}
\end{equation*}
with strict inequality if any $a_i > 1$. But the sum of degrees of $E$ on the
curves $X_i$ is at least
\begin{equation*}
\begin{split}
\sum \deg(\pi_i^*E) \; &\ge \; \left((\sum X_i) \cdot E\right) \\
&= \; \left( (D - aE - \sum (a_i - 1)D_i ) \cdot E \right) \\
&= \;  (D \cdot E) + an - \sum a_i(D_i \cdot E) .
\end{split}
\end{equation*}
Comparing the number of points of the pullbacks of $E$ to the normalizations
$X_i^\nu$ with the degrees of these pullbacks, we conclude that there must be
multiplicities in these divisors: specifically, the sum $\sum (m_j^i - 1)$
of the
multiplicities minus one must be the difference of these numbers, so that
\begin{equation*}
\begin{split}
\sum (m_j^i - 1) \; &\ge \; \sum \deg \pi_i^*(E) - \e - (D \cdot E) \\
&\ge \; (D \cdot E) + an - \sum (a_i - 1)(D_i \cdot E) - a - \sum a_i + 1 -
(D \cdot
E) \\
&\ge \;  a(n-1) - \sum (a_i - 1)(D_i \cdot E)  - \sum a_i + 1 .
\end{split}
\end{equation*}
This in turn allows us to bound the number of degrees of freedom of the curves
$X_i$: we have
\begin{equation*}
\begin{split}
\sum \dim \tilde V_{{\ul m}^i}(D_i) \; &= \; \sum r_0(D_i) - \sum (m_j^i - 1) \\
&= \; \sum_{i=1}^k \left( (-K_S \cdot D_i) - 1 \right) - \sum (m_j^i - 1) \\
&\le \; \sum (-K_S \cdot D_i) - k - a(n-1) + \sum (a_i - 1)(D_i \cdot E) +
\sum a_i
- 1 .
\end{split}
\end{equation*}
On the other hand, this must be at least equal to the dimension of $V(D)$ minus
one, that is,
\begin{equation*}
\begin{split}
r_0(D) - 1 \; &= \; (-K_S \cdot D) - 2 \\
&= \; a(-K_S \cdot E) + \sum a_i(-K_S \cdot D_i) - 2 \\
&= \; a(n-2) + \sum a_i(-K_S \cdot D_i) - 2 .
\end{split}
\end{equation*}
In the end, then, we must have
\begin{equation*}
\begin{split}
a(n-2) + \sum a_i(-K_S \cdot D_i) - 2 \; &\le \; \sum (-K_S \cdot D_i) - k
- a(n-1)
\\
 &\quad + \sum (a_i - 1)(D_i \cdot E) + \sum a_i - 1 .
\end{split}
\end{equation*}
We can (partially) cancel the $a(n-1)$ and $a(n-2)$ terms, and combine the
terms involving $(-K_S \cdot D_i)$ to rewrite this as
$$
a + \sum (a_i-1)(-K_S \cdot D_i) - 1 \; \le \; \sum (a_i - 1)(D_i \cdot E) - k +
\sum a_i - 1
$$
or, in other words,
$$
a + \sum (a_i - 1) \left[ \left( (-K_S - E) \cdot D_i \right) - 1 \right] -
1 \; \le \; 0 .
$$

Now, we have already observed that
$
-K_S - E \; = \; C + 2F
$
meets every curve $X_i$ strictly positively, so that the sum in this last
expression is nonnegative. We conclude that $a=1$, and (since any $a_i > 1$
would have led to strict inequality) that all $a_i = 1$. Next, since there
is a unique component of $Y$ mapping to each $X_i$, each curve $X_i$ will have
at most one point of intersection multiplicity $m > 1$ with $E$. Thus, finally,
$X_i$ is a general member of the family $\tilde V_m(D_i)$ for some collection of
integers $m_1,\dots,m_k$ with $\sum (m_i - 1) = n-k$, completing the proof of
Proposition \ref{codimensionone}.
 \end{pf}

Note that we have not said here that every reducible curve satisfying the
conditions of the Proposition  in fact lies in the closure of the locus of
irreducible
rational curves. This is  true,
 and is not hard to see in the case of curves of types (1); but for
curves of type (2) it is a deeper fact, and we will require
the proof of Proposition~\ref{singtotalspace} to establish it.

Having characterized as a set the locus $\Gamma$ of curves  in $V(D)$
passing through
$q_1,\dots,q_{r_o (D) -1}$, we  now turn to a statement   about the local
geometry of
$\Gamma$ around each point.

We  introduce one bit of terminology here. Let $X$ be a
 fiber of $\X \to \Gamma$; and, in case
$\Gamma$ is locally reducible at the point
$[X] \in \Gamma$, pick a branch of $\Gamma$ at $[X]$ (that is, a point
$b$ of the normalization $\Gamma^\nu$ of $\Gamma$ lying over $[X]$). Let
$P$ be a
node of $X$. We  then make the following

\

\ni {\bf Definition}. If $P$ is a limit of nodes of fibers of $\X \to
\Gamma$ near
$X$ in the chosen branch---that is, if $(P,b)$ is in the closure of the
singular locus of the map $\X \times_\Gamma (\Gamma^\nu
\setminus \{b\}) \to \Gamma^\nu$---we will say that $P$ is an {\it old}
node of $X$.
If
$(P,b)$ is an isolated singular point of the map $X \times_\Gamma (\Gamma^\nu
\setminus \{b\}) \to \Gamma^\nu$ we will say that $P$ is a {\it new} node
of $X$.

\

\ni Equivalently, $P$ is an old node if the fiber
$X^\nu$ of $\X^\nu \to \Gamma^\nu$ over $b$ is smooth at
the (two) points  lying over $P$; if it is a new node, $X^\nu$ will have
a single point lying over $P$, which will be a node of $X^\nu$.

\

Note that if $P$ is a singular point of $X$ other than a node, the situation
is not so black-and-white. For example, if $P$ is an $m$-fold
tacnode---that is, if the curve $X$ has two smooth branches at
$P$ with contact of order $m$---then a priori, any number $n \le m$ of
nodes of nearby fibers may approach
$P$ along any branch of $\Gamma$ at $[X]$, with the result that the
fiber of $\X^\nu
\to \Gamma^\nu$ over the corresponding point $b \in \Gamma^\nu$ will have an
$(m-n)$-fold tacnode over
$P$, or will be smooth over $P$ if $n=m$. (The proof of  the
relevant case $n=m-1$ will emerge in the proof of
Proposition \ref{singtotalspace}.)

In these terms, we can state

\begin{prop}\label{describegamma} Let $X$ be a reducible fiber of the
family
$\X \to \Gamma$. Keeping the notations and hypotheses of Proposition
\ref{codimensionone},

$1$.  If $X = X_1 \cup X_2$ does not contain $E$, and $X_1$ and $X_2$ meet at
$(D_1 \cdot D_2) = \ell$ points $P_1,\dots,P_\ell$, then in a
neighborhood of
$[X]$  $\Gamma$ has $\ell$ smooth branches $\Gamma_1 ,\ldots ,\Gamma _{\ell }$; along
$\Gamma_i$ the point $P_i$ is new, and all other nodes of $X$ are old.

$2a$.  If $X = E \cup X_1 \cup \dots \cup X_k$, and  $X_i$ meets $E$
transversely in
$(D_i \cdot E) = \ell_i$ points $P_{i,1},\dots,P_{i,\ell_i}$, then in a
neighborhood of
$[X]$
$\Gamma$ consists of
$\prod \ell_i$ smooth branches $\Gamma_\a =
\Gamma_{(\a_1,\dots,\a_k)}$. Along
$\Gamma_\a$ the points $P_{1,\a_1},\dots,P_{k,\a_k}$ are new, and all
other nodes of
$X$ are old.

$2b$.  If $X = E \cup X_1 \cup \dots \cup X_k$, and  $X_i$ meets $E$
transversely in
$(D_i \cdot E) = \ell_i$ points $P_{i,1},\dots,P_{i,\ell_i}$ for $i =
2,\dots,k$, while $D_1$ has a
 point $P$ of intersection multiplicity $m \ge 2$ with $E$, then in a
neighborhood of
$[X]$
$\Gamma$ consists of
$\prod_{i=2}^k \ell _i$ smooth branches $\Gamma _\a =
\Gamma_{(\a_2,\dots,\a_k)}$. Along
$\Gamma_\a$ the points $P_{2,\a_2},\dots,P_{k,\a_k}$ are new; all
other nodes of
$X$ are old; and exactly $m-1$ nodes of nearby fibers will tend to $P$.
\end{prop}

\ni \underbar{Remark 1}. The proof of this Proposition will not be
complete until the end of the following section.  More precisely, we will
postpone the proof of the existence and smoothness of the branches
of $\Gamma$.
Actually, cases $1$ and $2a$ could very well be proved here, but it is more
convenient do it later (that is, at the beginning of the proof of Proposition
\ref{singtotalspace}).

\ni \underbar{Remark 2}.
We believe that an analogous description of the family
$\X \to \Gamma$ may be given without the assumption that the components of
the curve
$X$ other than
$E$ have altogether at most one point of tangency with $E$, and otherwise
intersect $E$
transversely in distinct points. The
restricted statement above will suffice for our present purposes. We hope
to prove the
general statement in  the future.

\ni \underbar{Remark 3}. The statement of Proposition \ref{describegamma}
can also be expressed in terms of the normalized family
$\X^\nu
\to
\Gamma^\nu$, and indeed that is how we will use it in the following chapter.
In these terms, the
statements are:

\

\ni $1$. \thinspace If $[X]$ is a point of
$\Gamma$ corresponding to a curve $X$ in our family not containing
$E$, then there will be
$(D_1 \cdot D_2) = \ell$ points of $\Gamma^\nu$ lying over $[X]$, corresponding
naturally to the nodes of $X$.  The fibers of
$\X^\nu
\to \Gamma^\nu$ over these points will be the normalizations of $X$ at all the
nodes of $D_1$ and $D_2$ and at all but one of the
$\ell$ points of intersection of $D_1$ with $D_2$.

\

\ni $2a$. \thinspace  If $X = E + D_1 + \dots + D_k$ contains $E$ and
the components
$D_i$ intersect
$E$ transversely, then the fibers of $\X^\nu
\to \Gamma^\nu$ over points lying over $[X] \in \Gamma$ are  the curves obtained
by normalizing
$X$ at all nodes of the $D_i$, at all the points of pairwise intersection of
the
$D_i$, and at all but one of the points of intersection of $E$ with each of
the components
$D_i$. In other words, the  fibers consist of the disjoint union of the
normalizations
$\tilde D_i$ of the curves $D_i$, each attached to $E$ at one point.

\

\ni $2b$. \thinspace  If $X = E + D_1 + \dots + D_k$ as before
and one of the components
$D_1$ of $X$ has a smooth point $P$ of intersection multiplicity $m \ge
2$ with $E$, then the fibers $X^\nu$ of $\X^\nu
\to \Gamma^\nu$ corresponding tor $[X] \in \Gamma$ are
the curves obtained by normalizing
$X$ at all nodes of the $D_i$, at all the points of pairwise intersection of
the
$D_i$, at all but one of the points of intersection of $E$ with each of the
components
$D_i$ for $i = 2,\dots,k$, at all the transverse points of intersection of
$D_1$ with
$E$, and finally taking the partial normalization of $X$ at $P$ having an
ordinary node over $P$.  (The fact that each fiber of $\X^\nu \to
\Gamma^\nu$ lying
over $X$ has an ordinary node over $P$ follows either from the  fact
that the
$\delta$-invariant of the singularity $P \in X$ is $m$ and that, along
each branch,
$m-1$ nodes of nearby fibers tend to $P$; or---what is essentially the
same thing---the fact that the arithmetic genus of the fibers of $\X^\nu
\to \Gamma^\nu$ are zero. This will
 be verified independently in the course of the proof of
Proposition
\ref{singtotalspace}.) The picture is  therefore
similar to the preceding case: the fibers consist of the
disjoint union of the normalizations $\tilde D_i$ of the curves
$D_i$, each attached to $E$ at one point. The one difference is that, while
for $i = 2,\dots,k$ the point of attachment of the normalizations $\tilde
D_i$ with $E$ can lie over any of the points of intersection of $D_i$ with
$E$,   the point of intersection of the normalization of $D_1$ with $E$ can
only be the point lying over $P$.

\

A typical picture of the original curve $X$ and its partial normalization
$X^\nu$ is this:

\

\vspace*{5.3in}

\ni \special {picture basicdiagram}

\

\begin{pf} Consider first of all a reducible curve $X$ in our family
that does not contain $E$. By Proposition \ref{codimensionone}, this must
be  of the form $X = X_1 \cup X_2$ where
$X_i$ is a general member of the family $V(D_i)$ with
$D_1 +D_2 =D$. In particular, $X_i$ is an irreducible
rational curve with
$p_a(D_i)$ nodes, and $X_1$ and $X_2$ intersect transversely in $(D_1
\cdot D_2)$ points.  Note that
$$ p_a(D_i) = \frac{(D_i \cdot D_i) + (D_i \cdot K_S)}{2} + 1
$$ so that the total number of nodes of $X$ will be
 \begin{equation*}
p_a(D_1) + p_a(D_2) + (D_1 \cdot D_2) \; = \; p_a(D) + 1 .
 \end{equation*} In other words,  along any branch of $\Gamma$, all but one
of the
nodes of $X$ will be limits of nodes of nearby fibers
 (that is, will be old nodes), while one node of $X$ will be a
new node. Note also that not any node of $X$ can be the new node:
that must be one of the points of intersection of the two components
$X_1$ and $X_2$; otherwise the fiber of the normalization $\X^\nu$
would be disconnected.

In case $X$ contains $E$,
the analogous computation yields that   $X$ has
$p_a(D)+k$ nodes (or $p_a(D)+k-m$ nodes and one tacnode of order $m$ in
case $2b$); hence
$X$ has $k$ new nodes (or, $k-1$ in $2b$). Then
the
analysis in the proof of Proposition \ref{codimensionone} shows that
 in the normalization of the
total space of the family, the corresponding fiber will consist of a curve
$\tilde E$ mapping to
$E$, plus the normalizations $\tilde X_i$ of the curves $X_i$, each meeting
$\tilde E$ in one
point and disjoint from each other. In particular, all the nodes of $X$
arising from points of
pairwise intersection of the components
$X_i$ are  old. As for the points of intersection of the
components $X_i$ with $E$, there are two cases. First, if a component $X_i$
has a point of
contact of order $m > 1$ with $E$, that must be the image of the point
$\tilde X_i \cap
\tilde E \in \X^\nu$; and all the other points of $X_i \cap E$ will be old
nodes of $X$ on
any branch. On the other hand, if a component $X_i$ intersects $E$
transversely, any
one of its points of intersection with $E$ can be a new node.

\end{pf}

\subsection{Singularities of the total space}

We come finally to the fourth result, in which we will describe the
singularities of the total space of the normalized family
$\X^\nu
\to \Gamma^\nu$ along a given fiber $X^{\nu }$.  (Given a fiber $X$ over
$\Gamma$,
we will fix a corresponding fiber  $X^{\nu }$ throughout.)

We will keep a simplified form of the notation introduced in the
statement of Proposition
\ref{describegamma}:  we will denote by
$P_1,\dots,P_\ell$ the new nodes of $X$
 along $E$, coming from transverse points of intersection of other
components of
$X$ with
$E$; and by $P$ (if it exists) one double point of $X$  other than a
node, coming from a point of contact of order $m\geq 2$ of
$E$ with another component of $X$. We recall that  the nearby fibers
 of our family are smooth near $P_i$, there will be one point $p_i$
of
$\X^\nu$ lying over each $P_i$, which will be a node of $X^\nu$, while
the nearby fibers have $m-1$ nodes tending to the
point $P$, so that the partial normalization $X^\nu \to X$
 will again have one point $p$  lying over
$P$, and that point will be a node of
$X^\nu$.
With all this said, we have
\begin{prop}\label{singtotalspace}

1. If $X$ does not contain $E$, or if $X$ contains $E$ and the closure of $X
\setminus E$ intersects $E$ transversely, then
$\X^\nu$ is smooth along $X^\nu$.

2. In case $X$ does contain $E$  and the closure of $X
\setminus E$ has a point $P$ of intersection multiplicity $m \ge 2$ with
$E$, the point  $p$ of $X^\nu$ lying over
$P$ is a smooth point of
$\X^\nu$;  the other nodes  $p_i$ of $X^\nu$ will be singularities of type
$A_{m-1}$ of
$\X^\nu$ .
\end{prop}
\begin{pf}  We start with the first statement, which is by far the
easier.
Recall that by the  two previous propositions $X$,
being a general point on a codimension-one locus in $V(D)$, will have
$p_a(D)+k$ or
$p_a(D)+1$ nodes, depending whether
$X$ does or doesn't contain $E$. Of these, $p_a(D)$
will be old nodes
 and the remaining ones are new nodes;  if $E$ is contained in $X$, then the
new nodes all lie on $E$.
Let
$r_1,\dots,r_{p_a(D)}$ be the old nodes of
$X$ and let
$P$ be  any fixed new node. The fiber $X^\nu$ of
$\X^\nu$ lying over $X$ will be the partial normalization of $X$ at
$r_1,\dots,r_{p_a(D)}$, so that $\X^\nu$ will certainly be smooth there,
and we need only concern ourselves with the point of $\X^\nu$ lying
over $P$.

Consider, in an analytic neighboroohd of $[X]$ in $|D|$,  the locus $W$ of
curves
that pass  through the base points $q_1,\dots,q_{r_0(D)-1}$  and that
preserve all
of the old nodes of $X$.
The  projective tangent space to $W$ at $[X]$
 will be contained in the sub-linear series of $|D|$
of curves passing through the $p_a(D)$ old
nodes of $X$ and through
$q_1,\dots,q_{r_0(D)-1}$. This gives a total of
$r_0(D)-1 + p_a(D)=r(D)-1$ points which,
by an argument  analogous to the proof of
Lemma~\ref{indcons},
  impose independent conditions on the linear series
$|D|$. We only exhibit the proof in case $E$ is a component of $X$, the
other case being similar and easier.
Let $\cal H$ be the ideal sheaf of the subscheme of $S$ given by
the
old nodes  $r_1,\dots,r_{p_a(D)}$, and let $\nu :\tilde X\la X$ be the
normalization map.
 We have to show that $r_1,\dots,r_{p_a(D)}$ impose independent
conditions on
$|D|$, which will follow (cf. Lemma~\ref{indcons}) from
$$
H^1\bigl( \tilde X, \nu ^* ( \O _S(X)\otimes \cal H )\bigr) =0.
$$
This, by the adjunction formula, is equivalent to
$$
H^0 \bigl( \tilde X , \nu ^* ( K_S \otimes \I )
\otimes (\nu ^*\cal H )^{-1}\bigr) =0
$$
where $\I$ is the adjoint ideal of $X$.
Now notice that the line bundle $\nu ^* (\I )\otimes \nu ^* (\cal H
)^{-1}$ has degree $-k$ on the component of $\tilde X$ lying over $E$, and
degree
$-1$ on every other component.
Since $ K_S$ has degree $n-2 = k-2$ on $E$ and negative degree on $X_i$,
the line
bundle
$\nu ^* ( K_S \otimes \I )
\otimes (\nu ^*\cal H )^{-1}$
cannot have any
sections.

We  conclude that $W$ is smooth of dimension 1.
Notice that this  completes the proof  of Proposition~\ref{describegamma},
parts $1$
and $ 2a$.

To anlyze the total space of $\X ^\nu $ we consider the
 map from $W$ to
the versal deformation space of the node $(X,P)$.
This  has nonzero differential
 because $P$ is not a base point of
the linear series of curves   passing through
$q_1,\dots,q_{r_0(D)-1}$ and
through the $p_a(D)$ old nodes of $X$
(to see this, the argument above applied to the ideal sheaf of
the union of the old nodes of $X$ and $P$ will work).
  Thus the family
$\X^\nu \to \Gamma^\nu$ has local equation $xy-t=0$ near $p$; in
particular, it is
smooth at $p$.

We turn now to the second part, which will occupy us for the remainder of this
chapter. We will start by carrying out a global analysis of the family in a
neighborhood of the whole fiber
$X$, and then proceed to a local analysis around the point $P$
specifically. From the global
picture we will establish that, for some integer $\gamma$, the point
$P$ will be a singularity of type $A_\gamma$ and the points $P_i$ all
singularities
of type
$A_{\gamma m}$. The local analysis will then show that in fact we have
$\gamma=1$.

To carry out the global analysis, we use the family $\Y \to
\Gamma^\nu$ and the map $\pi :\Y \to \F _n$ (cf. section~\ref{setup}),
where
$\Y$ is the minimal desingularization of the surface $\X^\nu$. Since the
singularities of the fiber of $\X^\nu$ are all nodes, the total space
$\X^\nu$ will have singularities of type $A_{\beta}$ at each; let us say the
point $p$ is an
$A_\gamma$ singularity of $\X^\nu$, and the point $p_i$ an $A_{\gamma_i}$
singularity. When we resolve the singularity at $p$ we get a chain
$G_1,\dots,G_{\gamma-1}$ of smooth rational curves; likewise,
$p_i$  is replaced by a chain
$G_{i,1},\dots,G_{i, \gamma_i-1}$ of smooth rational curves. Denoting the
component of $X$
meeting $E$ at $P_i$ by
$D_i$ and the component meeting $E$ at $P$ by $D$ (we are not assuming
here that these are distinct irreducible components of $X$), we arrive at
a picture of the relevant part of the fiber
$Y$ of $Y$:

\vspace*{3.4in}

\ni \special {picture Fnmult}

\

\

(We hope that such a notation will not be too confusing!)

We now look at the pull-back of $E$ from $\F _n$ to $\Y$.
 We can write it as
$$
\pi^*(E) = k \cdot E + \sum a_i \cdot G_i +  \sum a_{i,j} \cdot G_{i,j} + E'
$$
where $E'$ is a curve in $\Y$ that meets the fiber $Y$ only along   $D_i$ and
$D$, with $(E'\cdot D_i) = (E\cdot \pi (D_i)) -1$ and $(E'\cdot D) =
(E\cdot \pi (D)) -m$.

We can use what we know about the degree of this divisor on the
various components of $Y$ to impose conditions on the
coefficients $k$,
$a_i$ and $a_{i,j}$.
First, since $\pi$ maps
components
$G_i$ and
$G_{i,j}$  to points in $\F_n$,
$$
\deg_{G_i}(\pi^*(E)) = \deg_{G_{i,j}}(\pi^*(E)) = 0 .
$$ Now, each of the  curves $G_i$ and $G_{i,j}$ has self-intersection $-2$;
so, setting
$a_\gamma = a_{i,\gamma_i} = 0$ and
$a_0 = a_{i,0} = k$, we get
$$ a_{i-1} - 2a_i + a_{i+1} = 0
$$ for each $i = 1,\dots,\gamma-1$; and similarly
$$ a_{i,j-1} - 2a_{i,j} + a_{i,j+1} = 0
$$ for each $j = 1,\dots,\gamma_i-1$---in other words,
 the sequences $a_0,\dots,a_\gamma$ and
$a_{i,0},\dots,a_{i,\gamma_i}$ are arithmetic progressions. On the
other hand, the map $\pi$ restricted to the component $D_i$ is transverse
to $E$ at $P_i =
\pi(p_i)$; so the multiplicity at $p_i$ of the restriction to $D_i$ of the
divisor
$\pi^*(E) - E'$ is one. This says that $a_{i,\gamma_i-1} = 1$; and similarly
$a_{\gamma -1} = m$. Following the arithmetic progression
$a_0,\dots,a_\gamma$ up from
$D$ to $E$, we arrive at
$$ k = \gamma \cdot m
$$ and hence
$$
\gamma_i = \gamma \cdot m .
$$

The proof of the Proposition will be completed once we show that $\gamma
=1$, that is,
that $p$ is a smooth point of $\X ^{\nu }$.

Note that this part of the analysis did not rely, except notationally, on the
hypothesis that all but one point of intersection of $E$ with the remaining
components of $X$ are transverse. If the points $P_i$ were  points of
intersection multiplicity $m_i$ of $E$ with other components $D_i$ of
$X$, we could (always assuming that $m_i-1$ nodes of the general fiber
of our family approach $P_i$) carry out the same analysis and deduce that
for some integer $k$, the point $p_i$ was a singularity of type
$A_{k/m_i}$---loosely speaking, the singularity of $\X^\nu$ at $p_i$ is
``inversely proportional" to the order of contact of
$D_i$ with $E$ at $P_i$. The remaining question then
would be, is the number $k$ as small as possible, that is, the least common
multiple of the
$m_i$? That is what we will establish with the following local analysis,
which does ultimately rely on the hypothesis that all but one of the  $m_i$
are one.

\subsubsection{The versal deformation space of the tacnode}
\label{tacnode}
We now carry out the analysis around
the point
$P$.  The  versal deformation of $P\in X\subset \F _n$ has
the vector space $\O_{\F_n, P}/{\cal J}$ as base, where
${\cal J}$ is the Jacobian ideal of $X$ at $p$.
Choose local coordinates $x,y$ for $\F_n$
centered at $P$, so that the curve $E$ is given as $y=0$ and the
equation of $X$ is
$$ y(y+x^m) = y^2 + yx^m = 0
$$
The Jacobian ideal of this polynomial is ${\cal J} = (2y+x^m, yx^{m-1})$.
The monomials   $y, xy, x^2y,\dots,x^{m-2}y$ and $1, x, x^2,\dots,x^{m-1}$
form a basis for  $\O_{\F_n, P}/{\cal J}$, so that we can
write down explicitly a versal deformation space: the base $\Delta $ will be
an analytic  neighborhood
of the origin in affine space
${\Bbb A}^{2m-1}$ with coordinates  $\alpha_0,
\alpha_1,\dots,\alpha_{m-2}$ and $\beta_0,
\beta_1,\dots,\beta_{m-1}$, and the deformation space will be the family
${\cal S} \to \Delta$, with ${\cal S} \subset \Delta \times {\Bbb A}^2$,
 given by the equation
$$
 y^2 +yx^m   + \alpha_0y + \alpha_1xy
+ \dots +
\alpha_{m-2}x^{m-2}y + \beta_0 +
\beta_1x + \beta_2x^2 + \dots + \beta_{m-1}x^{m-1}  = 0
$$
Inside $\Delta$ we look closely  at the closures $\Delta_{m-1}$ and
$\Delta_m$ of the loci corresponding to curves with
$m-1$ and $m$ nodes, respectively.
We have
\begin{lm}
\label{versal}

1. $\Delta _m$ is given in
$\Delta$
by the equations $\beta _0=\ldots =\beta _{m-1} = 0$; in particular it
is smooth of dimension $m-1$.

2. $\Delta _{m-1}$ is irreducible of dimension $m$,
with $m$ sheets  crossing transversely at a general point of $\Delta_m$.
\end{lm}

 \begin{pf}
We introduce the {\em discriminant} of the
polynomial
$f$ above, viewed as a quadratic polynomial in $y$:
$$
\delta = \delta_{\a,\b}(x) = (x^m + \alpha_{m-2}x^{m-2} + \dots  +
\alpha_1x +
\alpha_0)^2 - 4(\beta_{m-1}x^{m-1} + \dots + \beta_1x + \beta_0)
$$ Note that the map $\delta : \Delta \to V$  to the space
$V$ of monic polynomials of degree $2m$ in $x$ with vanishing
$x^{2m-1}$ term is an isomorphism of $\Delta$
with a neighborhood of the origin in
$V$: given an equation
\begin{equation*}
\begin{split} (x^m + \alpha_{m-2}x^{m-2} + \dots  + \alpha_1x +
\alpha_0)^2 - 4(\beta_{m-1}x^{m-1} + \dots + \beta_1x + \beta_0)  \\
=x^{2m} + c_{2m-2}x^{2m-2} + \dots + c_1x + c_0
\end{split}
\end{equation*} we can write
$$ \a_{m-2} = {c_{2m-2} \over 2}
$$
$$ \a_{m-3} = {c_{2m-3} \over 2}
$$
$$ \a_{m-4} = {c_{2m-4} - \a_{m-2}^2 \over 2} = {4c_{2m-4} - c_{2m-2}^2
\over 8}
$$ and so on, recursively expressing the coefficients $\a_i$ as polynomials
in  the coefficients
$c_{2m-2},
\dots,c_m$. We can then solve for the $\beta_i$ in terms of the remaining
coefficients
$c_{m-1},\dots,c_0$, thus obtaining a polynomial inverse to the map
$\delta$.

Now, since the equation
$f$ above for ${\cal S}$ is  quadratic  in  $y$,  the fibers of ${\cal S} \to
\Delta$ are expressed as double covers of the $x$-line. The discriminant
$\delta$ is a polynomial of degree $2m$ in $x$, so that the general fiber of
${\cal S} \to \Delta$, viewed as a double cover of the $x$-axis, will have
$2m$ branch points near
$P$. To say that any fiber $S_{\a,\b}$ has $m$ nodes is thus tantamount to
saying that
$\delta_{\a,\b}(x)$ has $m$ double roots--- that $\delta_{\a,\b}(x)$
is the square of a polynomial of degree
$m$. The locus of squares being smooth of dimension $m-1$ in $V$, we see
that
$\Delta_m$ is smooth of dimension $m-1$; indeed, it is given simply by the
vanishing
$\b_0 = \dots = \b_{m-1} = 0$.

Similarly, to say that a fiber $S_{\a,\b}$ has $m-1$ nodes amounts to
saying that
$\delta_{\a,\b}(x)$ has $m-1$ double roots, i.e., that it can be written as a
quadratic polynomial in $x$ times the square of a polynomial of degree
$m-1$:
$$
\delta _{\alpha ,\beta} (x) = (x^{m-1} + \l_{m-2}x^{m-2} + \dots + \l_1x +
\l_0)^2(x^2+\mu_1x+\mu_0).
$$
The Lemma is then proved.
\end{pf}

Now we consider the natural map from a suitable analytic neighborhood
$W$ of $[X]$ to $\Delta$.
 To set this up, let
$r_1,\dots, r_k$ be the old nodes of $X$; since
all the  singularities of
$X$ other than $P$ are nodes, this will consist of
$b$ nodes on
$E$ and
$k-b$ nodes lying off
$E$ where $b = (D\cdot E)$. Since
$m-1$ nodes of the general curve of our family tend to $P$, we have
$k=p_a(D)-m+1$.  Now consider, in an analytic  neighborhood of
the point
$[X]
\in |D|$, the locus $W$ of curves passing through the $r_0(D)-1$ assigned
points
$q_1,\dots,q_{r_0(D)-1}$
 and preserving the nodes $r_1,\dots,r_k$ of $X$---that is, such that the
restriction of the family of curves $\{D_\lambda\}_{\lambda \in |D|}$ to
$W$ is equisingular at each point $r_i$ of $X$. Since this is a total of
$r_0(D)-1 + p_a(D)-m+1=r(D)-m$ points
 and they impose
independent conditions on the linear series $|D|$, we see that $W$ is
smooth of dimension $m$ at $[X]$.

We then get a natural map $\phi : $W$ \to \Delta$
such that $\phi ([X])=0$.
We will prove that $\phi $ is an immersion and that the intersection of
$\phi (W)$ with
 $\Delta _{m-1}$
 is the union $\Delta _m$ with a smooth curve $\Psi$; moreover $\Psi$ and
$\Delta _m$
will have contact of order $m$ at the origin.
 This will conclude the proof of Proposition~\ref{describegamma}; in fact
 the original family $\X \to \Gamma$
 will be  the pullback to $W$ of the restriction to
$\Psi$ of the  versal deformation
${\cal S} \to \Delta$.

To illustrate, here is a representation of  the simplest case
$m=2$. This does not convey the general picture, because $\phi(W) \cap
\Delta_{m-1}$ happens to be proper.  Also,  the picture is inaccurate in at
least one respect: the actual surface
$\Delta_1$ in the deformation space of a tacnode is also singular along the
locus of curves
$S_{\a,\b}$ with cusps.

\newpage

\vspace*{7.3in}

\ni \special {picture deformation}

\

\ni (Note that we see again locally the picture that
we have already observed
globally in the linear series $|D|$: the closure of the variety $V(D)$ of
irreducible rational curves has the expected dimension; but the locus of
rational curves has another component of equal or larger dimension.)

\

Here is the outline of our argument.

First in \ref{phi} we will establish that
$\phi$ is an immersion, and identify in part its tangent space.

Second,  in \ref{example} we treat a special case.
We prove the Proposition by direct calculation when $\phi (W)$ is the
linear subspace  given by equations $\beta _1=\ldots =\beta _{m-1} =0$. The
results of \ref{example} also appear in
\cite{R}; we include our proof for the sake of completeness.

Then we use the action of the automorphism group of the singularity $(X,P)$
(cf. Lemma~\ref{transitivity})
to deduce the statement for any smooth, $m$-dimensional subvariety of
$\Delta$ containing $\Delta _m$ whose
tangent plane at the origin is not contained in the hyperplane $\beta_0 =0$ .

The proof of our Proposition (and of Proposition~\ref{describegamma} ) is
then completed in
the remaining part of \ref{end}.

\

\subsubsection{The deformations coming from $V(D)$}
\label{phi}

Let $\phi :W\to \Delta$ be as before,  denote by $H$ the subspace of $\Delta $
given by $\b _0=0$.
Then we have

\begin{lm}
\label{firstsinglemma}
The map $\phi$ is  an immersion; the
tangent space to the image at the origin  contains the plane
$\b_0 = \dots = \b_{m-1} = 0$ but is not contained in $H$.
\end{lm}

\ni \underbar{Remark}. It is important to note here, and throughout the
following argument,
that while the loci $\Delta_m$ and $\Delta_{m-1}$ are well-defined subsets
of the base $\Delta$
of the deformation space of our tacnode,  $H$ is not; it depends on the
choice of coordinates. It  is well-defined, however, as a hyperplane in the
tangent space
$X(\Delta)= \O_{\F_n, P}/{\cal J}$ to $\Delta$ at the origin: it
corresponds  to the
 quotient ${\bf m}/{\cal J} \subset \O_{\F_n, P}/{\cal J}$ of the maximal
ideal ${\bf m}
\subset \O_{\F_n, P}$.

\

\begin{pf} The projective tangent space to $W$ at
the point $[X]$ is  the sublinear series of $|D|$ of  curves passing
through the points $r_1,\dots,r_k$ and
$q_1,\dots,q_{r_0(D)-1}$. The kernel of the differential at $[X]$
of the map $\phi$ is thus the vector space of sections of the line bundle
${\cal L} = \O_{\F_n}(D)$ vanishing at $r_1,\dots,r_k$ and
$q_1,\dots,q_{r_0(D)-1}$ and lying in the subsheaf ${\cal L}
\otimes {\cal J}$, where ${\cal J} \subset \O_{\F_n,P}$ is as before the
Jacobian ideal of $[X]$ at $P$. The zero locus of such a section will be a
curve in the linear series
$|D|$ containing
$r_1,\dots,r_k, q_1,\dots,q_{r_0(D)-1}$  and  $P$ and so must
contain $E$, that is, must be of the form $E+G$ with $G \in |D-E|$. Moreover,
from the description above of
${\cal J}$ we see that $G$ must also have contact of order at least $m$ with
$E$ at $P$ as well as pass through the $k-b$ nodes of
$X$ lying off $E$ and the  assigned points
$p_1,p_2,q_3,\dots,q_{r_0(D)-1}$. This represents a total of
\begin{equation*}
\begin{split} m + r_0(D)-1 + p_a(D)-m+1-b&=r(D)-b \\ &= r(D-E)+1
\end{split}
\end{equation*} conditions, so we need to show that they are
independent to conclude that no such curve exists. But they are also a
subset of the adjoint conditions of $X$, hence impose independent
conditions on the series
$|D+K_{\F_n}| = |D-C-E-2F|$, and  hence on the series $|D-E|$.

The remaining statements of the lemma,
that the tangent space to the image contains the plane $\b_0 = \dots =
\b_{m-1} = 0$ but is not contained in the hyperplane $\b_0 = 0$, follow
from the facts that the image contains the subvariety $\Delta_m$ and that
not every curve in the linear series $|D|$ containing
$r_1,\dots,r_k, q_1,\dots,q_{r_0(D)-1}$  contains  $P$
\end{pf}

 \

\subsubsection{A special case}
\label{example}

Next, having seen that $\phi(W)$ is a smooth, $m$-dimensional variety of
$\Delta$ we will consider the intersection of
$\Delta_{m-1}$ with the simplest possible space satisfying the statement of
the previous Lemma, the plane
$\Lambda$ given by $\b_1 =
\dots =
\b_{m-1} = 0$. We obtain
\begin{lm}\label{linearcase} The intersection of $\Delta_{m-1}$ with
$\Lambda$ consists of the union of $\Delta _m$ with multiplicity
$m$ and a smooth curve $\Psi$ having contact of order $m$ with $\Delta _m$
at the origin.
\end{lm}
\begin{pf} Restricting to  $\Lambda$, we can rewrite the equation
of the family more simply as
$$ y^2+yx^m + \alpha_0y + \alpha_1xy + \dots + \alpha_{m-2}x^{m-2}y +
\beta = 0
$$ and the discriminant as
$$
\delta(x) = (x^m + \alpha_{m-2}x^{m-2} + \dots  + \alpha_1x + \alpha_0)^2
- 4\beta
$$ We need  now  to express the condition that $\delta$ has $m-1$ double
roots. One obviously sufficient condition is that $\b=0$, so that
$\delta$ is a square. If we assume $\b \ne 0$, however, things get more
interesting. To see the locus of $(\a_0,\dots,\a_{m-2},\b)$ that satisfy this
condition,  set
$$
\nu(x) = x^m + \alpha_{m-2}x^{m-2} + \dots  + \alpha_1x + \alpha_0
$$
and write
$$
\delta(x) = \nu(x)^2-4\b = (\nu(x)+2\sqrt{\b})\dot(\nu(x)-2\sqrt{\b}) .
$$
Now, if $\b \ne 0$, the two factors in this last expression have no
common factors; so if their product has $m-1$ double roots, each must
have a number of double roots itself:
$\nu(x)+2\sqrt{\b}$ and $\nu(x)-2\sqrt{\b}$ are polynomials of degree
$m$ with a combined total of $m-1$ double roots. In fact, this  uniquely
characterizes $\nu$ and $\b$ up to a one-parameter group of
automorphisms of $\P^1$, as we will prove in the following.
\begin{lm}
\label{sublemma}
Let $\gamma$ be a nonzero scalar, and let
$m$ be a positive integer. There is a
 polynomial $\nu(x)$ of degree $m$, monic with no $x^{m-1}$ term, such
that

1.  if $m$ is odd, the polynomials $\nu(x)+\gamma$ and $\nu(x)-\gamma$ each have
$(m-1)/2$ double roots; $\nu(x)$ is unique up to replacing $\nu(x)$ by
$-\nu(-x)$;

2.  if $m$ is even,  $\nu(x)+\gamma$ has
$m/2$ double roots and $\nu(x)-\gamma$ has
$(m-2)/2$ double roots; in this case $\nu$ is unique.
\end{lm}
\begin{pf} Suppose that $\nu(x)$ is a polynomial satisfying the conditions
of the lemma. Take first the case of $m = 2\ell+1$ odd, and consider the
map
$\nu :
\P^1
\to
\P^1$ given by  $\nu(x)$. This is a map of degree $m$, sending the point
$\infty$ to
$\infty$, and totally ramified there. In addition the hypotheses assert that
over the points $\pm\gamma$ in the target we have $\ell$ ramification
points. The point is, this accounts for a total of $(m-1) + 2(\ell-1) = 2m-2$
ramification points, and these are all a map of degree $m$ from $\P^1$ to
$\P^1$ will have. We have thus specified the covering
$\nu$ up to a finite number of coverings, and our principal  claim is that in
fact we have described $\nu$ uniquely, up to automorphisms of the
domain.

This is combinatorial. The monodromy permutation $\sigma$ around the
point $\infty$ is cyclic, while the the monodromy permutations $\tau$ and
$\mu$ around $\gamma$ and
$-\gamma$ are each products of $\ell$ disjoint transpositions. Our claim
that there is a unique such covering of $\P^1$ by $\P^1$ amounts then to
the assertion that, up to the action of the symmetric group $\frak{S}_m$
by conjugation, there is a unique pair of permutations $\tau$ and $\mu$,
each a product of $\ell$ disjoint transpositions, whose product $\tau \circ
\mu$ is cyclic of order $m$.

To see this, start with the unique element of the set on which $\tau$ and
$\mu$ act that is fixed by $\tau$ , and label it 1.  This element cannot also
be fixed by $\mu$; give the element exchanged with it by
$\mu$ the label 2, and let the element exchanged with 2 by $\tau$ be
labelled 3. This also cannot be the fixed point of $\mu$, or else the subset
$\{1,2,3\}$ would be fixed by both $\tau$ and $\mu$; let 4 be the element
exchanged with it by $\mu$ and 5 the element exchanged with 4 by
$\tau$. We can continue in this way until we have exhausted all the
elements of the set; and so we see that we can label the elements of the set
$\{1,2,\dots,m\}$ so that
$$
\mu = (1,2)(3,4)\dots(2l-1,2l)
$$ and
$$
\tau = (2,3)(4,5)\dots(2l,2l+1)
$$ This establishes the uniqueness of the covering $\nu$ up to
automorphisms of the domain in case $m$ is odd.

The case of $m = 2\ell$
even is similar; we see that we can always label the sheets of $\nu$ so that
the monodromy permuations  $\tau$ and $\mu$ around
$\gamma$ and
$-\gamma$ have the form
$$
\mu = (1,2)(3,4)\dots(2l-1,2l)
$$ and
$$
\tau = (2,3)(4,5)\dots(2l-2,2l-1)
$$

To complete the proof of Lemma \ref{sublemma}, consider the effect on
$\nu$  of automorphisms of the domain. The requirement that $\nu(\infty)
= \infty$---that is, that
$\nu(x)$ is a polynomial!---restricts us to the group of automorphisms $x
\mapsto ax+b$; the requirement that
$\nu(x)$ have no
$x^{m-1}$ term limits us to automorphisms of the form $x \mapsto ax$;
and the fact that $\nu(x)$ is monic says that $a$ must be an $m$-th root of
unity. Finally, the fact that the ramification points map to $\pm\gamma$
determines $\nu(x)$ completely in case
$m$ is even (where the two branch points $\pm\gamma$ have different
multiplicity), and up to the automorphism $x \mapsto -x$ in case $m$ is
odd.
\end{pf}

\

Back to the proof of  Lemma \ref{linearcase}.  Note first that, by
uniqueness, $\nu(x)$ will be even when
$m$ is even and odd when
$m$ is odd. Note also that if we do not specify the value of $\gamma$ the
polynomial
$\nu(x)$ will not be unique; we can replace it with $u^m\nu(x/u)$ for any
nonzero scalar $u$. Now, suppose first that
$m = 2\ell$ is even. Choose
$\gamma = 1$, and let
$$
\nu(x) = x^m + c_{m-2}x^{m-2} +  c_{m-4}x^{m-4} + \dots +c_0
$$ be the polynomial satisfying the conditions of the lemma. Then any
collection
$(\alpha_0,
\alpha_1, \dots,\alpha_{m-2},\beta)$ with $\b \ne 0$ such that the
discriminant
$$
\delta(x) = (x^m + \alpha_{m-2}x^{m-2} + \dots  + \alpha_1x + \alpha_0)^2
- 4\beta
$$ has $m-1$ double roots must be of the form
\begin{equation*}
\begin{split}
\a_0&= t^\ell \cdot c_0 \\
\a_1 &= 0 \\
\a_2 &= t^{\ell-1} \cdot c_2 \\
\a_3 &= 0 \\
\a_4 &= t^{\ell-2} \cdot c_4
\end{split}
\end{equation*} and so on, ending with
$$
\a_{m-2} = t \cdot c_{m-2};
$$ with finally
$$
\b = {t^m \over 4}.
$$ This is then a parametric representation of the closure $\Psi$ of the
intersection
$\Lambda \cap (\Delta_{m-1} \setminus \Delta_m)$. It is obviously a
curve; the fact that it is smooth is visible from the coordinate $
\a_{m-2} = t \cdot c_{m-2}$; and we see that it has contact of order $m$
with
$H$ from the exponent in the expression for $\b$.

Finally, in case $m = 2\ell+1$ is even we get a similar expression. Let
$$
\nu(x) = x^m + c_{m-2}x^{m-2} +  c_{m-4}x^{m-4} + \dots +c_1x
$$ be the polynomial satisfying the conditions of the lemma for
$\gamma = 1$. Then any collection
$(\alpha_0,
\alpha_1, \dots,\alpha_{m-2},\beta)$ with $\b \ne 0$ such that the
discriminant
$$
\delta(x) = (x^m + \alpha_{m-2}x^{m-2} + \dots  + \alpha_1x + \alpha_0)^2
- 4\beta
$$ has $m-1$ double roots must be of the form
\begin{equation*}
\begin{split}
\a_0 &= 0 \\
\a_1 &= t^\ell \cdot c_1 \\
\a_2 &= 0 \\
\a_3 &= t^{\ell-1} \cdot c_3
\end{split}
\end{equation*} and so on, ending with
$$
\a_{m-2} = t \cdot c_{m-2};
$$ again we have
$$
\b = {t^m \over 4}.
$$ So once more we see that  $\Psi$  is a smooth curve having contact
of order $m$ with
$H$ at the origin.
\end{pf}

\

Let us now prove Proposition
\ref{singtotalspace}  in this special case. First, in the case $m =
2\ell$ even,  the restriction
${\cal S}_\Psi \to
\Psi$ of the family ${\cal S}
\to \Delta$ to $\Psi$ has equation
$$
y^2 + y(x^m + tc_{m-2}x^{m-2} + t^2 c_{m-4}x^{m-4} + \dots +
t^{m/2}c_0) + {t^m
\over 4}=0
$$
We can think of the total space ${\cal S}_\Psi$ of this family as a
double cover of the $(x,t)$-plane, with branch divisor the zero locus of the
discriminant
$$
\delta = (x^m + tc_{m-2}x^{m-2} + t^2 c_{m-4}x^{m-4} + \dots +
t^{m/2}c_0)^2 - t^m
$$
By hypothesis, for each value of $t$ the polynomial $\delta$ is the product
of the square of a polynomial $g_t(x)$ of degree $m-1$ and a quadratic
polynomial $h_t(x)$. Since $\d$ is even, $g^2$ and $h$ must each be; and
given the homogeneity of $\d$ with respect to $t$ and $x$ we see that
we can write
$$
\d = x^2(x^2-\l_1t)^2(x^2-\l_2t)^2 \dots (x^2-\l_{\ell-1}t)^2 \cdot
(x^2-\mu t)
$$
for suitable constants $\l_1,\dots,\l_{\ell-1}$ and $\mu$. For example, in
case
$m=2$, the equation of
${\cal S}_\Psi$ is simply
$$
y^2 + y(x^2 + t) + {t^2
\over 4}=0
$$
and the discriminant is just $\d = x^2(x^2-2t)$. In general, the branch
divisor of ${\cal S}_\Psi$ over the $(x,t)$-plane will be simply a union
of the $t$-axis, with multiplicity 2; $\ell-1$ parabolas tangent to the
$x$-axis at the origin, each with multiplicity 2; and one more parabola
tangent to the
$x$-axis at the origin and appearing with multiplicity 1. The double cover
${\cal
S}_\Psi$ will thus be nodal over the double components of this branch
divisor, and smooth
elsewhere.

Finally, the  normalization ${\cal S}_\Psi^\nu$ of the total space ${\cal
S}_\Psi$ will  be  the double cover of the $(x,t)$-plane
branched over the single component of multiplicity 1 in the branch
divisor; that is, it will have equation
$$
y^2 = x^2-\mu t
$$
and in particular, since the component $(x^2-\mu t)$ is smooth, ${\cal
S}_\Psi^\nu$ will be smooth as well, establishing Proposition
\ref{singtotalspace} for this particular family.

The picture in case $m = 2\ell+1$ is odd is exactly the same: here ${\cal
S}_\Psi$ has equation
$$
y^2 + y(x^m + tc_{m-2}x^{m-2} + t^2 c_{m-4}x^{m-4} + \dots +
t^{m/2}c_1x) - {t^m
\over 4}=0
$$
with discriminant
\begin{equation*}
\begin{split}
\delta &= (x^m + tc_{m-2}x^{m-2} + t^2 c_{m-4}x^{m-4} + \dots +
t^{m/2}c_1x)^2 - t^m \\
 &= (x^2-\l_1t)^2(x^2-\l_2t)^2 \dots (x^2-\l_{\ell}t)^2 \cdot
(x^2-\mu t)
\end{split}
\end{equation*}
for suitable constants $\l_1,\dots,\l_{\ell}$ and $\mu$. For example, in
case
$m=3$, the equation of
${\cal S}_\Psi$ will be
$$
y^2 + y(x^3 -3tx) - t^3=0
$$
(we are scaling $t$ here to make the coefficients nicer), and the
discriminant is just
\begin{equation*}
\begin{split}
\d &= (x^3 -3tx)^2 +4t^3\\
&=x^6 -6tx^4+9t^2x^2+4t^3 \\
&=(x^2+t)^2(x^2+4t)
\end{split}
\end{equation*}
In general, for $m$ odd the branch divisor of ${\cal S}_\Psi$
over the
$(x,t)$-plane will be simply a union of $\ell$
parabolas tangent to the
$x$-axis at the origin, each with multiplicity 2; and one more parabola
tangent to the
$x$-axis at the origin and appearing with multiplicity 1. As before, the
normalization ${\cal S}_\Psi^\nu$ of the total space ${\cal S}_\Psi$
will  be simply the double cover of the $(x,t)$-plane branched over the
single component $(x^2-\mu t)$ of multiplicity 1 in the branch divisor; and
as before, since this component  is smooth, ${\cal
S}_\Psi^\nu$ will be smooth as well, establishing Proposition
\ref{singtotalspace} in this case.
\end{pf}

\subsubsection{The geometry of the locus $\Delta_{m-1}$}
\label{end}
In order to focus on the essential  aspects of the
geometry of
$\Delta_{m-1}$, and in particular to remove the excess intersection of
$\phi(W) \cap
\Delta_{m-1}$,  we will work on the blow-up
$\tau : \tilde\Delta = {\operatorname {Bl}}_{\Delta_m}\Delta \to \Delta$ of
$\Delta$ along $\Delta_m$. To express our results, we have to introduce
some notation. We will denote by
$Z = \tau^{-1}(\Delta_m)$ the exceptional divisor of the blow up, and by
$\tilde\Delta_{m-1}$ and
$\tilde W$ the proper transforms of
$\Delta_{m-1}$  and $\phi(W)$  in
$\tilde\Delta$.

Our  goal will be to describe the intersection $Z_{m-1} : =
\tilde\Delta_{m-1}
\cap Z$.  The fibers of
$Z$ over
$\Delta_m$ are projective spaces
$\P^{m-1}$with homogeneous coordinates $\b_0,\dots,\b_{m-1}$; we will
denote  the fiber $\tau^{-1}(0)$ of
$Z$ over the origin by $\Phi$, by $\Phi_0 \subset \Phi$ the open set given
by  $\b_0
\ne 0$, and by $Q$ the point of $\Phi$ with coordinates $[1,0,\dots,0]$ (this
is the point of intersection of $\tilde W$ with $\Phi$ in the example
above).

Note that there is  a more intrinsic
characterization of the  $\Phi$: the tangent space to $\Delta_m$ at the origin
is the subspace of $\O_{\F_n, P}/{\cal J}$ of polynomials divisible by $y$,
so that
$\Phi$---the projectivization of the normal space---is just the space of
polynomials
in
$x$ modulo those vanishing to order $m$ at $P = (0,0)$ and modulo scalars.
In these terms,
$\Phi_0$ is simply the subspace of polynomials not vanishing at the origin
and $Q$
the point corresponding to constants.

To study  $\Delta_{m-1}$  we use  the action of the
automorphism group of the  deformation space ${\cal S} \to \Delta$. We have
many automorphisms of the germ of the singularity $(X, P)$: for example, for
any power series
$$c(x) = c_1x + c_2x^2 + c_3x^3 +
\dots$$  with $c_1 \ne 0$ we can define an automorphism of the germ by
\begin{equation*}
\begin{split}
\gamma_c : (x,y) &\mapsto  \bigl(c(x), {c(x)^m \over x^m} \cdot y) \\
&=\bigl(c_1x
+ c_2x^2 +\dots, (c_1 + c_2x + c_3x^2 + \dots)^m \cdot y \bigr) \\
&=\bigl(c_1x +
c_2x^2 +\dots, c_1^my + mc_1^{m-1}c_2xy + (mc_1^{m-1}c_3 + {m(m-1)
\over 2}c_1^{m-2}c_2^2)x^2y + \dots \bigr).
\end{split}
\end{equation*}

Let $G$ be the group of automorphisms of the germ $(X, P)$.  By the
naturality
of the versal deformation space, $G$ acts as well on it, that is, $G$ acts
equivariantly on $\cal S$ and $\Delta$. Since the action on
$\Delta$ preserves the subvariety
$\Delta_m$, it lifts to an action on the blow-up $\tilde\Delta$; and since
the action
on
$\Delta$ preserves
$\Delta_{m-1}$ the lifted action will preserve $\tilde\Delta_{m-1}$. We can read
off from the above expression the action of the automorphism $\gamma_c$ on the
tangent space to $\Delta$ at the origin, and thereby on the fiber $\Phi$ of
$Z$ over
the origin: taking as basis for
$X(\Delta) =
\O_{\F_n, P}/{\cal J}$ the monomials $1, x, x^2,\dots,x^{m-1}$ and $y, xy,
x^2y,\dots,x^{m-2}y$, we can express the relevant part of this action as

\begin{equation*}
\begin{split} 1 &\mapsto  1 \\ x &\mapsto  c_1x + c_2x^2 + \dots +
c_{m-1}x^{m-1}
\\ x^2 &\mapsto c_1^2x^2 + 2c_1c_2x^3 + \dots \\ &\vdots \\ x^{m-1} &\mapsto
c_1^{m-1}x^{m-1}
\end{split}
\end{equation*}
 The monomials $x^ky$ are carried into linear combinations of other such
monomials; the exact linear combinations will not concern us. The key fact about
this action, for our present purposes, follows immediately from the description
above:

\begin{lm}\label{transitivity}
 Every orbit of the action of
$G$ on
$\tilde\Delta$ that intersects $\Phi_0$ contains the point $Q$ in its closure.
\end{lm}

We are now prepared to state and prove our main lemma on the geometry of
$Z_{m-1}$ and $\Delta_{m-1}$.

\begin{lm}\label{Zdescription}
1. The fibers of $Z_{m-1}$ over $\Delta_m$ are unions of linear
spaces.

2. For any arc $\a(t)$ in $\Delta_m$ tending to the origin, the limiting
position
of the fiber $Z_{\a(t)}$ of $Z_{m-1}$ over $\a(t)$ is contained in the
complement of
$\Phi _0$.

3.  $\Phi$ itself is contained in (and hence an irreducible component of)
$Z_{m-1}$.
\end{lm}
\begin{pf}
The proof is by induction on $m$,  using Lemma
\ref{linearcase}.

First we  introduce a natural stratification of the locus $\Delta_m$.
Identifying
$\Delta_m$ with the space of monic polynomials of degree $m$ in $x$ with no
$x^{m-1}$ term,
we  look at the loci of polynomials with roots of given multiplicity: for
any partition $m =
m_1 + m_2 + \dots + m_k$ we define the locus $
\Delta\{m_1,\dots,m_k\} \subset
\Delta$ by
\begin{equation*}
\begin{split}
\Delta\{m_1&,\dots,m_k\}  : = \bigl\{(\a_0,\dots,\a_{m-2},0,\dots,0) : \\ &x^m +
\a_{m-2}x^{m-2} + \dots + \a_0
 =  (x-\l_1)^{m_1}(x-\l_2)^{m_2}\dots(x-\l_k)^{m_k} \\
 &\qquad \qquad \qquad \text{for some distinct}
\;
\l_1,\dots,\l_k \bigr\}
\end{split}
\end{equation*}  Note that the codimension of
$\Delta\{m_1,\dots,m_k\}$ in
$\Delta_m$ is $\sum (m_\a-1)$.

 Suppose $\a$ is any point of $\Delta_m$
 other than the origin. Say $\a$ lies in the
stratum
$\Delta\{m_1,\dots,m_k\}$, and write the corresponding polynomial as
$$
(x-\l_1)^{m_1}(x-\l_2)^{m_2}\dots(x-\l_k)^{m_k}
$$
with $\l_1,\dots,\l_k$ distinct. The fiber $S_\a$ of ${\cal S} \to \Delta$ over
$\a$ is a reducible curve consisting of two branches, the $x$-axis $(y=0)$ and
the curve
$y = (x-\l_1)^{m_1}(x-\l_2)^{m_2}\dots(x-\l_k)^{m_k}$, which meet at the
$k$ points
$r_1=(\l_1,0),\dots,r_k=(\l_k,0)$ with multiplicities $m_1,\dots,m_k$.

Let $\Delta(i)$ be
 the versal deformation spaces $ \Delta(S_\a, r_i)$ of the singular
points $r_i
\in S_\a$.
By the
openness of versality the natural map $\sigma$ from a neighborhood $U$ of $\a$
in
$\Delta$ to the product $\prod
\Delta(i)$  has surjective differential at $\a$ (the fibers are the equisingular
deformations of $S_\a$, in which only the locations of the points
$r_i$ on the $x$-axis vary). Let $\Delta_{m_i-1}$ and $\Delta_{m_i} \subset
\Delta(i)$ be the loci in $\Delta(i)$ analogous to $\Delta_{m-1}$ and
$\Delta_{m}$ in $\Delta$, that is, the closures of the loci of deformations
of the
singular points $r_i \in S_\a$ with $m_i-1$ and $m_i$ nodes near $r_i$
respectively. Then in the neighborhood $U$ of $\a$, we have
$$
\Delta_m = \sigma^{-1}\left(\Delta_{m_1} \times \Delta_{m_2} \times \dots \times
\Delta_{m_k}\right)
$$
and
$$
\Delta_{m-1} = \bigcup_{i=1}^k\sigma^{-1}\left(\Delta_{m_1} \times \dots \times
\Delta_{m_i-1}
\times
\dots
\times
\Delta_{m_k}\right)
$$
In other words, the locus $\Delta_{m-1}$ will have $k$ branches
in a neighborhood of
$\a$, each containing $\Delta_m$, along the $i$th of which the fibers of
${\cal S}
\to
\Delta$ will have $m_j$ nodes tending to $r_j$ for each $j \ne i$ and $m_i-1$
nodes tending to $r_i$.

We can use this description to give a more intrinsic
characterization of the fiber $Z_\a = \tau^{-1}(\a)$ of $Z$ over the point $\a$,
analogous to the one given above for $\Phi$. Briefly, $Z_\a$ is the
projectivization
of the normal space to $\Delta_m$ in $\Delta$ at $\a$, which is the product
of the
normal spaces to the $\Delta_{m_i}$ in $\Delta(i)$ at the origin; this is
just the
space of polynomials on the $x$-axis modulo those vanishing to order $m_i$ at
$r_i$ for each $i$.

We may now
apply the induction hypothesis to describe, in these terms, the fiber of
$Z_{m-1}$over $\a$.
By the statement of the Lemma for $m = m_i$, the proper transform of the
$i$th branch of
$\Delta_{m-1}$ will intersect
$Z_\a$ in the linear subspace of $Z_\a$ corresponding to polynomials
vanishing to
order $m_j$ at
$r_j$ for each $j \ne i$; the intersection with $Z_\a$ with the proper
transform of
$\Delta_{m-1}$ itself will be  the union of these linear subspaces.

This
establishes part (1) of the Lemma. Now say that $\a(t)$ is any arc  in
$\Delta_m$
tending to the origin; $\a(t)$ will lie in some stratum
$\Delta\{m_1,\dots,m_k\}$ for all small $t \ne 0$.
As $t$ goes to zero, the singular points $r_i(t)$ of  $S_{\a(t)}$ approach
the point
$P$, so that the limiting position of the intersection with $Z_{\a(t)}$ of
the proper
transform of the
$i$th branch of
$\Delta_{m-1}$ will be simply the linear space of polynomials whose
restriction to the $x$-axis vanishes to order
$m-m_i$ at $P$; in particular, it is contained in the hyperplane $(\b_0 = 0)
\subset \Phi$ of polynomials vanishing at $P$. We have thus proved parts (1) and
(2) of the Lemma, given part (3) for all $m_i < m$.

\

Finally,  we need to prove for each new value of $m$ that $\Phi$ is
contained in (and
hence an irreducible component of) $Z_{m-1}$. Now,  by Lemma \ref{linearcase},
the point $Q = [1,0,\dots,0] \in \Phi$ lies in $Z_{m-1}$. But we have completely
described the closure in $Z_{m-1}$ of the inverse image $\tau^{-1}(\Delta_m
\setminus \{0\})$ of the complement of the origin, and $Q$ is not on it.
$Q$ must
thus lie on an irreducible component of $Z_{m-1}$ not meeting
$\tau^{-1}(\Delta_m \setminus \{0\})$, that is to say, an irreducible
component of
$Z_{m-1}$ contained in $\Phi$; since $Z_{m-1}$ has pure dimension $m-1$, this
irreducible component must be $\Phi$ itself.
\end{pf}

\

For example, here is a picture of $Z_1$ in the case $m=2$. In this case
$Z_1$ has only two components, $\Phi$ and a component finite of degree
$2$ over $\Delta_2$.

\vspace*{3.5in}

\ni \hskip.4in \special {picture blownupdeformation}

\

Next, we  deduce:

\begin{lm}
\label{multiplicitym}
1. $\tilde\Delta_{m-1}$ is smooth everywhere along $\Phi_0$

\ni
2.  The intersection multiplicity of
$\tilde\Delta_{m-1}$ and
$ Z$ along
$\Phi$ is
$m$.
\end{lm}
\begin{pf} We use the analysis carried out in Lemma \ref{linearcase}. Let
$\tilde\Lambda$ be the proper transform of the linear space
$\Lambda$ in $\tilde\Delta$. Since no
 component of $Z_{m-1}$ other than $\Phi$ passes through $Q$, the only
component of the intersection $\tilde\Lambda \cap \tilde\Delta_{m-1}$ containing
$Q$ will be the proper transform $\tilde\Psi$ of the curve $\Psi \subset
\Delta$ described in Lemma
\ref{linearcase}. Since this is smooth, and the intersection $\tilde\Lambda \cap
\tilde\Delta_{m-1}$ is proper in a neighborhood of $Q$ ($\tilde\Lambda$ and $
\tilde\Delta_{m-1}$ each have dimension
$m$ in the $(2m-1)$-dimensional $\tilde\Delta$, and their intersection is
locally a
curve) it follows that $
\tilde\Delta_{m-1}$ must be smooth at $Q$. By Lemma
\ref{transitivity}, then, it must be smooth at every point of $\Phi_0$.

For the second statement, notice that
 Lemma \ref{linearcase} asserts that this is true when restricted to the
proper transform $\tilde\Lambda$, and it follows that it is true on
$\tilde\Delta_{m-1}$
\end{pf}

\ni
{ \it End of the proof of Proposition~\ref{singtotalspace}}.
We shall now conclude that the intersection of
$\phi (W)$ with $\Delta _{m-1}$ is the union of
$\Delta _m$ and a smooth curve $\Psi$, such that $\Psi $ has contact of
order $m$
with $\Delta _m$ at the origin.  Notice that this will conclude the proof of
Proposition~\ref{describegamma} as well.
We know
 from Lemma \ref{firstsinglemma} that $\phi(W)$ is smooth, so that its
proper transform
$\tilde W$ intersects
$Z$ in a section, crossing $\Phi$ at some point $R$; we likewise have from Lemma
\ref{firstsinglemma} that $R \in \Phi_0$. $ \tilde\Delta_{m-1}$ is then
smooth at
$R$. Since the tangent space to $ \tilde\Delta_{m-1}$ at $R$ contains the
tangent
space to $\Phi$ and the tangent space to $\tilde W$ at $R$ is complementary to
the tangent space to $\Phi$, $ \tilde\Delta_{m-1}$ and $\tilde W$ intersect
transversely in a smooth curve in a neighborhood of $R$; since  that curve
is not
tangent to $\Phi$ at $R$, its image $\Phi \subset \Delta_{m-1} \cap \phi(W)$
is again a smooth curve. Finally, the intersection number of $\Psi$ with
$\Delta_m$ in
$\phi(W)$ will be the intersection number of $ \tilde\Delta_{m-1}$, $\tilde W$
and $Z$ at $R$; which by Lemma \ref{multiplicitym} will be $m$. We have thus
completed the proof of Proposition \ref{singtotalspace}.

Now, the inverse image of $\Psi$ in $W$ is an analytic neighborhood of $\Gamma$;
therefore to conclude the proof we need to show that the total space $\cal S
_{\Psi }$
of the versal deformation
over $\Psi$ is smooth at the point corresponding to $p$.
This follows as in the end of \ref{example}.

\section{Formulas}

Before we prove our  formulas, we need a simple result on the  order of
zeroes and poles of the
cross-ratio function $\phi$.

\subsection{A remark on the cross-ratio function} Suppose we are given a
family $ f : \X \to  B$
over a smooth one-dimensional base $B$, whose restriction $\tilde f :
\tilde {\X }=
f^{-1}(\tilde B) \to \tilde B$ to the complement
$\tilde B = B
\setminus
\{b_0\}$ of a point $b_0 \in B$ is a family of smooth rational curves; and
four sections
$p_i :
\tilde B
\to \X$, disjoint over $\tilde B$. We get a  map
$
\tilde \phi :
\tilde B
\to
\mg$, which then extends over $B$; and the problem is to determine the
coefficient of the point
$b$ in the pullback via $\tilde \phi$ of the boundary components of $\mg$. To
put it another way, the cross-ratio  of the four sections $p_1,p_3,p_2,p_4$
defines a rational function on
$\tilde B$ and hence on $B$; and we ask simply for the order of zero or
pole of this function at $b_0$.

We will answer this in terms of any completion of our family to a family of
nodal rational curves. Recall first of all the set-up of
section~\ref{setup}: we have a resolution
of singularities $\Y \to
B$ of the total space of our family, such that $\Y \to B$ is a family of nodal
curves and the extensions of the sections $p_i$ to $\Y$  are disjoint. We
then proceed to blow
down ``extraneous" components of
$Y$ to arrive at  the minimal smooth
semistable model of our family: that is, a family
${\cal Z}
\to B$ such that ${\cal Z}$ is smooth, the fibers $Z_b$ are nodal, the sections
$p_i$ are disjoint and ${\cal Z} \to B$ is minimal with respect to these
properties.  Finally, we  blow down the intermediate
components
 in this chain to arrive at a family ${\cal W}
\to B$ of 4-pointed stable curves. The special fiber $W$ of this family
will have just two components (or one, if $\ell = 0$), with a singularity of
type $A_\ell$ at the point of their intersection.

In these terms we  prove

\begin{lm}\label{crossratiomult} If the sections $p_1$ and $p_2$
(respectively, $p_1$ and $p_3$) meet
the same component of
$Y$, then the point $b_0$ is a zero (respectively, pole) of multiplicity
$\ell$ of
the function $\phi$.
\end{lm}
\begin{pf} We will consider the case where $p_1$ and $p_2$ meet the
same component of
$W$. Note first that if we blow down the component of $W$ meeting
$p_1$ and
$p_2$, we arrive at a smooth family, that is (replacing $B$ if necessary by a
neighborhood of
$b_0$ in
$B$), a product $B \times \P^1$. (Equivalently, we could arrive at this
family by blowing down the component of $Z$ meeting $p_1$ and
$p_2$, then doing the same thing on the resulting surface, and so on $\ell$
times.) $p_3$ and
$p_4$ will remain disjoint from each other in this process, and disjoint
from
$p_1$ and
$p_2$; but $p_1$ and
$p_2$ will meet each other with contact of order $\ell$: in other words,
we can choose an affine coordinate $z$ on $\P^1$ and a local coordinate $t$
on $B$ centered around $b_0$ so that the sections $p_i$ are given by
$$
p_1(t)  = t^\ell; \quad p_2(t) \equiv 0; \quad p_3(t) \equiv 1; \quad {\rm
and} \quad p_4(t) \equiv \infty .
$$ The cross-ratio function is then $\phi(t) = 1-t^\ell$, which takes on the
value 0 with multiplicity $\ell$ at $t=1$
\end{pf}

\

\subsection{The recursion for $\F_2$}

Let $D$ be any effective divisor class other than $E$ on the ruled surface
$S = \F_2$.
We are going to find a formula for the degree $N(D)$ of the variety $V(D)
\subset |D|$. To set
this up, we start by  choosing as usual $r_0(D) - 1$ general points on $S$,
which we label
$p_1, p_2, q_3,\dots,q_{r_0(D) - 1}$, and consider the one-parameter family
$\X \to
\Gamma$ of curves
$X
\in V(D) \subset |D|$ passing through $\{p_1, p_2, q_3,\dots,q_{r_0(D) - 1}\}$.
As before, we let $ \Gamma^\nu$ be the normalization of $\Gamma$ and $\X^\nu \to
\Gamma^\nu$ the normalization of the pullback family. Next, we fix general
curves
$C_3$ and
$C_4
\in |C|$ in the linear series $|C|$, and adopt as usual the convention that we
will choose points $p_3$ and $p_4$ on the curves $X$ of our family lying
on $C_3$ and $C_4$ respectively.  Making the corresponding base change, we
arrive
at a family $\X \to B$; as before, we will denote  by
$\Y$ the minimal desingularization of $\X$ and by ${\cal Z} \to B$ the
smooth semistable model.

Then we calculate the degree of the cross-ratio map
$\phi : B
\to \mg \cong \P^1$ in two ways by  equating the number of zeroes and poles of
$\phi$. We get one contribution to the degree of
$\phi^*(0)$ immediately from the curves $X$ in our family that happen to
pass through either of the two points of intersection of $C_3$ with $C_4$; this
gives a total contribution of $2 \cdot N(D)$ to the degree of $\phi^*(0)$.

The remaining zeroes and poles of
$\phi$ necessarily
correspond to reducible curves in the family $\{X\}$. There are
two types of these: those that contain $E$ and
those that don't.
Consider first  a reducible curve $X$ in our family that does not contain
$E$. By
Proposition \ref{codimensionone}, this must be  of the form $X = X_1 + X_2$
where $X_i$
is a general member of the family $V(D_i)$ for some pair of divisor classes
$D_1$
and
$D_2$ adding up to $D$. In particular, $X_i$ is an irreducible rational
curve with
$p_a(D_i)$ nodes, and $X_1$ and $X_2$ intersect transversely in $(D_1 \cdot
D_2)$
points.  Moreover, by Proposition \ref{describegamma}, the curve $\Gamma$
will  consist
of $(D_1 \cdot D_2)$ smooth branches near the point  $[X]$, corresponding
to the
points of intersection of $X_1$ and $X_2$; thus there are $(D_1
\cdot D_2)$ points in the normalization $\Gamma^\nu$ lying over each such
point $[X]
\in \Gamma$.

How does such a fiber of the family $\X \to B$ contribute to the degrees of
either
$\phi^*(0)$ or
$\phi^*(\infty)$? It depends on how the points $p_i$ are distributed. If
three or four
lie on one component, it does not contribute to either, but if there are
two on each it
may: for example, if
$p_1$ and
$p_2$ lie on the same component---say $X_1$---of $X$,  and $p_3$ and
$p_4$ on the other, we get a zero of $\phi$. Now, as we observed in the proof of
Proposition \ref{codimensionone}, each component $X_i$ of $X$ must contain
exactly
$r_0(D_i)$ of the points $p_1, p_2, q_3,\dots,q_{r_0(D) - 1}$. If $X_1$ is
to contain
$p_1$ and
$p_2$, it will contain $r_0(D_1)-2$ of the points $q_\a$, and $X_2$ will
contain the
remaining $r_0(D)-r_0(D_1)+1 = r_0(D_2)$. Thus, to specify such a fiber, we have
first to break the $r_0(D)-3$ points $q_\a$ into disjoint sets of
$r_0(D_1)-2$ and
$r_0(D_2)$. The curve $X_1$ can then be any of the $N(D_1)$ irreducible rational
curves in the linear series $|D_1|$ passing through $p_1$, $p_2$ and the
first set,
while $X_2$ can then be any of the $N(D_2)$ irreducible rational
curves in the linear series $|D_2|$ passing through the second set.
Altogether, then,
we see that there will be
$$N(D_1)N(D_2) \binom{r_0(D)-3}{r_0(D_1)-2}$$
points in $\Gamma$ of this type, and correspondingly
 $$N(D_1)N(D_2) (D_1 \cdot D_2) \binom{r_0(D)-3}{r_0(D_1)-2}$$
such points in the
normalization $\Gamma^\nu$. Finally, if a fiber of $\X \to B$ lying
over such a point of
$\Gamma^\nu$ is to contribute to $\phi^*(0)$, we have to choose
$p_3$ and $p_4$ to lie on $X_2$, that is, to be any of the $(D_2 \cdot C)$
points of
intersection of $X_2$ with $C_3$ and $C_4$ respectively. There are thus a
total of
$(D_2 \cdot C)^2$
 fibers of $\X \to B$ of this type lying
over each such point of
$\Gamma^\nu$.

To complete the calculation of the contribution of fibers of this type to
the degree of
$\phi^*(0)$, we observe that the fiber of the normalization $\X^\nu$ over such a
point will have two components, the normalizations of the curves $X_i$,
meeting at
one point (the point of each lying over the new node). Moreover, by
Proposition\ref{singtotalspace}, the total space $\X^\nu$ will be smooth at
such a
point; and it follows by Lemma \ref{crossratiomult} that the
corresponding point of
$B$ will be a simple zero of $\phi$. In sum, then, fibers of $\X \to B$ of
this type
contribute a total of
$$N(D_1)N(D_2) (D_1 \cdot D_2)
\binom{r_0(D)-3}{r_0(D_1)-2}(D_2 \cdot C)^2$$
to the degree of $\phi^*(0)$.

The contribution of such fibers to the degree of the divisor
$\phi^*(\infty)$ is found
analogously, the only difference being that, in order to get a pole of the
cross-ratio,
the points $p_1$ and $p_3$ must lie on one component---say $X_1$---of $X$, while
$p_2$ and $p_4$ will lie on the other. Thus, instead of breaking the $r_0(D)-3$
points $q_\a$ into subsets of $r_0(D_1)-2$ and $r_0(D_2)$, we divide them into
subsets of $r_0(D_1)-1$ and $r_0(D_2)-1$; and instead of $N(D_1)N(D_2)
\binom{r_0(D)-3}{r_0(D_1)-2}$ such points in $\Gamma $ of this
type we have $N(D_1)N(D_2) \binom{r_0(D)-3}{r_0(D_1)-1}$. Similarly, instead of
choosing
$p_3$ among the $(D_2 \cdot C)$ points of $X_2 \cap C_3$, we choose it among the
$(D_1
\cdot C)$ points of $X_1 \cap C_3$; so that instead of $(D_2 \cdot C)^2$
 zeroes of the cross-ratio lying
over each such point of
$\Gamma^\nu$ there will be $(D_1 \cdot C)(D_2 \cdot C)$. Again,
each pole of the cross-ratio corresponding to a fiber of this type will have
multiplicity one; so the total contribution to the degree of $\phi^*(\infty)$ is
$$N(D_1)N(D_2) (D_1 \cdot D_2)
\binom{r_0(D)-3}{r_0(D_1)-1}(D_1 \cdot C)(D_2 \cdot C)$$

It remains to add up the number of zeroes and poles of $\phi$
coming from members of our family containing $E$. Proposition
\ref{codimensionone}
describes all such curves, and the description is particularly simple,
given that we
are on the surface $\F_2$. There are only two types: a degenerate member
$X$ of our
family must consist either of

\

\ni 1.  the union of $E$ and an irreducible rational nodal curve $X_1 \in
|D-E|$,
simply tangent at one point (which will be a smooth point of $X_1$) and meeting
transversely elsewhere;

\vspace*{1.7in}

\hskip.5in \special {picture F2picture1}

\ni or

\

\ni 2.  the union of $E$ and two curves $X_i \in |D_i|$, which will
correspond to
general points of the varieties $V(D_i)$ for some pair of divisor classes
$D_1$ and
$D_2$ with $D_1 + D_2 = D-E$. In particular, $X_1$ and $X_2$ will intersect each
other and $E$ transversely.

\vspace*{2.5in}

\special {picture F2picture2}

\

Now,  we can forget about curves of the first type;
in fact, since $E$ cannot contain any
of the points
$p_1,\dots,p_4$, these will be distinct points of $X_1$. Hence the
cross-ratio function will not be
zero or infinite at such a point of $B$.
On the other hand,
fibers of the second type may contribute. To see what our family looks like in a
neighborhood of such a curve, recall first that by Proposition
\ref{describegamma}, as we approach $X$ along any branch of $\Gamma$, all the
points of intersection of $X_1$ and $X_2$, as well as all but one of the
points of
intersection of each curve $X_i$ with $E$, will be old nodes; exactly one
of the points of
intersection of each $X_i$ with $E$ will be new. The fiber of the
normalized family $\X^\nu \to
\Gamma^\nu$ will thus consist of the normalizations of $X_1$ and $X_2$,
each meeting a copy of
$E$ in one point and disjoint from each other:

\vspace*{3.5in}

\hskip.5in \special {picture F2semistable}

\ni Recall also that the total space of $\X^\nu$ will be smooth along such
a fiber.

Again,
$E$ can't contain any of the points
$p_i$, and if three or four lie on either curve $X_i$ the corresponding point of
$B$ will be neither a zero or a pole of $\phi$; but we may get a contribution if
two are on each $X_i$. Specifically, if $p_1$ and $p_2$ lie on one
component---say
$X_1$---and $p_3$ and $p_4$ on the other, we get a zero of $\phi$; while if
$p_1$ and $p_3$ lie on a component---again, call this one
$X_1$---and $p_2$ and $p_4$ on the other, we get a pole of $\phi$. That said,
we can count the number of such fibers exactly as in the preceding case.

We do the zeroes first. We begin by specifying  a point $[X]$ in
$\Gamma$---that is,
we break the points
$q_\a$ into subsets of size $r_0(D_1)-2$ and
$r_0(D_2)$ respectively, and  choose $X_1$ among the $N(D_1)$ irreducible
rational curves in $|D_1|$ through $p_1$, $p_2$ and the first set and $X_2$ among
the $N(D_2)$ irreducible rational curves in $|D_2|$ through  the
second set. Next, a point in $\Gamma^\nu$: we can take any of the $(D_1
\cdot E)(D_2
\cdot E)$ points of $\Gamma^\nu$ lying over $[X] \in \Gamma$. Lastly, we have to
choose
$p_3$ and $p_4$ among the $(D_2 \cdot C)$ points of intersection of $X_2$ with
$C_3$ and $C_4$ respectively. We have, in sum,
$$
N(D_1)N(D_2) \binom{r_0(D)-3}{r_0(D_1)-2}(D_1 \cdot E)(D_2
\cdot E)(D_2 \cdot C)^2
$$
zeroes of $\phi$ of this type.

The poles of the cross-ratio coming from such  are counted in the
same way; the differences being exactly as in the preceding case: in specifying
the point
$[X]
\in
\Gamma$ we have to choose a subset of $r_0(D_1)-1$ rather than $r_0(D_1)-2$
of the
points $q_\a$; and $p_3$ must be chosen among the $(D_1 \cdot C)$ points of $X_1
\cap C_3$. There are thus a total of
$$
N(D_1)N(D_2) \binom{r_0(D)-3}{r_0(D_1)-1}(D_1 \cdot E)(D_2
\cdot E)(D_1 \cdot C)(D_2 \cdot C)
$$
poles of this type.

There is one important difference between this case and the previous, however:
here, the fiber of the normalization $\X^\nu \to \Gamma^\nu$ has three
components,
with the components $X_1$ and $X_2$ containing the points $p_i$ separated by the
component $E$. Since by Proposition \ref{singtotalspace} the total space
$\X^\nu$ is
smooth, we see by Lemma \ref{crossratiomult} that such points will be double
zeroes and poles of $\phi$. The contribution to the degrees of these
divisors coming
from fibers of this type is thus twice the number of such fibers.

We can now  calculate the  degree of the divisors  $\phi^*(0)$ and
$\phi^*(\infty)$. We have
\begin{equation*}
\begin{split}
\deg(&\phi^*(0)) \\
&= 2 \cdot N(D) \\
& \quad + \sum_{D_1+D_2=D \atop D_1, D_2 \ne E} N(D_1)N(D_2)
\binom{r_0(D)-3}{r_0(D_1)-2}(D_1
\cdot D_2) (D_2 \cdot C)^2 \\
& \quad + 2 \cdot \sum_{D_1+D_2=D-E \atop D_1, D_2 \ne E} N(D_1)N(D_2)
\binom{r_0(D)-3}{r_0(D_1)-2}(D_1 \cdot E)(D_2
\cdot E)(D_2 \cdot C)^2
\end{split}
\end{equation*}
Similarly,
\begin{equation*}
\begin{split}
&\deg(\phi^*(\infty)) \\
& \quad = \sum_{D_1+D_2=D \atop D_1, D_2 \ne E} N(D_1)N(D_2)
\binom{r_0(D)-3}{r_0(D_1)-1}(D_1
\cdot D_2)(D_1 \cdot C)(D_2 \cdot C) \\
& \quad \quad + 2 \cdot \sum_{D_1+D_2=D-E \atop D_1, D_2 \ne E} N(D_1)N(D_2)
\binom{r_0(D)-3}{r_0(D_1)-1}(D_1 \cdot E)(D_2
\cdot E)(D_1 \cdot C)(D_2 \cdot C)
\end{split}
\end{equation*}
To express the final result we introduce    the
notation:
 \begin{equation*}
 \begin{split}
  \gamma (D_1,D_2):=N(D_1)
N(D_2) &\biggl[ {r_0(D)-3 \choose r_0(D_1) - 1} (D_1\cdot C)(D_2\cdot C) \\
& \quad -
{r_0(D)-3 \choose r_0(D_1) - 2} (D_2\cdot C)^2
 \biggr]
 \end{split}
 \end{equation*}
We now write $\deg(\phi^*(0))=\deg(\phi^*(\infty )))$ and solve the
resulting equation
for $N(D)$ to arrive at the recursion formula for
$N(D)$ on
$\F_2$:

\begin{thm}
\label{F2}
Let $D\in \Pic(\F _2)$ and let $\ND$ be the number of irreducible rational
curves in the linear series $|D|$ that pass through $\rD$ general points of
$\F _2$; then we have
\begin{equation*}
\begin{split}
N(D) \quad = \qquad &\frac{1}{2}
\sum _{D_1+D_2=D \atop D_1, D_2 \ne E} \gamma (D_1,D_2) (D_1 \cdot D_2) \\
& \quad  +\sum _{D_1+D_2=D-E \atop D_1, D_2 \ne E} \gamma (D_1,D_2) (D_1 \cdot
E)(D_2\cdot E).
\end{split}
\end{equation*}

\end{thm}

\subsection{The class $2C$ on $\F_n$ }

We will now analyze the linear series $|2C|$ on the ruled surface $\F_n$ for
any $n$.  By restricting ourselves to this linear series we
will arrive at a closed-form expression for $N(D)$ rather than a
recursion. This is clear: since every linear series $|D|$ on $\F_n$ with $D
< 2C$ that
actually contains irreducible curves has arithmetic genus 0, we can say
immediately how many degenerate fibers of each type there are in our
one-parameter family of curves in $|2C|$.

The dimension of the
linear series $|2C|$ is $ 3n+2$. The arithmetic genus of the curves in
the series is $n-1$, so that the expected dimension of the Severi variety is
$r_0(2C) = 2n+3$. This is in fact the actual dimension: any irreducible nodal
curve $D \in |2C|$ will be disjoint from $E$ (if it met $E$, it would contain
it, having intersection number 0 with it); so that the nodes of $D$ will
impose independent conditions on $|2C|$.

So, we choose as usual $2n+2$ general points on $\F_n$, which we label
$p_1, p_2, q_3,\dots,q_{2n+2}$ and consider the one-parameter family of
curves
$X
\in |2C|$ passing through $\{p_1, p_2, q_3,\dots,q_{2n+2}\}$; we will
denote this
family $\X \to \Gamma$. As before, we let $ \Gamma^\nu$ be the normalization of
$\Gamma$ and $\X^\nu \to
\Gamma^\nu$ the normalization of the pullback family. Next, we fix general
curves
$C_3$ and
$C_4
\in |C|$ in the linear series $|C|$, and adopt the convention that we
will choose points $p_3$ and $p_4$ on the curves $X$ of our family lying
on $C_3$ and $C_4$ respectively.  Making the corresponding base change, we
arrive
at a family $\X \to B$; as before, we will denote  by
$\Y$ the minimal desingularization of $\X$ and by ${\cal Z} \to B$ the
smooth semistable model.

Now we consider the cross-ratio map
$\phi : B
\to \mg \cong \P^1$ as before; we shall obtain a formula for $N(D)$ from
$$
\deg \phi ^*(0) = \deg \phi ^*(\infty).
$$
Of course, we get one
contribution to the degree of
$\phi^*(0)$  from the curves  in our family that
pass through any of the $n$ points of intersection of $C_3$ with $C_4$;
this gives a
total contribution of $n \cdot N(2C)$ to the degree of $\phi^*(0)$.

The remaining zeroes and poles of $\phi$
correspond to reducible curves in the family $\{X\}$. As before we look
first at curves that
do not contain $E$.
They can only be   of the form $X = D_1 +
D_2$ where $D_1$ and $D_2$ are each linearly equivalent to $C$.

Such a fiber of the family $\X \to B$ can be  either a pole or a zero of $\phi$,
depending of course on how the points
$p_i$ are distributed. Namely, if
$p_1$ and
$p_2$ lie on the same component $D_i$ and $p_3$ and
$p_4$ on the other, we get a zero. To specify such a fiber, we
simply have to break the $2n$ points $q_\a$ into disjoint sets of $n-1$ and
$n+1$. The two components
$D_i$ of the curve $X$ will be the (unique) curve in the series $|C|$
containing
$p_1$, $p_2$ and the first subset; and the unique curve in the series $|C|$
containing the second subset.

\vspace*{.9in}
\ni \hskip.3in \special {picture FnC+C}

\

Next, we have to count the number of points of $B$ lying over each point of
$\Gamma$ corresponding to such a curve.
Since the general curve of our family
has
$n-1$ nodes, and the
$n$ nodes of
$X = D_1
\cup D_2$ impose independent conditions on the series $|2C|$, the curve
$\Gamma$ will have $n$ (smooth) branches at the point
$[X]$; thus the normalization of $\Gamma$ will have $n$
points lying over $[X]$ (cf. Proposition~\ref{describegamma}).
Moreover, for each of these points there will be a point of
$B$ for every possible choice of points $p_3$ and $p_4$;  these can be any
of the
$(C
\cdot C) = n$ points of intersection of the component $D_i$ not containing
$p_1$ or
$p_2$ with
$C_3$ and $C_4$ respectively. There are thus a total of ${2n \choose n-1}
\cdot n
\cdot n^2$ fibers of $\X \to \Gamma$ of this type.

Now, for each such fiber of $\X \to \Gamma$, the fiber of $\X
\to B$ will be simply the normalization of $X$ at the $n-1$  old nodes:

\newpage

\

\vspace*{.9in}

\hskip.3in \special {picture Fnfinal0}

In particular, it has just two components and is stable, and as we have seen
$\X^\nu$ already is smooth at the node of such a fiber. Each such point  is
thus a simple zero of $\phi$; so the total
contribution to the degree of $\phi^*(0)$ of such curves is

$$ {2n \choose n-1}
\cdot n^3
$$

Similarly, we could have
$p_1$ and
$p_3$  on the same component $D_i$ and $p_2$ and
$p_4$ on the other; in this case, we get a point of $\phi^*(\infty )$. The only
difference in this case is that to specify such a fiber, we  have to break the
$2n$ points
$\{q_\a\}$ into two disjoint sets of
$n$ points apiece. The two components
$D_i$ of the curve $X$ will be the (unique) curves in the series $|C|$
containing
$p_1$ and the first subset; and the unique curve in the series $|C|$
containing $p_2$ and the second subset.  The rest of the analysis is exactly
the same---for each such curve, we get $n^3$ points of $B$, each of which
is a simple pole of the cross-ratio function $\phi$---so the total
contribution to the degree of $\phi^*(\infty)$ of such curves is

$$ {2n \choose n}
\cdot n^3
$$

\

The remainder of the calculation will be spent evaluating the contributions
to the degree of the pullbacks of the boundary components of $\mg$
coming from the curves in our original family containing $E$. As we have
indicated, these curves are of $n-1$ types: for each $k = 1,\dots,n-1$ we
will have a finite number of curves in our family consisting of the sum of
$E$,
$k$ fibers $F_1,\dots,F_k$ of $\F_n \to \P^1$ and an irreducible curve $D$
linearly equivalent to
$C+(n-k)F$, with $D$ having a single point of $(n-k)$-fold intersection with
$E$:

\newpage

\

\vspace*{2.3in}

\ni \special {picture Fn}

\

For each of these types, there are a number of possibilities for the
distribution of the points $p_1,\dots,p_4$ on the various components. For
each such distribution corresponding to points $b$ in the inverse image of
a boundary component
$\Delta$ of
$\mg$, we will consider the number of fibers $X_b$ of that type and the
coefficient with which the corresponding points $b \in B$ appear in the
divisor
$\tilde \phi^*(\Delta)$; in the end we will sum up the contributions to
arrive at an
expression for $N(2C)$ on $\F_n$ .

\

\

\ni $\bullet_1$  \thinspace $p_1, p_2 \in D$; $p_3, p_4 \in F_i$.  \enspace
In such a curve, the fiber components
$F_j$ must each contain one of the points $q_\alpha$. To specify such a
curve, then, we must first choose a subset of $k$ of the $2n$ points $q_\a$
and take $F =
\cup F_i$ the unique curve in the linear series $|k\cdot F|$ containing
them. Next, we have to single out one of these $k$ points, and label the
corresponding fiber
$F_i$. At this point $p_3$ and $p_4$ will be determined, as the unique
points of intersection of $F_i$ with the curves $C_3$ and $C_4$
respectively. Finally, we choose a curve $D \in |C+(n-k)F|$ passing through
the remaining $2n-k$ of the points $q_\a$ and having a point of
$(n-k)$-fold tangency with
$E$. (Note that the ordering of the $k$ points $q_\a$ chosen to lie on fibers
does not matter; all that counts is which one is chosen to lie on $F_i$.)

Now, the linear series $ |C+(n-k)F|$ cuts on the curve $E \cong \P^1$ the
complete linear series
$|\O_{\P^1}(n-k)|$. This linear series is parametrized by the space
$\P^{n-k}$ of polynomials of degree $n-k$ on $\P^1$ modulo scalars; and in
that projective space the divisors consisting simply of $n-k$ times a single
point---that is,
$(n-k)$th powers of linear forms---form a rational normal curve. In the
linear series $ |C+(n-k)F|$, then, the locus of curves $D$ having a single
point of $(n-k)$-fold intersection with $E$ is a cone over a rational normal
curve in $\P^{n-k}$ (with vertex the subseries $E + |(2n-k)F| \subset
|C+(n-k)F|$ of curves containing $E$); in particular, it has degree $n-k$.
There are thus exactly $n-k$ curves $D$ linearly equivalent to $|C+(n-k)F|$
passing through $p_1$, $p_2$ and the remaining
$2n-k$ of the points
$q_\a$ and having a point of
$(n-k)$-fold tangency with
$E$. In sum, the number of
fibers $X$ of this type in our family is ${2n \choose k} \cdot k \cdot
(n-k)$.

It remains to determine, for each such fiber  of our family, the
coefficient with which the corresponding point
$b\in B$ appears in the pullback divisor $\phi^*(0)$. To do this, we need to
know the local geometry of the family near $b \in B$; in particular, we
need to have the picture of the corresponding  fibers of the families $\X
\to B$ and $Y \to B$. For the first, the only thing we need to know is which
of the singular points of the fiber $X$ are limits of nodes of nearby
fibers and (in the case of the point of intersection of $D$ with $E$) how
many. The answer, as provided in Proposition
\ref{singtotalspace}, is that the points of intersection of $D$ with the fibers
$F_i$ are all limits of nodes on nearby curves; and the remaining
$(n-k-1)$ nodes of the general fiber of the family tend to the point of
intersection of $D$ with $E$.  When we normalize the total space of the
family, then, the curves $D$ and $F_i$ are pulled apart; and the point of
intersection of
$D$ with $E$ becomes a node, so that the fiber of $\X \to B$ over $b$ is

\newpage

\vspace*{4in}

\hskip.8in \special {picture Fnnormal}

\

But as we also saw in Proposition \ref{singtotalspace}, $\X$ will not be
smooth: at the point lying over each point of intersection of $E$ with a fiber
$F_i$, $\X$ will have a singularity of type $A_{n-k-1}$. Resolving each of
these, we arrive at this picture of the fiber of $\Y \to B$ over $b$:

\newpage

\vspace*{3.2in}

\ni \special {picture Fnsmooth}

\

Finally, we can blow down the extraneous curves $F_j$ and $G_{j,*}$ for $j
\ne i$ to arrive at the picture of the fiber $Z$  of the family $Z \to B$ of
semistable 4-pointed curves with smooth total space:

\vspace*{2.1in}

\ni \special {picture Fnfinal1}

\

Inasmuch as there are $(n-k)$ rational curves in the chain connecting the
components $D$ and $F_i$ containing the pairs $\{p_1, p_2\}$ and $\{p_3,
p_4\}$, each such fiber represents a point of multiplicity $n-k+1$ in the
divisor
$\phi^*(0)$. In sum, then, the fibers of this type contribute a total of
$$\sum_{k=1}^{n-1}{2n
\choose k} \cdot k \cdot (n-k) \cdot (n-k+1)$$ to the degree of $\phi^*(0)$.

\

\ni $\bullet_2$  \thinspace $p_1, p_3 \in D$; $p_2, p_4 \in F_i$ or $p_2,
p_4 \in D$;
$p_1, p_3 \in F_i$.  \enspace  In the first of these cases we are simply
exchanging the locations of
$p_2$ and
$p_3$ to obtain a fiber $X$ corresponding to a point $b$ in the inverse
image
$\phi^*(\infty)$; this will affect the count of the number of such fibers,
but not
the final configuration on the semistable model with smooth total space, so
the multiplicity of each such point in the divisor $\phi^*(0)$ will be as in
the preceding case $n-k+1$.

The difference here is that, because the fixed point $p_2$ lies on one of the
fiber components, we can put the remaining fiber components through
only $k-1$ of the points $q_\a$; at the same time, we can put the curve
$D$ through $p_1$ and the remaining $2n-k+1$. To specify such a curve,
then, we must first choose a subset of
$k-1$ of the
$2n$ points
$q_\a$ and take $F =
\cup F_j$ the unique curve in the linear series $|k\cdot F|$ containing them
and
$p_2$; the component of $F$ containing $p_2$ we will call $F_i$. As in the
preceding case, there will be exactly $n-k$ curves in the linear series
$|C+(n-k)F|$ passing through the remaining $2n-k+1$ points $q_\a$ and the
point $p_1$ and having a point of intersection multiplicity $n-k$ with $E$;
the curve
$D$ can be any of these. At this point
 $p_4$ will be determined, as the unique point of intersection of $F_i$ with
the curve  $C_4$; while $p_3$ can be taken to be any of the

$$ (D \cdot C_3) = \bigl( (C+(n-k)F) \cdot C \bigr) = 2n-k
$$

\noindent points of intersection of $D$ with $C_3$. The number of fibers
$X$ of this type in our family is thus ${2n \choose k-1} \cdot (n-k) \cdot
(2n-k)$.

As we said, each such fiber $X$ of our family is a pole of order $n-k+1$
of the cross-ratio function $\phi$. Finally, since exchanging $p_1$ with
$p_4$ (as in the second case above) yields an identical result, the total
contribution of the fibers of these types to the poles of $\phi$ is
$$ 2\cdot \sum_{k=1}^{n-1} {2n \choose k-1} \cdot (n-k) \cdot (2n-k) \cdot
(n-k+1).
$$

\

\ni $\bullet_3$  \thinspace $p_1, p_2 \in D$; $p_3 \in F_i$ and $p_4 \in
F_j$, $i \ne j$.
\enspace  This case is very similar to the first; again, we have first to select
a subset of $k$ of the $2n$ points $q_\a$ and take $F = \cup F_i$ the
unique curve in the linear series $|k\cdot F|$ containing them. We then
have to single out two of these $k$ points, and label the corresponding
fibers $F_i$ and $F_j$, which in turn determines the points $p_3 = F_i \cap
C_3$ and $p_4 = F_j \cap C_4$.  Finally, we take as before $D$ to be any of
the $n-k$ curves in $|C+(n-k)F|$ passing through
$p_1$ and $p_2$ and the remaining
$2n-k$ of the points $q_\a$ and having a point of
$(n-k)$-fold tangency with
$E$. Thus, the number of fibers $X$ of this type in our family is ${2n
\choose k} \cdot k \cdot (k-1) \cdot (n-k)$.

At this point, we see another difference with the preceding case: here, to
arrive at the fiber of the family of semistable 4-pointed curves with
smooth total space we
 blow down the curves $F_m$ and $G_{m,*}$ for all $m$ including
$i$ and $j$, to arrive at the simpler fiber:

\vspace*{.9in}

\ni \hskip.5in \special {picture Fnfinal2}

\

Since this is already semistable, each such fiber represents a simple zero
of $\phi$.  In
sum, then, the fibers of this type contribute a total of
$$\sum_{k=1}^{n-1}{2n
\choose k} \cdot k \cdot (k-1) \cdot (n-k)$$ to the degree of $\phi^*(0)$.

\

\ni $\bullet_4$  \thinspace $p_1, p_3 \in D$; $p_2 \in F_i$ and $p_4 \in
F_j$, $i \ne j$; or $p_2, p_4 \in D$; $p_1 \in F_i$ and $p_3 \in F_j$, $i \ne
j$.
\enspace  This case bears the same relation to the preceding as the second
did to the first: we are simply exchanging $p_2$ and $p_3$ (or $p_1$ and
$p_4$), so that the fibers $X$ will correspond to poles
rather zeroes of $\phi$; the multiplicity will be 1 as in the last case,
but the number of
such fibers will be different.

Take the case $p_1, p_3 \in D$ first. Such a fiber may be specified by
choosing first a subset of
$k-1$ of the
$2n$ points
$q_\a$ and taking $F = \cup F_i$ the unique curve in the linear series
$|k\cdot F|$ containing them and $p_2$; the component containing $p_2$
will will label $F_i$. We then have to single out one of these
$k$ points, and label the corresponding fiber
$F_j$ , which in turn determines the point $p_4 = F_j \cap C_3$. As before
we take
$D$ to be any of the $n-k$ curves in
$|C+(n-k)F|$ passing through
$p_1$ and the remaining
$2n-k+1$ of the points $q_\a$ and having a point of
$(n-k)$-fold tangency with
$E$.; the point $p_3$ may be any of the $2n-k$ points of $D
\cap C_3$.   Thus, the number of fibers $X$ of this type in our family is
${2n
\choose k-1} \cdot (k-1) \cdot (n-k) \cdot (2n-k)$.

As we said, each such fiber represents a point of multiplicity $1$ in the
divisor
$\phi^*(\infty )$; and since the case $p_2, p_4 \in D$ contributes an equal
number, the fibers of this type contribute a total of
$$2 \cdot \sum_{k=1}^{n-1}{2n
\choose k-1}  \cdot (k-1) \cdot (n-k) \cdot (2n-k) $$ to the degree of
$\phi^*(\infty)$.

\

\ni $\bullet_5$  \thinspace $p_3, p_4 \in D$; $p_1 \in F_i$ and $p_2 \in
F_j$, $i \ne j$.
\enspace  This case is also a variant of case (3) above: here we are
exchanging both
$p_2$ for
$p_3$ and
$p_1$ for
$p_4$. The result is that the fibers $X$ will once more correspond to
zeroes of $\phi$, with multiplicity 1 as in the last two cases, but again
there will
be a different number of such fibers.

To evaluate this number, note that this time such a fiber may be specified
by choosing first a subset of
$k-2$ of the
$2n$ points
$q_\a$ and taking $F = \cup F_i$ the unique curve in the linear series
$|k\cdot F|$ containing them and both $p_1$ and $p_2$;
 $F_i$ will be the component of $F$ containing $p_1$ and $F_j$  the
component  containing $p_2$.  We choose
$D$  any of the
$n-k$ curves in $|C+(n-k)F|$ passing through  the remaining
$2n-k+2$ of the points $q_\a$ and having a point of
$(n-k)$-fold tangency with
$E$; and then $p_3$  and $p_4$ may be any of the $2n-k$ points of $D
\cap C_3$ and $D
\cap C_4$ respectively. Thus, the number of fibers $X$ of this type in
our family is ${2n
\choose k-2}  \cdot (n-k) \cdot (2n-k)^2$; and since  each corresponding $b
\in B$
is a simple zero of $\phi$, such fibers contribute a total of
$$ \sum_{k=1}^{n-1}{2n
\choose k-2}  \cdot (n-k) \cdot (2n-k)^2 $$ to the degree of $\phi^*(0)$.

\

\ni $\bullet_6$  \thinspace $p_1 \in D$, $p_2 \in F_i$ and $p_3, p_4 \in
F_j$, $i \ne j$; or $p_1 \in F_i$, $p_2 \in D$ and $p_3, p_4 \in F_j$, $i \ne
j$. \enspace  These again give zeroes of $\phi$: to arrive at the
semistable model, in
the end we will blow down $D$ and the chains  $F_m$ and $G_{m,*}$ for all $m \ne j$.
Consider first the case $p_1 \in D$, $p_2 \in F_i$.  To specify such a
fiber we have to
select a subset of $k-1$ of the points $q_\a$: $F = \cup F_i$ will then be
the unique curve in the linear series
$|k\cdot F|$ containing them and
$p_2$ (and $F_i$ the component of $F$ containing $p_2$). Then we have to
choose which of the remaining $k-1$ components of $F$ is to be $F_j$, and
take $p_3 = F_j
\cap C_3$ and $p_4 = F_j \cap C_4$. Finally, we choose $D$  any of the
$n-k$ curves in $|C+(n-k)F|$ passing through $p_1$ and the remaining
$2n-k+1$ of the points $q_\a$ and having a point of
$(n-k)$-fold tangency with
$E$; there are thus a total of ${2n \choose k-1} \cdot (k-1) \cdot (n-k)$
such fibers in our family.

As for multiplicity, as we said the semistable model with smooth total
space is obtained from $Y$ by blowing down $D$ and the chains  $F_m$
and $G_{m,*}$ for all $m \ne j$ to arrive at a fiber of the form

\vspace*{1.2in}

\ni \special {picture Fnfinal3}

\

Since there are $(n-k-1)$ rational curves in the chain connecting the
components $E$ and $F_j$ containing the pairs $\{p_1, p_2\}$ and $\{p_3,
p_4\}$, each such fiber represents a zero of multiplicity $n-k$ of the function
$\phi$. Finally, the case $p_1 \in F_i$, $p_2 \in D$ contributes an equal
number; so that the fibers of this type contribute a total of
$$2 \cdot \sum_{k=1}^{n-1}{2n
\choose k-1}  \cdot (k-1) \cdot (n-k)^2$$ to the degree of $\phi^*(0)$.

\

\ni $\bullet_7$  \thinspace $p_1 \in D$, $p_3 \in F_i$ and $p_2, p_4 \in
F_j$, $i \ne j$; or $p_2 \in D$, $p_4 \in F_i$,  and $p_1, p_3 \in F_j$, $i \ne
j$. \enspace  These are obtained by exchanging  $p_2$ and $p_3$ in the
first case immediately above or
$p_1$ and
$p_4$ in the second (the remaining two possible switches will be
considered next), so that the fibers
$X$ will correspond to poles rather than zeroes of $\phi$ and will have the
same multiplicity $n-k$.  We thus simply have to count the number of
such fibers, which is straightforward: for example, in the first case ($p_1
\in D$, $p_3
\in F_i$), to specify such a fiber we have to select first a subset of
$k-1$ of the points $q_\a$ and take $F = \cup F_i$  the unique curve in the
linear series
$|k\cdot F|$ containing them and
$p_2$;  $F_j$ will be the component of $F$ containing $p_2$. Then we have
to choose which of the remaining $k-1$ components of $F$ is to be $F_i$,
and take
$p_3 = F_i
\cap C_3$ and $p_4 = F_j \cap C_4$. Finally, we choose $D$  any of the
$n-k$ curves in $|C+(n-k)F|$ passing through $p_1$ and the remaining
$2n-k+1$ of the points $q_\a$ and having a point of
$(n-k)$-fold tangency with
$E$; there are thus a total of ${2n \choose k-1} \cdot (k-1) \cdot (n-k)$
such fibers in our family, so that the fibers of this type contribute a total of
$$2 \cdot \sum_{k=1}^{n-1}{2n
\choose k-1}  \cdot (k-1) \cdot (n-k)^2$$ to the degree of $\phi^*(\infty)$.

\

\ni $\bullet_8$  \thinspace $p_3 \in D$, $p_1 \in F_i$ and $p_2, p_4 \in
F_j$, $i \ne j$; or $p_4 \in D$, $p_2 \in F_i$,  and $p_1, p_3 \in F_j$, $i \ne
j$. \enspace  These are the remaining two cases obtained by switching
points in case (6) above; as opposed to the immediately preceding case
these are obtained by exchanging
$p_1$ and
$p_4$ in the first case of (6) or
$p_2$ and
$p_3$ in the second. Thus the fibers
 will again correspond to poles rather than zeroes of $\phi$ and will again
appear with multiplicity
$n-k$ in $\phi^*(\infty)$, but the number of such fibers will be different.
To compute
it, consider the first case ($p_1 \in D$,
$p_3
\in F_i$). To specify such a fiber we have to select first a subset of
$k-2$ of the points $q_\a$ and take $F = \cup F_i$  the unique curve in the
linear series
$|k\cdot F|$ containing them and both $p_1$ and
$p_2$;  $F_i$ will be the component of $F$ containing $p_1$ and $F_j$  the
component  containing $p_2$.  This then fixes the point $p_4 = F_j \cap
C_4$. We choose
$D$  any of the
$n-k$ curves in $|C+(n-k)F|$ passing through  the remaining
$2n-k+2$ of the points $q_\a$ and having a point of
$(n-k)$-fold tangency with
$E$; and $p_3$ can be any of the $2n-k$ points of intersection of $D$ with
$C_3$.  There are thus a total of
${2n
\choose k-2}  \cdot (n-k) \cdot (2n-k)$ fibers of this type in our family; so
that the fibers of this type contribute a total of
$$2 \cdot \sum_{k=1}^{n-1}{2n
\choose k-2}  \cdot (2n-k) \cdot (n-k)^2$$ to the degree of $\phi^*(\infty)$.

\

We come now to the last three cases, those in which none of the four points
$p_i$ lie on $D$. The next one is the last to contribute to the degree of
$\phi^*(0)$.

\

\ni $\bullet_9$  \thinspace $p_1 \in F_i$, $p_2 \in F_j$ and $p_3, p_4 \in
F_\ell$, $i \ne j \ne \ell \ne i$. \enspace  To determine a fiber of this type
we have to specify first $k-2$ of the points $q_\a$, and let $F \in |kF|$
contain
$p_1$, $p_2$ and these $k-2$. We then have to pick a component of $F$
among those not containing $p_1$ or $p_2$, and call it $F_\ell$; $p_3$ and
$p_4$ will then be the points of intersection of $F_\ell$ with $C_3$ and
$C_4$ respectively. As always,
$D
\in |C+(n-k)F|$ can be any of the $n-k$ curves passing through the
remaining
$2n-k+2$ of the points $q_\a$ and having a point of $(n-k)$-fold tangency
with
$E$; so there are a total of ${2n
\choose k-2}  \cdot (k-2) \cdot (n-k) $ such fibers in our family. Finally,
for each such fiber, after blowing down  $D$ and the chains $F_m$ and
$G_{m,*}$ for all $m \ne \ell$ we arrive at the smooth semistable model,
whose special fiber has the form

\vspace*{2.4in}

\ni \hskip-.5in \special {picture Fnfinal4}

There being $n-k-1$ intermediate curves in this chain, each such fiber
corresponds to a zero of order
$n-k$  of $\phi$; the total contribution of such fibers is
thus

$$ \sum_{k=1}^{n-1}{2n
\choose k-2}  \cdot (k-2) \cdot (n-k)^2$$

\

\ni $\bullet_{10}$  \thinspace $p_1 \in F_i$  $p_3 \in F_j$ and $p_2, p_4
\in F_\ell$; or $p_2 \in F_i$  $p_4 \in F_j$ and $p_1, p_3 \in F_\ell$, $i
\ne j \ne \ell \ne i$.
\enspace  These are the two cases obtained from the preceding by
exchanging either $p_2$ and $p_3$ or $p_1$ and $p_4$; each such fiber
thus represents a pole of multiplicity $n-k$ of $\phi$. To
count the number, consider first the case $p_1 \in F_i$  $p_3 \in F_j$ and
$p_2, p_4 \in F_\ell$. To determine such a fiber we specify
$k-2$ of the points $q_\a$, and let $F \in |kF|$ contain
$p_1$, $p_2$ and these $k-2$; the component of $F$ containing $p_1$ we
will call
$F_i$, and then $p_3$ will be determined as the point $F_i \cap C_3$ . We
then have to pick a component of
$F$ among those not containing $p_1$ or $p_2$, and call it $F_\ell$;
$p_4$ will  be the point of intersection of $F_\ell$ with $C_4$. As always,
$D
\in |C+(n-k)F|$ can be any of the $n-k$ curves passing through the
remaining
$2n-k+2$ of the points $q_\a$ and having a point of $(n-k)$-fold tangency
with
$E$; so there are a total of ${2n
\choose k-2}  \cdot (k-2) \cdot (n-k) $ such fibers in our family. Counting
both possible exchanges, we see that the total contribution of such fibers to
the degree of $\phi^*(\infty)$ is

$$ 2 \cdot \sum_{k=1}^{n-1}{2n
\choose k-2}  \cdot (k-2) \cdot (n-k)^2$$

\

\ni $\bullet_{11}$  \thinspace $p_1, p_3 \in F_i$ and $p_2, p_4 \in F_j$,
$i
\ne j$ \enspace  This is our final case; note that there is no analogous
source of
zeroes of
$\phi$, since $p_1$ and $p_2$ do not lie on the same fiber of $\F_n$. First, to
count the number: such fibers are determined first by choosing
$k-2$ of the points $q_\a$, and letting $F$ be the union of fibers
containing them and $p_1$ and $p_2$; $F_i$ and $F_j$ will be the fibers
containing $p_1$ and
$p_2$ and $p_3$ and $p_4$ the points of intersection of these fibers with
$C_3$ and $C_4$ respectively, so nothing more need by specified. We only
have to choose $D$ among the $n-k$ curves in $|C+(n-k)F|$ passing through
the other $2n-k+2$ of the points $q_\a$ and meeting
$E$ in only one point, so that there are a total of just ${2n
\choose k-2}   \cdot (n-k) $ such fibers in our family.

To arrive at the smooth semistable model near such a fiber, we have to
blow down the curve $D$ and the chains $G_{\ell,*}$ for all $\ell \ne i, j$;
we arrive at the fiber

\vspace*{2.4in}

\ni \hskip-.5in \special {picture Fnfinal5}

Since there are $2n-2k-1$ rational curves in the chain connecting $F_i$ and
$F_j$, each such fiber gives a pole of multiplicity $2n-2k$; the total
contribution of
such fibers to the degree of
$\phi^*(\infty)$ is thus
$$ \sum_{k=1}^{n-1}{2n
\choose k-2}   \cdot (n-k) \cdot (2n-2k)$$

We are now ready to add up all the contributions to the degrees of
$\phi^*(0)$ and $\phi^*(\infty)$, equating the results and solving for
$N(2C)$.  We have
\begin{equation*}
\begin{split}
\deg(\phi^*(0)) = n \cdot N(2C) + n^3{2n \choose n-1} \\
 + \sum_{k=1}^{n-1} (n-k)
\biggl[ &{2n \choose k} \bigl( k(n-k+1) + k(k-1) \bigr) \\ +
 &{2n \choose k-1} \bigl( 2(k-1)(n-k) \bigr) \\ +
 &{2n \choose k-2} \bigl( (2n-k)^2 + (k-2)(n-k) \bigr) \biggr] .
\end{split}
\end{equation*} while on the other hand
\begin{equation*}
\begin{split}
\deg(\phi^*(\infty)) =  n^3{2n \choose n}  \\ +
\sum_{k=1}^{n-1} (n-k)
\biggl[ &{2n \choose k-1} \bigl( 2(2n-k)(n-k+1) + 2(2n-k)(k-1) +
2(k-1)(n-k)
\bigr) \\ +
 &{2n \choose k-2} \bigl( 2(2n-k)(n-k) +2(k-2)(n-k) + 2(n-k) \bigr)  \biggr].
\end{split}
\end{equation*} Combining these, we arrive at the expression
\begin{equation*} n \cdot N(2C) = n^3 \bigl({2n \choose n} - {2n \choose
n-1} \bigr) + S
\end{equation*} where
\begin{equation*}
\begin{split}
\!\!\!\!\!\! S =  \sum_{k=1}^{n-1} (n-k)
\biggl[ &{2n \choose k} \bigl( -k(n-k+1) - k(k-1) \bigr) \\ +
 &{2n \choose k-1} \bigl( 2(2n-k)(n-k+1) + 2(2n-k)(k-1) +
2(k-1)(n-k)-2(k-1)(n-k)
\bigr) \\ +
 &{2n \choose k-2} \bigl( 2(2n-k)(n-k) +2(k-2)(n-k) + 2(n-k) -(2n-k)^2 -
(k-2)(n-k)
\bigr) \biggr] .
\end{split}
\end{equation*} The expression for $S$ may be reduced immediately to
\begin{equation*}
\begin{split} S =  \sum_{k=1}^{n-1} (n-k)
\biggl[ &{2n \choose k} \cdot (-kn) \\ +
 &{2n \choose k-1}  \cdot 2(2n-k)n \\ +
 &{2n \choose k-2}  \cdot (-kn) \biggr] .
\end{split}
\end{equation*} (Note that we can now enlarge the formal limits of
summation to include $k=0$ without affecting the sum; this will be
convenient in the following calculations.) To reduce this further, we
separate it into two terms: we write
$S = S' - S''$, where
$$ S' = \sum_{k=0}^{n-1} 4n^2(n-k) {2n \choose k-1}
$$ and
$$ S'' =  \sum_{k=0}^{n-1} (n-k) \cdot kn \cdot
\biggl[ {2n \choose k}  +
 2{2n \choose k-1}  + {2n \choose k-2} \biggr] .
$$ The expression for $S''$ telescopes nicely: since
$$ {2n \choose k-1}  +
 {2n \choose k-2}  = {2n+1 \choose k-1},
$$
$$ {2n \choose k}  +
 {2n \choose k-1}  = {2n+1 \choose k}
$$ and
$$ {2n+1 \choose k}  +
 {2n+1 \choose k-1}  = {2n+2 \choose k},
$$ we have simply
$$ S'' =  \sum_{k=0}^{n-1} (n-k) \cdot kn \cdot {2n+2 \choose k}
$$ As for $S'$, we can combine that with the remaining two terms in the
expression for $N(2C)$, and together they simplify. To start with,
observe that
\begin{equation*}
\begin{split} {2n \choose n} - {2n \choose n-1} &= {(2n)! \over n!n!} -
{(2n)! \over (n-1)!(n+1)!}\\ &= {(2n)! \over n!(n+1)!} \cdot \bigl( (n+1) - n
\bigr) \\ &= {(2n)! \over n!(n+1)!} \\ &= {1 \over n}{2n \choose n-1}
\end{split}
\end{equation*} so
$$ n^3 \bigl({2n \choose n} - {2n \choose n-1} \bigr) = n^2{2n \choose n-1}
$$ Now, combining this with the expression for $S'$ above, we have
\begin{equation*}
\begin{split} n^2{2n \choose n-1} &+ \sum_{k=0}^{n-1} 4n^2(n-k) {2n
\choose k-1} \\ &= n^2 \biggl( {2n \choose n-1} + 4{2n \choose n-2} + 8{2n
\choose n-3} + \dots + (4n-4){2n \choose 0} \biggr)
\end{split}
\end{equation*} To reduce this, we use the relation
$$ {2n \choose n-1} + {2n \choose n-2} = {2n+1 \choose n-1}
$$ to absorb the first term; then 3 times the relation
$$ {2n \choose n-2} + {2n \choose n-3} = {2n+1 \choose n-2}
$$ to absorb the rest of the second; then 5 times the relation
$$ {2n \choose n-3} + {2n \choose n-4} = {2n+1 \choose n-3}
$$ to absorb the rest of the third, and so on; ultimately, we arrive at
\begin{equation*}
\begin{split} n^2{2n \choose n-1} &+ \sum_{k=0}^{n-1} 4n^2(n-k) {2n
\choose k-1} \\ &= n^2 \biggl( {2n+1 \choose n-1} + 3{2n+1 \choose n-2} +
5{2n+1 \choose n-3} +
\dots + (2n-1){2n+1 \choose 0} \biggr)
\end{split}
\end{equation*} Now we play the same game again: using the relation
$$ {2n+1 \choose n-1} + {2n+1 \choose n-2} = {2n+2 \choose n-1}
$$ to absorb the first term, then twice the relation
$$ {2n+1 \choose n-2} + {2n+1 \choose n-3} = {2n+2 \choose n-2}
$$ to absorb the remainder of the second, and so on, we may re-express
this as
\begin{equation*}
\begin{split} n^2{2n \choose n-1} &+ \sum_{k=0}^{n-1} 4n^2(n-k) {2n
\choose k-1} \\ &= n^2 \biggl( {2n+2 \choose n-1} + 2{2n+1 \choose n-2} +
3{2n+1 \choose n-3} +
\dots + n{2n+1 \choose 0} \biggr) \\ &= n^2 \sum_{k=0}^{n-1} (n-k)  {2n+2
\choose k}
\end{split}
\end{equation*} Finally, we can combine this and the expression above for
$S''$: we have
\begin{equation*}
\begin{split} n \cdot N(2C) &= n^3 \bigl( {2n \choose n} - {2n \choose
n-1} \bigr) + S' - S'' \\  &= n^2 \sum_{k=0}^{n-1} (n-k)  {2n+2 \choose k} - n
\sum_{k=0}^{n-1} k(n-k) {2n+2 \choose k} \\ &= n \sum_{k=0}^{n-1}
(n-k)^2  {2n+2 \choose k}.
\end{split}
\end{equation*}
We have therefore proved the
\begin{thm}
\label{2C}
Let $N(2C)$ be the number of irreducible rational curves in the linear
series $|2C|$ on $\F _n$ passing through
 $2n+3$ points, then
$$ N(2C) = \sum_{k=0}^{n-1} (n-k)^2  {2n+2 \choose k} .
$$
\end{thm}

For example, on $\F_2$ we have
$$ N(2C) = {6 \choose 1} + 4{6 \choose 0} = 6+4 =10 ;
$$ on
$\F_3$ we have
$$N(2C) = {8 \choose 2} + 4{8 \choose 1}  + 9{8 \choose 0} = 28 + 32 + 9 =
69 ;
$$ and on
$\F_4$ we have
$$N(2C) = {10 \choose 3} + 4{10 \choose 2}  + 9{10 \choose 1} + 16{10
\choose 0}= 120 + 180 + 90 + 16 = 406
$$ and so on.

\

We will now show how to arrive at an expression of $N(2C)$ on
$\F_n$ as a coefficient of a simple generating function. We simply write out
the sum involved, and then telescope it using the standard binomial
relations as before: that is, we write
$$ N(2C)=  {2n+2 \choose n-1} + 4{2n+2 \choose n-2} + 9{2n+2 \choose
n-3} +
\dots + n^2{2n+2 \choose 0}
$$ and use the relations ${2n+2 \choose n-1} + {2n+2 \choose n-2} = {2n+3
\choose n-1}$, $3{2n+2 \choose n-2} + 3{2n+2 \choose n-3} = 3{2n+3
\choose n-2}$, and so on to rewrite this as
\begin{equation*}
\begin{split} N(2C) &=  {2n+3 \choose n-1} + 3{2n+3 \choose n-2} +
6{2n+3 \choose n-3} +
\dots + {n(n+1) \over 2}{2n+3 \choose 0} \\ &= \sum_{k=0}^{n-1} {n-k+1
\choose 2}{2n+3 \choose k}.
\end{split}
\end{equation*} We can also think of this as the coefficient of $t^n$ in the
product of the power series
$$
\sum {2n+3 \choose k} t^k = (1+t)^{2n+3}
$$ and
$$
\sum {\ell+2 \choose 2} t^\ell = {1 \over (1-t)^3}
$$ so that we can write $N(2C)$ on
$\F_n$ as the coefficient
$$ N(2C)  =  \bigg[ {(1+t)^{2n+3} \over (1-t)^3} \biggr]_{t^n} .
$$

\

\

\subsection{A formula for $\F _n$}
\label{Fn}
We conclude our paper with a formula for the general ruled surface $\F _n$.
 Here we  define
 the function $\g_{i_1,\dots,i_t}(D_{i_1},\dots,D_{i_t})$  giving the
contribution to the cross-ratio
corresponding to  a given decomposition
$D = D_1+D_2+\dots+D_t+E$  or $D=D_1+D_2$, with $D_j\in \VDj$.
Recall that the variety $\VDi$ is the closure in $|D|$
of the locus of irreducible rational curves that have a point of  contact
of order $i$ with the
exceptional curve $E$ (cf. \ref{term} ).
We define

$$
\g_{i_1,\dots,i_t}(D_{i_1},\dots,D_{i_t}):= \prod (i_j \NDj ) \cdot
$$
$$
\cdot \Bigl[ {\rD -3 \choose \rDu -1,\rDt -1,
\rDT,....}  [\sum _{j\geq 3} {(C\cdot D_j) \over i_j}\bigl( {(C\cdot D_1)
\over i_1} +  {(C\cdot
D_2)
\over i_2}\bigr) - \sum _{j\geq 3} {(C\cdot D_j)^2 \over i_j}]
$$
$$ - {\rD -3 \choose \rDu -2,\rDt ,
\rDT,....}  [\sum _{j\geq 2} (C\cdot D_j)^2 ({1 \over i_j }+ {1 \over i_1
}) + {1 \over i_1
}\sum _ {2\leq j<k \leq t}(C\cdot D_j)(C\cdot D_k) ] \Bigr]
$$

\ni In these terms, we can state
\begin{thm}
\label{Fn} Let $D$ be a divisor on the surface $\F_n$.  Let $\ND$ be the
number of
irreducible rational curves in $|D|$ that pass through
$r_0(D)$ general points of $\F_n$. Then
$$ n\ND =
$$
$$\sum _{D_1 + D_2 = D} \Prod \g _{1,1}(D_1,D_2) +
$$
$$
 + \sum _{t=2}^{n} \; \sum _{D_1+D_2+\dots+D_t=D-E} \; \sum
_{i_1,\dots,i_t}  \prod _{j
:i_j=1}(E\cdot D_{i_j})
\g_{i_1,\dots,i_t}(D_{i_1},\dots,D_{i_t})
$$
\end{thm}

\

\

\ni {\bf{Acknowledgments}}. \thinspace Our interest in these questions grew
out of conversations with Enrico
Arbarello, Ciro Ciliberto and Bill Fulton, to whom we are very grateful.

\

\end{document}